\documentclass[10pt]{article}
\pdfoutput=1

\usepackage{amsmath}
\usepackage{mathtools}
\usepackage{stackrel}
\usepackage{scalerel}
\usepackage{leftidx}
\usepackage{xfrac}

\usepackage[charter]{mathdesign}
\usepackage{euscript}
\let\mathbb\undefined
\usepackage{bbold}
\DeclareSymbolFont{usualmathcal}{OMS}{cmsy}{m}{n}
\DeclareSymbolFontAlphabet{\mathcal}{usualmathcal}

\usepackage{authblk}

\usepackage[shortlabels]{enumitem}

\usepackage{hyperref}
\usepackage{cleveref}
\hypersetup{colorlinks=true, linkcolor=red!50!black, citecolor=green!50!black, urlcolor=blue!80!black}

\usepackage{xcolor}
\usepackage{tcolorbox}
\usepackage{graphicx}
\usepackage{subfigure}

\usepackage[all,2cell,arrow,matrix]{xy} \UseAllTwocells \SilentMatrices
\usepackage{tikz-cd}
\usepackage{tikz}
\usetikzlibrary{arrows,arrows.meta,%
decorations.markings,decorations.pathreplacing,%
decorations.pathmorphing,calc}

\usepackage{verbatim}
\usepackage{comment}

\usepackage[nottoc]{tocbibind} % References in ToC
\usepackage{appendix}

\usepackage{geometry}
\tikzset{->-/.style={decoration={markings,mark=at position #1 with {\arrow{Stealth}}},postaction={decorate}},->-/.default=0.55}
\usepackage[amsmath,amsthm,thmmarks]{ntheorem}
\theoremstyle{definition}

\newtheorem{thm}{Theorem}[section]

\newtheorem{prop}[thm]{Proposition}

\newtheorem{cor}[thm]{Corollary}
\newtheorem{lem}[thm]{Lemma}

{
\theoremsymbol{\mbox{${\blacksquare}$}}
\newtheorem{defn}[thm]{Definition}
}
{
\theoremsymbol{\mbox{$\heartsuit$}}
\newtheorem{expl}[thm]{Example}
}
{
\theoremsymbol{\mbox{$\diamondsuit$}}
\newtheorem{rem}[thm]{Remark}
}
\qedsymbol{\mbox{$\square$}}

\numberwithin{equation}{section}
\numberwithin{thm}{section}

\newcommand\be            {\begin{equation}}
\newcommand\ee            {\end{equation}}
\newcommand\bea           {\begin{eqnarray}}
\newcommand\eea         {\end{eqnarray}}
\newcommand\bnu          {\begin{enumerate}}
\newcommand\enu          {\end{enumerate}}
\newcommand\bit          {\begin{itemize}}
\newcommand\eit          {\end{itemize}}

\newcommand{\pf}{\begin{proof}}
\newcommand{\epf}{\qed\end{proof}}

\makeatletter
\providecommand{\leftsquigarrow}{%
  \mathrel{\mathpalette\reflect@squig\relax}%
}
\newcommand{\reflect@squig}[2]{%
  \reflectbox{$\m@th#1\rightsquigarrow$}%
}
\makeatother

\newcommand\Zb			{\mathbb{Z}}

\newcommand\bk			{\mathbb{k}}

\newcommand\CA			{\EuScript{A}}
\newcommand\CB			{\EuScript{B}}
\newcommand\CC			{\EuScript{C}}
\newcommand\CD			{\EuScript{D}}
\newcommand\CE			{\EuScript{E}}

\newcommand\CG			{\EuScript{G}}

\newcommand\CM			{\EuScript{M}}

\newcommand\CP			{\EuScript{P}}

\newcommand\CV			{\EuScript{V}}
\newcommand\CW			{\EuScript{W}}
\newcommand\CX			{\EuScript{X}}

\newcommand{\FZ}			{\text{\usefont{U}{euf}{m}{n}Z}}

\newcommand\SC			{\mathsf{C}}
\newcommand\SD			{\mathsf{D}}

\DeclareMathOperator{\Hom}{Hom}
\DeclareMathOperator{\End}{End}
\DeclareMathOperator{\Coend}{Coend}
\DeclareMathOperator{\Aut}{Aut}

\DeclareMathOperator{\id}{id}

\DeclareMathOperator{\tr}{tr}

\DeclareMathOperator{\ev}{ev}
\DeclareMathOperator{\coev}{coev}

\DeclareMathOperator{\nat}{Nat}
\DeclareMathOperator{\fun}{Fun}

\DeclareMathOperator{\Mod}{Mod}
\DeclareMathOperator{\LMod}{LMod}
\DeclareMathOperator{\RMod}{RMod}
\DeclareMathOperator{\BMod}{BMod}

\DeclareMathOperator{\map}{Map}
\DeclareMathOperator{\Kar}{Kar}

\newcommand{\op}			{\mathrm{op}}
\newcommand{\rev}			{\mathrm{rev}}

\newcommand{\one}			{\mathbb{1}}

\newcommand\forget  {\text{\usefont{U}{euf}{m}{n}f}}

\newcommand\set			{\mathrm{Set}}

\newcommand\vect			{\mathrm{Vec}}

\newcommand\rep			{\mathrm{Rep}}

\newcommand\Irr			{\mathrm{Irr}}

\newcommand{\bscale}	{0.7}
\makeatletter
\newcommand{\ec}[2][]	{{\@ec{#1 |}{#2}}}
\newcommand{\bc}[2][]	{{\@ec{#1}{#2}}}
\newcommand{\@ec}[2]	{\mathchoice
  {\displaystyle \raise.9ex\hbox{$\scaleobj{\bscale}{#1}$} {#2}}%
  {\textstyle \raise.9ex\hbox{$\scaleobj{\bscale}{#1}$} {#2}}%
  {\scriptstyle \raise.55ex\hbox{$\scriptstyle \scaleobj{\bscale}{#1}$} {#2}}%
  {\scriptscriptstyle \raise.38ex\hbox{$\scriptscriptstyle \scaleobj{\bscale}{#1}$} {#2}}%
}
\makeatother

\newcommand{\dquotient} {/\mkern-6mu/}

\begin{document}

\title{The 2-character theory of finite 2-groups}

\author[a,b]{Mo Huang \thanks{Email: \href{mailto:jasmine.huang.527@gmail.com}{\tt jasmine.huang.527@gmail.com}}}
\author[c,d]{Hao Xu \thanks{Email: \href{mailto:haoxu@imada.sdu.dk}{\tt haoxu@imada.sdu.dk}}}
\author[e,f,g,h]{Zhi-Hao Zhang \thanks{Email: \href{mailto:zhangzhihao@bimsa.cn}{\tt zhangzhihao@bimsa.cn}}}
\affil[a]{School of Mathematical Sciences, East China Normal University, Shanghai, 200241, China}
\affil[b]{Department of Physics, The University of Hong Kong, Pokfulam Road, Hong Kong, China}
\affil[c]{Centre for Quantum Mathematics, University of Southern Denmark, Odense, 5230, Denmark}
\affil[d]{Mathematisches Institut, Georg-August-Universität Göttingen, 37073, Germany}
\affil[e]{Beijing Institute of Mathematical Sciences and Applications, Beijing, 101408, China}
\affil[f]{Wu Wen-Tsun Key Laboratory of Mathematics of Chinese Academy of Sciences, \authorcr
School of Mathematical Sciences, University of Science and Technology of China, Hefei, 230026, China}
\affil[g]{International Quantum Academy, Shenzhen, 518048, China}
\affil[h]{Guangdong Provincial Key Laboratory of Quantum Science and Engineering, \authorcr
Southern University of Science and Technology, Shenzhen, 518055, China}
\author{}
\date{\vspace{-5ex}}

\maketitle

\begin{abstract}
We generalize the notion of character for 2-representations of finite 2-groups. The properties of 2-characters bear strong similarities to those classical characters of finite groups, including additivity, multiplicativity, conjugation invariance, $S^1$-invariance and orthogonality. With a careful analysis using homotopy fixed points and quotients for categories with 2-group actions, we prove that the category of class functors on a 2-group $\CG$ is equivalent to the Drinfeld center of the 2-group algebra $\vect_\CG$, which categorifies the Fourier transform on finite abelian groups. After transferring the canonical nondegenerate braided monoidal structure from $\FZ_1(\vect_\CG)$, we discover that irreducible 2-characters of $\CG$ coincide with full centers of the corresponding 2-representations, which are in a one-to-one correspondence with Lagrangian algebras in the category of class functors on $\CG$. In particular, the fusion rule of $2\rep(\CG)$ can be calculated from the pointwise product of Lagrangian algebras as class functors. From a topological quantum field theory (TQFT) point of view, the commutative Frobenius algebra structure on a 2-character is induced from a 2D topological sigma-model with target space $\lvert \mathrm{B} \CG \rvert$.%, while the nondegeneracy is achieved only if the 2-character is irreducible. 
% Moreover, this HQFT yields the notion of joint 2-character on torus, and we find a natural interpretation of the modularity of joint 2-characters in terms of the action of the modular group $\mathrm{SL}(2,\mathbb{Z})$ on the underlying torus.
\end{abstract}

\tableofcontents

\section{Introduction}
Character is a basic notion and plays an important role in representation theory of finite groups, over an algebraically closed field $\Bbbk$ of characteristic $0$. For a finite group $G$ and a finite-dimensional representation $(V,\rho) \in \rep(G)$, its character is a function $\chi_V \colon G \to \bk$ defined by $\chi_V(g) \coloneqq \tr(\rho(g))$. Let us list some basic results in the character theory:
\bnu
\item The characters have conjugation invariance. So they are class functions, i.e., functions that are constant on conjugacy classes.
\item The representations are determined by their characters.
\item For $V,W \in \rep(G)$, we have $\chi_{V \oplus W} = \chi_V + \chi_W$ and $\chi_{V \otimes W} = \chi_V \cdot \chi_W$. It follows that the characters generate a ring which is equivalent to the Grothendieck ring of $\rep(G)$ with coefficient extended to $\Bbbk$. Therefore, the characters can be viewed as de-categorifications of representations and they are more computable, hence provide much information about representations.
\item There is an inner product on the space of class functions. For $V,W \in \rep(G)$, the inner product of their characters $\langle \chi_V,\chi_W \rangle$ is equal to the dimension of the hom space $\Hom_G(V,W)$. In particular, when $V,W$ are irreducible representations, Schur's lemma implies that the inner product $\langle \chi_V,\chi_W \rangle$ is either $1$ or $0$, depending on whether $V$ and $W$ are isomorphic. This is called the orthogonality of characters.
\enu

It is natural to categorify the character theory and the above results to finite 2-groups. To define the 2-character $\chi_\CV$ of a finite semisimple 2-representation $\CV$ of a finite 2-group $\CG$, we want to use an analogous definition and thus need the notion of the trace of a linear functor. There are at least two ways to define this trace. The first way is based on the following observation: for a finite-dimensional vector space $V$ and a linear map $f \colon V \to V$, the following composite map $\bk \to \bk$, viewed as a number, is equal to $\tr(f)$:
\[
\bk \xrightarrow{\coev} V^* \otimes V \xrightarrow{1 \otimes f} V^* \otimes V \xrightarrow{\ev} \bk ,
\]
where $\coev(1) = \sum_{i=1}^n v^i \otimes v_i$ for a basis $\{v_i\}$ of $V$ and its dual basis $\{v^i\}$, and $\ev(\alpha \otimes v) = \alpha(v)$. The coevaluation map $\coev$ and the evaluation map $\ev$ exhibit $V^*$ as a dual of $V$, and their existences are due to the finite-dimensionality of $V$. Similarly, for a finite semisimple category $\CV$ (which is equivalent to a Kapranov-Voevodsky 2-vector space \cite{KV94}, i.e., a finite direct sum of $\vect$), it is also dualizable and the dual is the opposite category $\CV^\op$ obtained from $\CV$ by reversing all morphisms. Thus for a linear functor $F \colon \CV \to \CV$ we can similarly define its trace, which is a linear functor from $\vect$ to $\vect$ and can be identified as a vector space. Using the graph calculus, this trace is usually depicted as a circle with a defect labeled by $F$. Another way is proposed by Ganter and Kapranov in any 2-category \cite{GK08}, which defines the trace of $F$ to be the hom space $\Hom_{\Hom_{2\vect}(\CV,\CV)}(1_\CV,F)$, that is, the space $\nat(1_\CV,F)$ of natural transformations from the identity functor $1_\CV$ to $F$. In our case, these two definitions are equivalent. Parallelly, Bartlett also developed the 2-character theory for finite groups \cite{Bar09, Bar09a} via extended topological quantum field theory (TQFT). Subsequent works, such as \cite{Oso10,RW18}, build upon Ganter-Kapranov by providing more explicit computations. Related categorical and generalized character theories appeared in several other settings. 
In derived algebraic geometry, Ben-Zvi, Francis, and Nadler \cite{BZFN10} related the Drinfeld center to loop spaces, providing a geometric analogue of the categorical trace. In homotopy theory, the Dennis trace from algebraic $K$-theory to topological Hochschild homology \cite{BHM93,DGM12} is a categorified Chern character, and the generalized character theory of Hopkins, Kuhn, and Ravenel \cite{HKR00} gives class functions on commuting $n$-tuples; at height $2$ this is a direct precursor of the ``categorical character'' of Ganter and Kapranov \cite{GK08}.

The 2-character $\chi_\CV$ and its properties can also be understood from a TQFT point of view. However, in this paper we do not aim to construct the relevant TQFT; all of our discussion of TQFT is intended only as motivation. The finite semisimple category $\CV \in 2\vect$ is not only dualizable, but also fully dualizable in the sense that the evaluation and coevaluation 1-morphisms are also dualizable in $2\vect$. Hence by the cobordism hypothesis \cite{BD95,Lur08}, $\CV$ defines a fully extended 2D TQFT. A fully extended 2D TQFT is a symmetric monoidal 2-functor $Z \colon \mathrm{Cob}_2^{\mathrm{fr}} \to 2\vect$, where $\mathrm{Cob}_2^{\mathrm{fr}}$ is the 2-category of 2-1-0 framed cobordisms with the tensor product being the disjoint union. The cobordism hypothesis states that for any fully dualizable object $\CV$ there is an essentially unique fully extended 2D TQFT $Z_\CV$ which maps the single point to $\CV$. Moreover, if we consider $\CV$ as a finite semisimple 2-representation of $\CG$, the TQFT $Z_\CV$ can be extended to a $\CG$-TQFT, that is, a symmetric monoidal 2-functor $Z_\CV \colon \mathrm{Cob}_2^\CG \to 2\vect$ where the cobordisms in $\mathrm{Cob}_2^\CG$ are all equipped with a continuous map to the classifying space $\lvert \mathrm B \CG \rvert$. Then the value of $Z_\CV$ on the circle $S^1$ with chosen maps to $\lvert \mathrm B \CG \rvert$ recovers the 2-character of $\CV$ (see Remark \ref{rem_G_TQFT}). By considering the value of $Z_\CV$ on certain cylinders, we obtain the conjugation invariance of the 2-character (see Remark \ref{rem_TQFT_conjugation_invariance}), which categorifies the classical conjugation invariance of the characters for finite groups. We can also consider the value of $Z_\CV$ on certain torus and obtain a numerical invariant of the 2-representation $\CV$ with an obvious modular invariance (see Section \ref{sec_numerical_2-character}), which recovers a classical notion of character \cite{GK08,Dav10a} when $\CG$ is a 1-group.

There is a ubiquitous phenomenon in 2D TQFT or CFT called the open-closed duality \cite{MS06,RFFS07,KR09}. It suggests that the closed string algebra should be the full center of the open string algebra, where the closed string algebra and the open string algebra are the values (state spaces) of the field theory on the closed string (circle) $S^1$ and the open string (interval) $[0,1]$, respectively. For the TQFT $Z_\CV$ discussed above, the closed string algebra is the 2-character $\chi_\CV$. So we expect that $\chi_\CV$ is a full center, which should be a Lagrangian algebra in $\FZ_1(\vect_\CG)$ by \cite{KR09}. We prove this result in Theorem \ref{thm_2-character_Lagrangian_algebra}, and it also implies that the 2-representations are determined by their 2-characters.

The last result of this work is to establish the orthogonality of 2-characters. We define an inner product for functors $\CG \to \vect$ equipped with a conjugation invariance data, which is a vector space. Then we show that the inner product of two 2-characters $\chi_\CV,\chi_\CW$ is isomorphic to the `dimension' of the hom category $\fun_\CG(\CV,\CW)$, where the dimension of a finite semisimple category $\CP$ is defined to be the endomorphism algebra $\End(1_\CP)$ of the identity functor. This gives a direct categorification of the orthogonality of characters for finite groups.

\medskip
The layout of this paper is as follows. In Section \ref{sec_preliminary}, we review some notions in tensor category theory and fix the conventions. We also review the basic results about 2-representations of finite 2-groups, which were mainly developed by two of the authors in \cite{HZ23}. In Section \ref{sec_homotopy_fixed_point_quotient}, we discuss the homotopy fixed points and homotopy quotients of categories equipped with 2-group actions. These notions are useful when we study the conjugation invariance. The homotopy fixed points (equivariantization) of categories equipped with 1-group actions and the homotopy quotients (action groupoids) of sets equipped with 1-group actions are well-known. Our constructions are direct generalizations. Section \ref{sec_2-character} is the main part of this work. For a finite 2-group $\CG$, we define the 2-character of a finite semisimple 2-representation of $\CG$ and discuss its conjugation invariance data, based on some intuitions from TQFT. We also show that the irreducible 2-characters are Lagrangian algebras in a certain category. In particular, the equivalence classes of 2-representations are one-to-one corresponding to the isomorphism classes of 2-characters. Finally, we establish the orthogonality of 2-characters, which categorifies the classical result in finite group representation theory.

\medskip
\noindent \textbf{Conventions and Notations}: In this work, the ground field $\bk$ is an algebraically closed field of characteristic zero. For a category $\CC$, its opposite category $\CC^\op$ is obtained from $\CC$ by reversing all morphisms. For a monoidal category $\CD$, its reversed category $\CC^\rev$ is the same underlying category $\CD$ equipped with the reversed tensor product $\otimes^\rev$ defined by $x \otimes^\rev y \coloneqq y \otimes x$. If $\CV$ is a finite semisimple category, the set of isomorphism classes of simple objects in $\CV$ is denoted by $\Irr(\CV)$. In a category $\CC$, the hom space between two objects $x,y \in \CC$ is denoted by $\CC(x,y)$. For a monoidal category $\CD$ and a left $\CD$-module category $\CM$, the module action functor is usually denoted by $\odot \colon \CD \times \CM \to \CM$, and for any $x \in \CM$, the right adjoint of $- \odot x \colon \CD \to \CM$ (if exists) is denoted by $[x,-] \colon \CM \to \CD$, called the internal hom.

\medskip
\noindent \textbf{Acknowledgements}: We would like to thank Ansi Bai, Hank Chen, Liang Kong and Tian Lan for helpful discussions, and also thank Ansi Bai and Tian Lan for finding some typos. Part of this work was conducted during MH and ZHZ's visit at BIMSA. We would like to thank Hao Zheng for his support throughout this period. HX would also like to thank the Max Planck Institute for Mathematics and Peter Teichner for their hospitality during his stay in Bonn for the \emph{Twinned workshop on ``Quantum Field Theory and Topological Phases via Homotopy Theory and Operator Algebras''}, during which he had many fruitful discussions with João Faria Martins, who pointed out several important references on representation theory of 2-groups. MH is supported
by Research Grants Council of Hong Kong under GRF 17311322. HX is supported by DAAD Graduate School Scholarship Programme (57572629) and DFG Project 398436923. ZHZ is supported by the start-up grant of BIMSA and the China Postdoctoral Science Foundation under Grant Number 2025M783072.

\section{Preliminaries} \label{sec_preliminary}

\subsection{Ends and coends}

In this subsection we briefly recall some basic facts about (co)ends. For more details, we refer readers to \cite{Lor21}.

\begin{defn}
Let $\CC,\CD$ be categories and $P,Q \colon \CC^\op \times \CC \to \CD$ be functors. A \emph{dinatural transformation} $\beta \colon P \xRightarrow{..} Q$ is a family $\{\beta_x \colon P(x,x) \to Q(x,x)\}_{x \in \CC}$ of morphisms in $\CD$, such that the following diagram commutes for every morphism $f \colon x \to y$ in $\CC$:
\[
\xymatrix{
P(y,x) \ar[r]^{P(f,1)} \ar[d]_{P(1,f)} & P(x,x) \ar[r]^{\beta_x} & Q(x,x) \ar[d]^{Q(1,f)} \\
P(y,y) \ar[r]^{\beta_y} & Q(y,y) \ar[r]^{Q(f,1)} & Q(x,y)
}
\]
The \emph{end} of $P \colon \CC^\op \times \CC \to \CD$ is an object $\End(P) \in \CD$ equipped with a dinatural transformation $\tau \colon \Delta_{\End(P)} \xRightarrow{..} P$ (here $\Delta$ denotes the constant functor) satisfying the following universal property: for any object $d \in \CD$ equipped with a dinatural transformation $\rho \colon \Delta_d \xRightarrow{..} P$, there exists a unique morphism $\bar \rho \colon d \to \End(P)$ such that $\tau_x \circ \bar \rho = \rho_x$ for all $x \in \CC$.

Dually, the \emph{coend} of $P$ is an object $\Coend(P) \in \CD$ equipped with a dinatural transformation $\tau \colon P \xRightarrow{..} \Delta_{\Coend(P)}$ satisfying the following universal property: for any object $d \in \CD$ equipped with a dinatural transformation $\rho \colon P \xRightarrow{..} \Delta_d$, there exists a unique morphism $\bar \rho \colon \Coend(P) \to d$ such that $\bar \rho \circ \tau_x = \rho_x$ for all $x \in \CC$.
\end{defn}

\begin{rem} \label{rem_end_limit}
Define the category $\mathrm{TA}(\CC)$ of \emph{twisted arrows} in $\CC$ as follows:
\bit
\item The objects are morphisms in $\CC$.
\item The morphisms from $f \colon x \to y$ to $g \colon z \to w$ are pairs of morphisms $(u \colon z \to x,v \colon y \to w)$ satisfying $v \circ f \circ u = g$.
\item The composition and identity morphisms are induced by those of $\CC$.
\eit
Then from a functor $P \colon \CC^\op \times \CC \to \CD$ we can define a functor $\bar P \colon \mathrm{TA}(\CC) \to \CD$ by
\[
\begin{array}{c}
\xymatrix{
(f \colon x \to y) \ar[d]^{(u,v)} \\
(g \colon z \to w)
}
\end{array}
\mapsto
\begin{array}{c}
\xymatrix{
P(x,y) \ar[d]^{P(u,v)} \\
P(z,w)
}
\end{array}
\]
Moreover, a dinatural transformation $\tau \colon \Delta_d \xRightarrow{..} P$ defines a natural transformation $\bar \tau \colon \Delta_d \Rightarrow \bar P$ by $\bar \tau_f \coloneqq P(1_x,f) \circ \tau_x = P(f,1_y) \circ \tau_y$. Conversely, a natural transformation $\bar \tau \colon \Delta_d \Rightarrow \bar P$ defines a dinatural transformation $\tau \colon \Delta_d \xRightarrow{..} P$ by $\tau_x \coloneqq \bar \tau_{1_x}$. These two constructions are mutually inverse. Hence we have a canonical isomorphism $\End(P) \simeq \lim \bar P$. Similarly, $\Coend(P)$ can also be written as a colimit.

More generally, (co)ends can be naturally written as weighted (co)limits. See Remark \ref{rem_end_as_weighted_limit}.
\end{rem}

We usually denote $\End(P)$ and $\Coend(P)$ by the following integral notation:
\[
\End(P) \eqqcolon \int_{x \in \CC} P(x,x) , \quad \Coend(P) \eqqcolon \int^{x \in \CC} P(x,x) .
\]
Here $x$ is only a `dummy index' and can be replaced by any other symbol. This integral notation is convenient for stating the `Fubini's theorem' for (co)ends. Suppose $\CC,\CD,\CE$ are categories and $Q \colon \CC^\op \times \CD^\op \times \CC \times \CD \to \CE$ is a functor. If the end $\int_{y \in \CD} Q(x,y,x,y) \in \CE$ exists for any $x \in \CC$, we have a functor
\[
\int_{y \in \CD} Q(-,y,-,y) \colon \CC^\op \times \CC \to \CE .
\]
Similarly, we have a functor
\[
\int_{x \in \CC} Q(x,-,x,-) \colon \CD^\op \times \CD \to \CE
\]
if these ends exist. Then one can easily show that there are canonical isomorphisms between the following three ends:
\[
\int_{x \in \CC} \int_{y \in \CD} Q(x,y,x,y) , \quad \int_{y \in \CD} \int_{x \in \CC} Q(x,y,x,y) , \quad \int_{(x,y) \in \CC \times \CD} Q(x,y,x,y)
\]
because they satisfy the same universal property. A similar statement for coends also holds. Here and in the following, when we refer to the `canonical isomorphism' between two (co)ends or (co)limits without assuming their existence, we implicitly mean that if either exists, then the other also exists and they are canonically isomorphic.

\begin{expl}
Let $\CC,\CD$ be categories and $F,G \colon \CC \to \CD$ be functors. Then there is a canonical isomorphism
\[
\nat(F,G) \simeq \int_{x \in \CC} \CD(F(x),G(x)) .
\]
Indeed, the morphism $\nat(F,G) \to \CD(F(x),G(x))$ sending a natural transformation $\alpha \colon F \Rightarrow G$ to its component $\alpha_x \colon F(x) \to G(x)$ is dinatural in $x$. Then one can easily verify that $\nat(F,G)$ satisfies the universal property of ends.
\end{expl}

\begin{expl}
Let $\CC$ be a $\bk$-linear category and $F \colon \CC \to \vect$ be a $\bk$-linear functor. Then for every $x \in \CC$ there is a canonical isomorphism
\[
\int^{y \in \CC} \CC(y,x) \otimes F(y) \simeq F(x) .
\]
Indeed, the morphism $\CC(y,x) \otimes F(y) \to F(x)$ sending $f \otimes a$ to $F(f)(a)$ is dinatural in $y$. Then one can easily verify that $F(x)$ satisfies the universal property of coends. Another proof is to use the Yoneda lemma:
\begin{multline*}
\vect \biggl(\int^{y \in \CC} \CC(y,x) \otimes F(y),V \biggr) \simeq \int_{y \in \CC} \vect(\CC(y,x) \otimes F(y),V) \\
\simeq \int_{y \in \CC} \vect(\CC(y,x),\vect(F(y),V)) \\
\simeq \nat(\CC(-,x),\vect(F(-),V)) \simeq \vect(F(x),V) .
\end{multline*}
Dually, we also have
\[
\int_{y \in \CC} \CC(x,y)^* \otimes F(y) \simeq \int_{y \in \CC} \vect(\CC(x,y),F(y)) \simeq \nat(\CC(x,-),F) \simeq F(x) .
\]
These isomorphisms are also called the \emph{co-Yoneda lemma} or the \emph{density formula} or the \emph{ninja Yoneda lemma} in \cite{Lor21}.
\end{expl}

\begin{lem} \label{lem_end_Karoubi_complete}
Let $\CC,\CD$ be $\bk$-linear categories. Suppose $\iota \colon \CC \hookrightarrow \Kar(\CC)$ is a Karoubi completion and $P \colon \Kar(\CC)^\op \times \Kar(\CC) \to \CD$ is a functor. Then there is a canonical isomorphism $\End(P) \simeq \End(P \circ (\iota \times \iota))$, or written in the integral notation:
\[
\int_{x \in \Kar(\CC)} P(x,x) \simeq \int_{x \in \CC} P(\iota(x),\iota(x)) .
\]
Similarly, there is a canonical isomorphism $\Coend(P) \simeq \Coend(P \circ (\iota \times \iota))$.
\end{lem}

\pf
We prove the canonical isomorphism between two ends. The proof for the canonical isomorphism of two coends is similar.

By abuse of notation, we omit the embedding $\iota \colon \CC \to \Kar(\CC)$ and view $\CC$ as a full subcategory of $\Kar(\CC)$. For any $d \in \CD$, we have
\be \label{eq_end_as_weighted_limit}
\CD \biggl(d,\int_{x \in \CC} P(x,x) \biggr) \simeq \int_{(y,x) \in \CC^\op \times \CC} \vect(\CC(y,x),\CD(d,P(y,x)) \simeq \nat(\Hom_\CC,\CD(d,P(-,-))) ,
\ee
where the right hand side is a hom space in the functor category $\fun(\CC^\op \times \CC,\vect)$. Similarly we have
\[
\CD \biggl(d,\int_{x \in \Kar(\CC)} P(x,x) \biggr) \simeq \nat(\Hom_{\Kar(\CC)},\CD(d,P(-,-))) ,
\]
where the right hand side is a hom space in $\fun(\Kar(\CC)^\op \times \Kar(\CC),\vect) \simeq \fun(\Kar(\CC^\op \times \CC),\vect)$. The universal property of Karoubi completion implies that these two functor categories are equivalent via the restriction
\[
F \mapsto F \circ (\iota \times \iota) \colon \fun(\Kar(\CC)^\op \times \Kar(\CC),\vect) \to \fun(\CC^\op \times \CC,\vect) .
\]
Hence we obtain canonical isomorphisms
\[
\CD \biggl(d,\int_{x \in \CC} P(x,x) \biggr) \simeq \CD \biggl(d,\int_{x \in \Kar(\CC)} P(x,x) \biggr)
\]
for all $d \in \CD$, which induce a canonical isomorphism $\int_{x \in \CC} P(x,x) \simeq \int_{x \in \Kar(\CC)} P(x,x)$.
\epf

\begin{rem} \label{rem_end_as_weighted_limit}
The isomorphism \eqref{eq_end_as_weighted_limit} means that an end is a weighted limit:
\[
\int_{x \in \CC} P(x,x) \simeq {\lim}^{\Hom_\CC} P .
\]
Then the above proof is equivalent to the fact that every (weighted) cone can be uniquely extended along the retraction of idempotents. Similarly, coends can be written as weighted colimits.
\end{rem}

%\pf
%We prove the canonical isomorphism between two ends. The proof for the canonical isomorphism of two coends is similar.
%
%By abuse of notation, we omit the embedding $\iota \colon \CC \to \Kar(\CC)$ and view $\CC$ as a full subcategory of $\Kar(\CC)$. By Remark \ref{rem_end_limit} we have
%\[
%\int_{x \in \Kar(\CC)} P(x,x) \simeq \lim_{f \in \mathrm{TA}(\Kar(\CC))} \bar P(f) , \quad \int_{x \in \CC} P(x,x) \simeq \lim_{f \in \mathrm{TA}(\CC)} \bar P(f) .
%\]
%From the construction of Karoubi completion it is easy to see that $\mathrm{TA}(\Kar(\CC)) = \Kar(\mathrm{TA}(\CC))$. Then by the universal property of Karoubi completion, any natural transformation $\Delta_d \Rightarrow \bar P \colon \mathrm{TA}(\CC) \to \CD$ can be uniquely extends to a natural transformation $\Delta_d \Rightarrow \bar P \colon \mathrm{TA}(\Kar(\CC)) \to \CD$. Thus any dinatural transformation $\Delta_d \xRightarrow{..} P \colon \CC^\op \times \CC \to \CD$ can be uniquely extends to a dinatural transformation $\Delta_d \xRightarrow{..} P \colon \Kar(\CC)^\op \times \Kar(\CC) \to \CD$. Hence we see that $\int_{x \in \Kar(\CC)} P(x,x)$ and $\int_{x \in \CC} P(x,x)$ satisfy the same universal property, thus are canonically isomorphic.
%\epf

\begin{cor} \label{cor_semisimple_end_direct_sum}
Let $\CC,\CD$ be $\bk$-linear categories. Assume that $\CC$ is finite semisimple and $\CC(x,x) \simeq \bk$ for every simple object $x \in \CC$. Then for any $\bk$-bilinear functor $P \colon \CC^\op \times \CC \to \CD$ we have canonical isomorphisms
\[
\int_{x \in \CC} P(x,x) \simeq \bigoplus_{x \in \Irr(\CC)} P(x,x) \simeq \int^{x \in \CC} P(x,x) .
\]
\end{cor}

\pf
Let $\CC_0 \subseteq \CC$ be a skeleton of the full subcategory consists of all simple objects in $\CC$. Then $\Kar(\CC_0) \simeq \CC$. By Lemma \ref{lem_end_Karoubi_complete} we have $\int_{x \in \CC} P(x,x) \simeq \int_{x \in \CC_0} P(x,x)$. Since $P$ is $\bk$-bilinear and $\CC(x,x) \simeq \bk$ for all simple objects $x$, the dinaturality of a family of morphisms $\{\tau_x \colon d \to P(x,x)\}_{x \in \CC_0}$ is trivial. So we have $\int_{x \in \CC} P(x,x) \simeq \bigoplus_{x \in \Irr(\CC)} P(x,x)$. The proof for the second isomorphism is similar.
\epf

Since we are working over an algebraically closed field $\bk$, the condition that $\CC(x,x) \simeq \bk$ for all simple object $x$ in Corollary \ref{cor_semisimple_end_direct_sum} always holds.

%\begin{rem}
%The condition that $\Hom_\CC(x,x) \simeq \bk$ for all simple object $x$ in Proposition \ref{cor_semisimple_end_direct_sum} is necessary. For example, let $\vect_\Cb$ and $\vect_\Rb$ be the category of finite-dimensional $\Cb$-vector spaces and the category of finite-dimensional $\Rb$-vector spaces, repectively, and let $F \colon \vect_\Cb \to \vect_\Rb$ be the forgetful functor. Then we have
%\[
%\int_{V \in \vect_\Cb} \Hom_\Rb(V,V) = Z(\Hom_\Rb(\Cb,\Cb)) .
%\]
%\end{rem}

\begin{rem}
A generalization of Lemma \ref{lem_end_Karoubi_complete} and Corollary \ref{cor_semisimple_end_direct_sum} for non-semisimple categories can be found in \cite[Section 5.1.3]{KL01}.
\end{rem}

The following useful proposition is also called the `change of variables' theorem for (co)ends.

\begin{prop} \label{prop_change_variable_end}
Let $\CC,\CD,\CE$ be categories and $F \colon \CC \to \CD , G \colon \CD \to \CC , P \colon \CC^\op \times \CD \to \CE$ be functors. If $F$ is left adjoint to $G$, then there are canonical isomorphisms
\[
\int_{x \in \CC} P(x,F(x)) \simeq \int_{y \in \CD} P(G(y),y) , \quad \int^{x \in \CC} P(F(x),x) \simeq \int^{y \in \CD} P(y,G(y)) .
\]
\end{prop}

\pf
Let us prove the first isomorphism. The proof of the second one is similar.

To construct a morphism $f \colon \int_{x \in \CC} P(x,F(x)) \to \int_{y \in \CD} P(G(y),y)$, by the universal property of ends, it suffices to give a family of morphisms
\[
\biggl\{f_y \colon \int_{x \in \CC} P(x,F(x)) \to P(G(y),y) \biggr\}_{y \in \CD}
\]
that is dinatural in $y$. We define
\[
f_y \coloneqq \biggl( \int_{x \in \CC} P(x,F(x)) \xrightarrow{\tau_{G(y)}} P(G(y),FG(y)) \xrightarrow{P(1,\varepsilon_y)} P(G(y),y) \biggr) ,
\]
where $\varepsilon$ is the counit for the adjunction $F \dashv G$. In other words, $f$ is the unique morphism rendering the following diagram commutative for all $w \in \CD$:
\[
\xymatrix@C=5em{
\int_{x \in \CC} P(x,F(x)) \ar@{-->}[r]^-{f} \ar[d]_{\tau_{G(w)}} & \int_{y \in \CD} P(G(y),y) \ar[d]^{\tau_w} \\
P(G(w),FG(w)) \ar[r]^-{P(1,\varepsilon_w)} & P(G(w),w)
}
\]
Similarly, we define a morphism $g \colon \int_{y \in \CD} P(G(y),y) \to \int_{x \in \CC} P(x,F(x))$ to be the unique morphism rendering the following diagram commutative for all $z \in \CC$:
\[
\xymatrix@C=5em{
\int_{y \in \CD} P(G(y),y) \ar@{-->}[r]^-{g} \ar[d]_{\tau_{F(z)}} & \int_{x \in \CC} P(x,F(x)) \ar[d]^{\tau_a} \\
P(GF(z),F(z)) \ar[r]^-{P(\eta_z,1)} & P(z,F(z))
}
\]
where $\eta$ is the unit for the adjunction $F \dashv G$. Then it is not hard to verify that $f$ and $g$ are mutually inverse to each other by the universal property of ends.

There is another simpler proof:
\begin{multline*}
\int_{x \in \CC} P(x,F(x)) \simeq \int_{(x,y) \in \CC \times \CD} \set(\CD(F(x),y),P(x,y)) \\
\simeq \int_{(x,y) \in \CC \times \CD} \set(\CC(x,G(y)),P(x,y)) \simeq \int_{y \in \CD} P(G(y),y) ,
\end{multline*}
whre the first and last isomorphisms are due to the Yoneda lemma.
\epf

\subsection{Drinfeld centers}

In this subsection we briefly recall some basic facts about Drinfeld center construction and fix the convention.

%\begin{defn} \label{defn:Drinfeld_center}
Let $\CA$ be a monoidal category. The \emph{Drinfeld center} or \emph{monoidal center} of $\CA$ is the braided monoidal category $\FZ_1(\CA)$ defined as follows:
\bit
\item The objects are pairs $(a,\gamma_{-,a})$, where $a \in \CA$ and $\gamma_{-,a} \colon {- \otimes a} \Rightarrow {a \otimes -} \colon \CA \to \CA$ is a natural isomorphism (called a half-braiding) such that the following diagram commutes for all $b,c \in \CA$,
\be \label{diag_Drinfeld_center_object}
\begin{array}{c}
\xymatrix{
 & b \otimes (a \otimes c) \ar[r]^{\alpha_{b,a,c}^{-1}} & (b \otimes a) \otimes c \ar[dr]^{\gamma_{b,a} \otimes \id_c} \\
b \otimes (c \otimes a) \ar[ur]^{\id_b \otimes \gamma_{c,a}} \ar[dr]^{\alpha_{b,c,a}^{-1}} & & & (a \otimes b) \otimes c \\
 & (b \otimes c) \otimes a \ar[r]^{\gamma_{b \otimes c,a}} & a \otimes (b \otimes c) \ar[ur]^{\alpha_{a,b,c}^{-1}}
}
\end{array}
\ee
where $\alpha$ is the associator of $\CA$.
\item The morphisms from $(a,\gamma_{-,a})$ to $(b,\gamma_{-,b})$ are morphisms $f \in \CA(a,b)$ such that the following diagram commutes for all $c \in \CA$.
\be \label{diag_Drinfeld_center_morphism}
\begin{array}{c}
\xymatrix{
c \otimes a \ar[r]^{\gamma_{c,a}} \ar[d]_{\id_c \otimes f} & a \otimes c \ar[d]^{f \otimes \id_c} \\
c \otimes b \ar[r]^{\gamma_{c,b}} & b \otimes c
}
\end{array}
\ee
\item The tensor product of two objects $(a,\gamma_{-,a})$ and $(b,\gamma_{-,b})$ is $(a \otimes b,\gamma_{-,a \otimes b})$, where $a \otimes b$ is the tensor product in $\CA$ and the half-braiding $\gamma_{c,a \otimes b} \colon c \otimes a \otimes b \to a \otimes b \otimes c$ is the composition
\[
c \otimes a \otimes b \xrightarrow{\gamma_{c,a} \otimes \id_b} a \otimes c \otimes b \xrightarrow{\id_a \otimes \gamma_{c,b}} a \otimes b \otimes c .
\]
Here we omit the associators of $\CA$.
\item The tensor unit is $(\one,\gamma_{-,\one})$, where $\one \in \CA$ is the tensor unit of $\CA$ and the half-braiding $\gamma_{a,\one} \colon a \otimes \one \to \one \otimes a$ is
\[
a \otimes \one \xrightarrow{r_a} a \xrightarrow{l_a^{-1}} \one \otimes a ,
\]
or simply $\gamma_{a,\one} = \id_a$ if we ignore the unitors.
\item The associator and left/right unitor of $\FZ_1(\CA)$ are induced by those of $\CA$.
\item The braiding of $(a,\gamma_{-,a})$ and $(b,\gamma_{-,b})$ is given by $\gamma_{a,b} \colon a \otimes b \to b \otimes a$. Note that it is independent of the half-braiding of $a$.
\eit
%\end{defn}

\begin{thm}[\cite{Mueg03a,DGNO10}]
Let $\CA$ be a fusion category. Then $\FZ_1(\CA)$ is a nondegenerate braided fusion category.
\end{thm}

Suppose $\CA$ is rigid. A half-braiding $\{\gamma_{b,a} \colon b \otimes a \to a \otimes b\}_{b \in \CA}$ can be equivalently defined by the mates
\[
\{\gamma'_{b,a} \colon (b \otimes a) \otimes b^R \to a\} .
\]
Then the commutative diagram \eqref{diag_Drinfeld_center_object} is equivalent to the following diagram:
\be \label{diag_Drinfeld_center_object_coend}
\begin{array}{c}
\xymatrix{
((b \otimes c) \otimes a) \otimes (b \otimes c)^R \ar[r]^-{\simeq} \ar[d]_{\gamma'_{b \otimes c,a}} & (b \otimes ((c \otimes a) \otimes c^R)) \otimes b^R \ar[d]^{(1 \otimes \gamma'_{c,a}) \otimes 1} \\
a & (b \otimes a) \otimes b^R \ar[l]^{\gamma'_{b,a}}
}
\end{array}
\ee
The naturality of $\gamma_{-,a}$ is equivalent to the dinaturality of $\gamma'_{-,a}$, i.e., the following diagram commutes for all morphisms $f \colon b \to c$ in $\CA$:
\[
\xymatrix{
(b \otimes a) \otimes c^R \ar[r]^{1 \otimes f^R} \ar[d]_{(f \otimes 1) \otimes 1} & (b \otimes a) \otimes b^R \ar[d]^{\gamma'_{b,a}} \\
(c \otimes a) \otimes c^R \ar[r]^-{\gamma'_{c,a}} & a
}
\]
The commutative diagram \eqref{diag_Drinfeld_center_morphism} for a morphism $f \colon (a,\gamma_{-,a}) \to (b,\gamma_{-,b})$ is equivalent to the following diagram:
\be \label{diag_Drinfeld_center_morphism_coend}
\begin{array}{c}
\xymatrix{
(c \otimes a) \otimes c^R \ar[r]^-{\gamma'_{c,a}} \ar[d]_{(1 \otimes f) \otimes 1} & a \ar[d]^{f} \\
(c \otimes b) \otimes c^R \ar[r]^-{\gamma'_{c,b}} & b
}
\end{array}
\ee

If the coend $T(a) \coloneqq \int^{b \in \CA} (b \otimes a) \otimes b^R$ exists in $\CA$ (for example, if $\CA$ is a fusion category), we obtain a functor $T \colon \CA \to \CA$ and the dinatural transformation $\gamma'_{-,a}$ can be equivalently defined by a morphism $T(a) \to a$. Moreover, $\int^{b \in \CA} (b \otimes a) \otimes b^R$ also has a canonical half-braiding structure by the change of variables (see Proposition \ref{prop_change_variable_end}), which is the unique morphism rendering the following diagram commutative for all $d \in \CA$:
\[
\xymatrix@C=5em{
c \otimes \int^{b \in \CA} b \otimes a \otimes b^R \ar@{-->}[r] & \int^{b \in \CA} b \otimes a \otimes b^R \otimes c \\
c \otimes d \otimes a \otimes d^R \ar[r]^-{1 \otimes \coev_c} \ar[u]_{1 \otimes \tau_d} & c \otimes d \otimes a \otimes d^R \otimes c^R \otimes c \ar[u]^{\tau_{c \otimes d} \otimes 1}
}
\]
where $\coev_c \colon \one \to c^R \otimes c$ is the coevaluation morphism for the right dual of $c$. Thus we obtain a functor $\bar T \colon \CA \to \FZ_1(\CA)$ lifting $T$ along the forgetful functor $\FZ_1(\CA) \to \CA$. This also equips the functor $T$ a canonical monad structure. The commutative diagram \eqref{diag_Drinfeld_center_object_coend} is the associativity of a $T$-module structure on $a$, and the commutative diagram \eqref{diag_Drinfeld_center_morphism_coend} means that $f$ is a $T$-module morphism. In this case, $\FZ_1(\CA)$ is equivalent to the Eilenberg-Moore category of $T$ (i.e., the $T$-module category). It also follows that the left adjoint of the forgetful functor $\FZ_1(\CA) \to \CA$ is $\bar T \colon \CA \to \FZ_1(\CA)$.

\begin{rem}
Moreover, the monad $T$ has a structure of a quasi-triangular Hopf monad, and the equivalence between $\FZ_1(\CA)$ and the Eilenberg-Moore category of $T$ is a braided monoidal equivalence \cite{DS07,BV12}.
\end{rem}

Equivalently, a half-braiding $\{\gamma_{b,a} \colon b \otimes a \to a \otimes b\}_{b \in \CA}$ can also be defined by another mates
\[
\{\gamma''_{b,a} \colon a \to (b^R \otimes a) \otimes b\} .
\]
Also, if the ends $S(a) \coloneqq \int_{b \in \CA} (b^R \otimes a) \otimes b$ exists in $\CA$, the above morphisms are dinatural in $b$ and thus equivalent to a single morphism $a \to S(a)$. Similarly, $S(a)$ has a canonical half-braiding for every $a \in \CA$ and thus $S$ can be lifted to be a functor $\bar S \colon \CA \to \FZ_1(\CA)$, which is right adjoint to the forgetful functor $\FZ_1(\CA) \to \CA$. This means that the functor $S$ has a canonical comonad structure, and $\FZ_1(\CA)$ is equivalent to the Eilenberg-Moore category of $S$ (i.e., the $S$-comodule category) if these ends exist.

\subsection{Full centers}

%Let. If $C \in \FZ_1(\CC)$ is a commutative algebra and $A \in \CC$ is an algebra, then $C \otimes A$ is also an algebra in $\CC$ with the multiplication defined by
%\[
%C \otimes A \otimes C \odot A \xrightarrow{1 \otimes \gamma_{A,C} \otimes 1} C \otimes C \otimes A \otimes A \xrightarrow{\mu_C \otimes \mu_A} C \otimes A ,
%\]
%where $\gamma_{-,C}$ is the half-braiding of $C$, and $\mu_C$ and $\mu_A$ are the multiplications of $C$ and $A$, respectively.

\begin{defn}[\cite{Dav10}]
Let $\CC$ be a monoidal category and $A$ be an algebra in $\CC$. For an object $X \in \FZ_1(\CC)$, a morphism $f \colon X \to A$ is called \emph{central} if the following diagram commutes:
\[
\xymatrix{
A \otimes X \ar[rr]^-{\gamma_{A,X}} \ar[d]_{1 \otimes f} & & X \otimes A \ar[d]^{f \otimes 1} \\
A \otimes A \ar[r]^-{\mu} & A & A \otimes A \ar[l]_-{\mu}
}
\]
where $\gamma_{-,X}$ is the half-braiding of $X$ and $\mu$ is the multiplication of $A$. A \emph{full center} of $A$ is an object $Z(A) \in \FZ_1(\CC)$ equipped with a central morphism $z \colon Z(A) \to A$ satisfying the following universal property: for every $X \in \FZ_1(\CC)$ and central morphism $f \colon X \to A$, there exists a unique morphism $\bar f \colon X \to Z(A)$ in $\FZ_1(\CC)$ such that $z \circ \bar f = f$.
\end{defn}

We list some basic facts about full centers.

\begin{thm}[\cite{Dav10}]
Let $\CC$ be a monoidal category and $A \in \CC$ be an algebra. The full center of $A$, if exists, is a commutative algebra in $\FZ_1(\CC)$. Moreover, the full center is a Morita invariant, i.e., only depends on the left $\CC$-module category $\RMod_A(\CC)$ of right $A$-modules in $\CC$.
\end{thm}

Davydov proved this theorem by defining the full center $Z(\CM)$ for left $\CC$-modules $\CM$ and showing that $Z(\RMod_A(\CC)) \simeq Z(A)$.

\begin{expl} \label{expl_full_center_end}
Let us recall an explicit construction of a full center in \cite{DKR15}. Let $\CC$ be a rigid monoidal category and $A \in \CC$ be an algebra. Then the full center $Z(A)$ (i.e., the full center $Z(\RMod_A(\CC))$) is the end
\[
\int_{M \in \RMod_A(\CC)} [M,M] \in \CC
\]
equipped with the a canonical half-braiding induced by the change of variables (see Proposition \ref{prop_change_variable_end}), which is the unique morphism rendering the following diagram commutative for all $P \in \RMod_A(\CC)$:
\[
\xymatrix@C=5em{
x \otimes \int_{M \in \RMod_A(\CC)} [M,M] \ar@{-->}[r] \ar[d]_{1 \otimes \tau_{x^R \odot P}} & \int_{M \in \RMod_A(\CC)} [M,M] \otimes x \ar[d]^{\tau_P \otimes 1} \\
x \otimes [x^R \odot P,x^R \odot P] \ar[d]_{\simeq} & [P,P] \otimes x \\
[x^R \odot P, (x \otimes x^R) \odot P] \ar[r]^-{[1,\ev_x \odot 1]} & [x^R \odot P,P] \ar[u]_{\simeq}
}
\]
Equivalently, this half-braiding can also be defined by its mate:
\be \label{diag_full_center_half_braiding}
\begin{array}{c}
\xymatrix@C=5em{
\int_{M \in \RMod_A(\CC)} [M,M] \ar@{-->}[r] \ar[d]_{\tau_{x^R \odot P}} & x^R \otimes \int_{M \in \RMod_A(\CC)} [M,M] \otimes x \ar[d]^{1 \otimes \tau_P \otimes 1} \\
[x^R \odot P,x^R \odot P] \ar[r]^-{\simeq} & x^R \otimes [P,P] \otimes x
}
\end{array}
\ee
The multiplication of $Z(A)$ is the unique morphism rendering the following diagram commutative for all $P \in \RMod_A(\CC)$:
\be \label{diag_full_center_multiplication}
\begin{array}{c}
\xymatrix@C=5em{
\int_{M \in \RMod_A(\CC)} [M,M] \otimes \int_{N \in \RMod_A(\CC)} [N,N] \ar@{-->}[r]^-{\mu} \ar[d]_{\tau_P \otimes \tau_P} & \int_{P \in \RMod_A(\CC)} [P,P] \ar[d]^{\tau_P} \\
[P,P] \otimes [P,P] \ar[r] & [P,P]
}
\end{array}
\ee
where the bottom arrow is the composition morphism for internal homs. Here the internal hom can be written more explicitly \cite{Ost03}: $[M,N] = (M \otimes_A N^R)^L \in \CC$ for $M,N \in \RMod_A(\CC)$.
\end{expl}

Recall that an algebra $A$ in a fusion category $\CC$ is \emph{connected} if $\CC(\one,A) \simeq \bk$, and is \emph{separable} if there exists an $A$-$A$-bimodule map $\delta \colon A \to A \otimes A$ such that $\mu \circ \delta = 1_A$, where $\mu \colon A \otimes A \to A$ is the multiplication of $A$. If $\CB$ is a braided fusion category and $A \in \CB$ is a commutative algebra, a right $A$-module $(M,\rho \colon M \otimes A \to M)$ is called \emph{local} if $\rho \circ c_{A,M} \circ c_{M,A} = \rho$. A \emph{Lagrangian algebra} \cite{DMNO13} in $\CB$ is a commutative connected separable algebra $A \in \CB$ such that the category $\Mod_A^{\mathrm{loc}}(\CB)$ of local $A$-modules in $\CB$ is equivalent to $\vect$.

\begin{thm}[\cite{KR09,DMNO13}]
Let $\CC$ be a fusion category and $A \in \CC$ be a connected separable algebra. Then the full center $Z(A) \in \FZ_1(\CC)$ is a Lagrangian algebra. Moreover, the category $\RMod_{Z(A)}(\FZ_1(\CC))$ of right $Z(A)$-modules in $\FZ_1(\CC)$ is equivalent to the category $\BMod_{A|A}(\CC)$ of $A$-$A$-bimodules in $\CC$.
\end{thm}

\begin{rem} \label{rem_full_center_Lagrangian}
More generally, for a multi-fusion category $\CC$, it can be decomposed as the direct sum of indecomposable multi-fusion categories $\CC = \bigoplus_{i = 1}^n \CC_i$. Then $\FZ_1(\CC) = \bigoplus_{i=1}^n \FZ_1(\CC_i)$ is the direct sum of nondegenerate braided fusion categories. Suppose $A \in \CC$ is an indecomposable separable algebra. It must be contained in $\CC_{i_0}$ for some $i_0$. It is easy to see that the full center $Z(A)$ of $A$ in $\FZ_1(\CC)$ is the same as the full center of $A$ in $\FZ_1(\CC_{i_0})$. Therefore, $Z(A)$ is a Lagrangian algebra in $\FZ_1(\CC_{i_0})$.
\end{rem}

\subsection{Finite 2-groups and their 2-representations} \label{sec_2-group}

Let us recall some basic facts about 2-groups and 2-representations. See \cite{HZ23} for more details.

Given a monoidal category $\CC$, an object $x \in \CC$ is called \emph{invertible} if there exists $y \in \CC$ such that $x \otimes y \simeq \one \simeq y \otimes x$. This inverse object $y$ is both the left and right dual of $x$.

\begin{defn}
A \emph{2-group} is a monoidal category in which every object is invertible and every morphism is an isomorphism.
\end{defn}

Given a 2-group $\CG$, we always fix a left dual $g^L$ for every object $g \in \CG$. This left dual is automatically a right dual, so we also use $g^*$ to denote this dual object. The \emph{inverse functor} $(-)^* \colon \CG \to \CG^\rev$ is defined by
\[
(a \colon g \to h) \mapsto (a^* \coloneqq (a^L)^{-1} \colon g^* = g^L \to h^L = h^*) .
\]
In other words, let $I \colon \CG^\op \to \CG$ be the \emph{inverse morphism functor} that sends each morphism to its inverse. Then we have $(-)^* = I \circ (-)^L$.

\begin{defn}
Let $\CG$ be a 2-group. Its \emph{first homotopy group} $\pi_1(\CG)$ is the group of isomorphism classes of objects of $\CG$ with the multiplication induced by the tensor product of $\CG$. Its \emph{second homotopy group} $\pi_2(\CG)$ is the endomorphism group $\CG(\one,\one)$ of the tensor unit $\one$.
\end{defn}

By the Eckmann-Hilton argument, the second homotopy group $\pi_2(\CG)$ of any 2-group $\CG$ is always an abelian group. Also, $\pi_1(\CG)$ naturally acts on $\pi_2(\CG)$ by `conjugation' and the associator $\alpha$ is a 3-cocycle lies in $Z^3(\pi_1(\CG);\pi_2(\CG))$. It is known \cite{Sin75} that a 2-group $\CG$ is completely determined by two homotopy groups $\pi_1(\CG),\pi_2(\CG)$, together with the $\pi_1(\CG)$-action on $\pi_2(\CG)$ and the cohomology class $[\alpha] \in H^3(\pi_1(\CG);\pi_2(\CG))$ represented by the associator.

\begin{defn}
A 2-group $\CG$ is called \emph{finite} if both $\pi_1(\CG)$ and $\pi_2(\CG)$ are finite.
\end{defn}

\begin{defn}
Let $\CG$ be a finite 2-group. The 2-category of \emph{finite semisimple 2-representations of $\CG$} is the 2-functor 2-category $\fun(\mathrm B \CG,2\vect)$, where $\mathrm B \CG$ is the one-point delooping 2-category of $\CG$, and $2\vect$ is the 2-category of finite semisimple categories. We also denote this 2-category $\fun(\mathrm B \CG,2\vect)$ by $2\rep(\CG)$ for simplicity.
\end{defn}

In other words, a finite semisimple 2-representations of $\CG$ is a finite semisimple left $\CG$-module.

Given a finite 2-group $\CG$, its \emph{2-group algebra} $\vect_\CG$ is obtained by first linearizing the hom sets and then taking the Karoubi completion. One can show that $\vect_\CG$ is a multi-fusion category and the number of its indecomposable components is equal to the number of $\pi_1(\CG)$-orbits in $\widehat{\pi_2(\CG)}$ \cite{HZ23}. By the universal property of Karoubi completion, we have $2\rep(\CG) \simeq \LMod_{\vect_\CG}(2\vect)$ as linear 2-categories. Hence $2\rep(\CG)$ is a finite semisimple 2-category in the sense of Douglas and Reutter \cite{DR18}.

The symmetric monoidal structure of $2\rep(\CG)$ is inherited from $2\vect$. In other words, the tensor product of $\CV,\CW \in 2\rep(\CG)$ is the Deligne tensor product $\CV \boxtimes \CW$ equipped with the diagonal $\CG$-action. Moreover, $2\rep(\CG)$ is rigid and hence a symmetric fusion 2-category. Indeed, the left and right dual of $\CV \in 2\rep(\CG)$ are both $\CV^\op \simeq \fun_\bk(\CV,\vect)$. The evaluation 1-morphism $\ev \colon \CV^\op \boxtimes \CV \to \vect$ is the hom functor of $\CV$, and the coevaluation 1-morphism $\coev \colon \vect \to \CV \boxtimes \CV^\op$ is defined by
\[
\bk \mapsto \bigoplus_{x \in \Irr(\CV)} x \boxtimes x \simeq \int_{x \in \CV} x \boxtimes x .
\]

\section{Homotopy fixed points and homotopy quotients} \label{sec_homotopy_fixed_point_quotient}

\subsection{Homotopy fixed points}

\begin{defn}
Let $\CG$ be a 2-group and $\CC$ be a category equipped with a $\CG$-action. The \emph{homotopy fixed points} (or \emph{equivariantization}) of $\CC$ is the category $\CC^\CG$ defined as follows:
\bit
\item The objects in $\CC^\CG$ are pairs $(x,u = \{u_g\}_{g \in \CG})$, where $x \in \CC$ and $u_g \colon g \odot x \to x$ is an isomorphism in $\CC$ such that for any $g,h \in \CG$ and morphism $a \colon g \to h$ in $\CG$ the following diagrams commute:
\be \label{diag_equivariantization_object}
\begin{array}{c}
\xymatrix{
(g \otimes h) \odot x \ar[r]^{\simeq} \ar[d]_{u_{g \otimes h}} & g \odot (h \odot x) \ar[d]^{1 \odot u_h} \\
x & g \odot x \ar[l]_{u_g}
}
\end{array}
\quad
\begin{array}{c}
\xymatrix{
g \odot x \ar[r]^-{u_g} \ar[d]_{a \odot 1} & x \\
h \odot x \ar[ur]_{u_h}
}
\end{array}
\ee
\item A morphism $f \colon (x,u) \to (y,v)$ is a morphism $f \colon x \to y$ in $\CC$ such that for any $g \in \CG$ the following diagram commutes:
\be \label{diag_equivariantization_morphism}
\begin{array}{c}
\xymatrix{
g \odot x \ar[r]^{1 \odot f} \ar[d]_{u_g} & g \odot y \ar[d]^{v_g} \\
x \ar[r]^{f} & y
}
\end{array}
\ee
\item The composition and identity morphisms are induced by those of $\CC$.
\eit
Note that there is a forgetful functor $\CC^\CG \to \CC$ sending $(x,u)$ to $x$.
\end{defn}

\begin{rem}
The isomorphisms $u_g \colon g \odot x \to x$ can be equivalently defined by their mates $x \to g^* \odot x$. In the following we use these two formulations interchangeably.
\end{rem}

\begin{rem}
When $\CC$ is a monoidal category and the 2-group $\CG$ monoidally acts on $\CC$, i.e., every object in $\CG$ acts on $\CC$ as a monoidal equivalence and every morphism in $\CG$ acts on $\CC$ as a monoidal natural isomorphism, then the equivariantization $\CC^\CG$ is also a monoidal category, and the forgetful functor $\CC^\CG \to \CC$ is a monoidal functor. Indeed, the tensor product of two objects $(x,u),(y,v) \in \CC^\CG$ is the object $x \otimes y \in \CC$ equipped with the isomorphisms $g \odot (x \otimes y) \simeq (g \odot x) \otimes (g \odot y) \xrightarrow{u_g \otimes v_g} x \otimes y$.
\end{rem}

\begin{expl}
Let $\CG$ be a finite 2-group and $\CG \in 2\rep(\CG)$ be a finite semisimple 2-category. The equivariantization $\CV^\CG$ is equivalent to the hom category $2\rep(\CG)(\vect,\CV) = \fun_\CG(\vect,\CV)$, where $\vect \in 2\rep(\CG)$ is the trivial 2-representation.
\end{expl}

\begin{expl} \label{expl_Drinfeld_center_equivariantization}
Let $\CG$ be a finite 2-group. The conjugation $\CG$-action on itself, i.e., $g \odot h \coloneqq (g \otimes h) \otimes g^*$ for $g,h \in \CG$, induces a $\CG$-action on $\vect_\CG$. Then the equivariantization $(\vect_\CG)^\CG$ is equivalent to $\FZ_1(\vect_\CG)$. Indeed, an object in $(\vect_\CG)^\CG$ is an object $V \in \vect_\CG$ equipped with isomorphisms $u_g \colon g \odot V = (g \otimes V) \otimes g^* \to V$ for all $g \in \CG$ such that the diagrams in \eqref{diag_equivariantization_object} commute. The first diagram in \eqref{diag_equivariantization_object} is the same as the diagram \eqref{diag_Drinfeld_center_object_coend}, and the second diagram in \eqref{diag_equivariantization_object} is the dinaturality. Therefore, an object in $(\vect_\CG)^\CG$ is the same as an object in $\vect_\CG$ equipped with half-braidings with all $g \in \CG$, which can be uniquely extends to a half-braiding with all objects in $\vect_\CG$ because $\vect_\CG$ is the Karoubi completion of $\bk_\ast \CG$. Similarly, for a morphism in $(\vect_\CG)^\CG$ the diagram \eqref{diag_equivariantization_morphism} commutes, but this diagram is the same as the diagram \eqref{diag_Drinfeld_center_morphism_coend}. Hence we conclude that $(\vect_\CG)^\CG \simeq \FZ_1(\vect_\CG)$. Furthermore, the conjugation $\CG$-action on $\vect_\CG$ is monoidal, and one can verify that the $(\vect_\CG)^\CG \simeq \FZ_1(\vect_\CG)$ is a monoidal equivalence.
\end{expl}

\begin{defn}
Let $\CG$ be a 2-group, $\CC,\CD$ be categories equipped with $\CG$-actions and $F \colon \CC \to \CD$ be a $\CG$-module functor. Then the \emph{homotopy fixed point} (or \emph{equivariantization}) of $F$ is the functor $F^\CG \colon \CC^\CG \to \CD^\CG$ defined as follows:
\bit
\item For an object $(x,u = \{u_g\}_{g \in \CG}) \in \CC^\CG$, the object $F^\CG(x,u)$ is the object $F(x) \in \CD$ equipped with the isomorphisms
\[
g \odot F(x) \simeq F(g \odot x) \xrightarrow{F(u_g)} F(x) , \quad \forall g \in \CG ,
\]
where the first isomorphism is the $\CG$-module structure of $F$.
\item For a morphism $f \colon (x,u) \to (y,v)$ in $\CC^\CG$, the morphism $F^\CG(f) \coloneqq F$.
\eit
By abuse of notation, we also use $F \colon \CC^\CG \to \CD^\CG$ to denote the equivariantization $F^\CG$.
\end{defn}

We can also define the equivariantization of $\CG$-module natural transformations.

\subsection{Homotopy quotients}

\begin{defn}
Let $\CG$ be a 2-group and $\CC$ be a category equipped with a $\CG$-action. The \emph{homotopy quotient} of $\CC$ by $\CG$ is is the 2-category $\CC \dquotient \CG$ defined as follows:
\bit
\item The objects of $\CC \dquotient \CG$ are the objects of $\CC$.
\item A 1-morphism $(g,r) \colon x \to y$ consists of an object $g \in \CG$ and a morphism $r \colon g \odot x \to y$ in $\CC$.
\item A 2-morphism $a \colon (g,r) \Rightarrow (h,s) \colon x \to y$ is a morphism $a \colon g \to h$ in $\CG$ such that the following diagram commutes:
\[
\xymatrix{
g \odot x \ar[d]_{a \odot 1} \ar[r]^-{r} & y \\
h \odot x \ar[ur]_{s}
}
\]
\item The vertical composition of 2-morphisms is given by the composition of morphisms in $\CG$.
\item The identity 2-morphism on a 1-morphism $(g,r)$ is the identity morphism $1_g$ in $\CG$.
\item The horizontal composition of two 1-morphisms $(g,r) \colon x \to y$ and $(h,s) \colon y \to z$ is
\[
\bigl( h \otimes g , (h \otimes g) \odot x \simeq h \odot (g \odot x) \xrightarrow{1 \odot r} h \odot y \xrightarrow{s} z \bigr) .
\]
\item The horizontal composition of 2-morphisms is given by the tensor product of morphisms in $\CG$.
\item The identity 1-morphism on $x \in \CC$ is $(e,e \odot x \simeq x)$ where $e \in \CG$ is the tensor unit.
\item The associator and left/right unitors in $\CC \dquotient \CG$ are induced by those of $\CG$ and the $\CG$-module structure on $\CC$.
\eit
When $\CC$ is a groupoid, $\CC \dquotient \CG$ is a 2-groupoid and also called the \emph{action 2-groupoid}.
\end{defn}

For a category $\CC$ equipped with a $\CG$-action, there is a forgetful functor $\CC^\CG \to \CC$ and a 2-functor $\CC \to \CC \dquotient \CG$ sending a morphism $f \colon x \to y$ to $(e,e \odot x \simeq x \xrightarrow{f} y)$.

\begin{expl}
Let $X$ be a set and $G$ be a group acting on $X$. They can be viewed as a discrete category and a discrete 2-group with only identity morphisms, respectively. Then $X \dquotient G$ is the usual action groupoid.
\end{expl}

\begin{expl}
Let $\CG$ be a 2-group and $\ast$ be the trivial category which has only one object $\ast$ and only identity morphism. There is a trivial $\CG$-action on $\ast$, i.e., $g \odot \ast = \ast$ for any $g \in \CG$. Now we consider the homotopy quotient $\ast \dquotient \CG$.
\bit
\item There is only one object $\ast$.
\item A 1-morphism is a pair $(g,1_\ast)$ where $g \in \CG$.
\item A 2-morphism $(g,1_\ast) \Rightarrow (h,1_\ast)$ is a morphism $a \colon g \to h$ in $\CG$.
\eit
It is clear that $\ast \dquotient \CG$ is equivalent to the one-point delooping $\mathrm B \CG$.
\end{expl}

\begin{expl}
Let $\CG$ be a 2-group, acting on itself by left translation. Now we consider the homotopy quotient $\CG \dquotient \CG$.
\bit
\item The objects are the same as those in $\CG$.
\item A 1-morphism $(k,c) \colon g \to h$ consists of an object $k \in \CG$ and a morphism $c \colon k \otimes g \to h$.
\item A 2-morphism $a \colon (k,c) \Rightarrow (l,d) \colon g \to h$ is a morphism $a \colon k \to l$ in $\CG$ such that the following diagram commutes:
\[
\xymatrix{
k \otimes g \ar[d]_{a \otimes 1} \ar[r]^-{c} & h \\
l \otimes g \ar[ur]_{d}
}
\]
\eit
Since $- \otimes g \colon \CG \to \CG$ is an equivalence for any $g \in \CG$, there is a unique 2-morphism between two 1-morphisms $(k,c) , (l,d) \colon g \to h$. Moreover, for any $g,h \in \CG$ there exists at least one 1-morphism $g \to h$, for example:
\[
\bigl( h \otimes g^* , h \otimes g^* \otimes g \xrightarrow{1 \otimes \ev_g} h \bigr) .
\]
Hence, the 2-groupoid $\CG \dquotient \CG$ is contractible. We denote this 2-groupoid by $\mathrm E \CG$.

Moreover, there is a right $\CG$-action on $\mathrm E \CG$ induced by the right translation of $\CG$ on itself. Then the homotopy quotient $\mathrm E \CG \dquotient \CG$ (or more precisely $\mathrm E \CG \dquotient \CG^\rev$) is equivalent to $\mathrm B \CG$.
\end{expl}

%\begin{rem}
%The geometric realization of $\mathrm B \CG$ is called the \emph{classifying space} of $\CG$. Its homotopy groups are the same as those of $\CG$. It may be interesting to consider the geometric realization of $\mathrm E \CG$. It should a contractible space with a `free right $\CG$-action' and should be viewed as a `universal pincipal $\CG$-bundle'. However, I do not know what is a 2-group action on a topological space, nor what is a principal $\CG$-bundle on a topological space.
%\end{rem}

\begin{expl} \label{expl_conjugation_2-groupoid_loop_2-groupoid}
Let $\CG$ be a 2-group, acting on itself by conjugation, i.e., $g \odot h \coloneqq (g \otimes h) \otimes g^*$ for any $g,h \in \CG$. Now we consider the homotopy quotient $\CG \dquotient \CG$.
\bit
\item The objects are the same as those in $\CG$.
\item A 1-morphism $(k,c) \colon g \to h$ consists of an object $k \in \CG$ and a morphism $c \colon (k \otimes g) \otimes k^* \to h$.
\item A 2-morphism $a \colon (k,c) \Rightarrow (l,d) \colon g \to h$ is a morphism $a \colon k \to l$ in $\CG$ such that the following diagram commutes:
\be \label{eq_conjugate_action_2-groupoid_2-morphism}
\begin{array}{c}
\xymatrix{
(k \otimes g) \otimes k^* \ar[d]_{(a \otimes 1) \otimes a^*} \ar[r]^-{c} & h \\
(l \otimes g) \otimes l^* \ar[ur]_{d}
}
\end{array}
\ee
\eit
We also consider the \emph{loop 2-groupoid} (also called the \emph{inertia 2-groupoid}) $\mathrm L \mathrm B \CG \coloneqq \fun(\mathrm B \Zb,\mathrm B \CG)$ consisting of 2-functors, 2-natural transformations and modifications.
\bit
\item A 2-functor $F \colon \mathrm B \Zb \to \mathrm B \CG$ is the same as a monoidal functor $F \colon \Zb \to \CG$, which consists of a set of objects $\{F(n) \in \CG\}_{n \in \Zb}$ and a 2-cocycle $(F^2_{n,m}) \in Z^2(\Zb;\pi_2(\CG))$. Here the $\Zb$-action on $\pi_2(\CG)$ is induced by the $\pi_1(\CG)$-action on $\pi_2(\CG)$ and the group homomorphism $F \colon \Zb \to \pi_1(\CG)$. Since $H^2(\Zb;M) = 0$ for any $\Zb$-module $M$, we can assume that $F$ is a strict monoidal functor and only depends on $F(1) \in \CG$. For every $g \in \CG$, we use $L^g \colon \Zb \to \CG$ to denote the strict monoidal functor with $L^g(n) = g^{\otimes n}$. By the above argument, every 2-functor $\mathrm B \Zb \to \mathrm B \CG$ is equivalent to $L^g$ for some $g \in \CG$.
\item A 2-natural transformation $\alpha \colon L^g \Rightarrow L^h$ for $g,h \in \CG$ consists of a 1-morphism $\alpha_\ast \colon \ast \to \ast$ in $\mathrm B \CG$ and a 2-morphism $\alpha_n$ in $\mathrm B \CG$ for any $n \in \Zb$ as depicted in the following diagram:
\[
\xymatrix{
\ast \ar[r]^{L^g(n)} \ar[d]_{\alpha_\ast} & \ast \ar[d]^{\alpha_\ast} \\
\ast \ar[r]_{L^h(n)} & \ast \ultwocell<\omit>{\alpha_n}
}
\]
By the naturality, $\alpha$ is determined by $\alpha_\ast \in \CG$ and the morphism $\alpha_1 \colon \alpha_\ast \otimes g \to h \otimes \alpha_\ast$ in $\CG$. Since $\alpha_\ast \in \CG$ is invertible, $\alpha_1$ is also equivalent to a morphism $\tilde \alpha_1 \colon (\alpha_\ast \otimes g) \otimes \alpha_\ast^* \to h$.
\item A modification $\Gamma \colon \alpha \Rrightarrow \beta \colon L^g \Rightarrow L^h$ consists of a 2-morphism $\Gamma_\ast \colon \alpha_\ast \Rightarrow \beta_\ast$ in $\mathrm B \CG$ satisfying a condition. So $\Gamma_\ast \colon \alpha_\ast \to \beta_\ast$ is a morphism in $\CG$. The modification axiom for $\Gamma$ means that the following diagram commutes:
\[
\xymatrix{
\alpha_\ast \otimes g \ar[r]^-{\alpha_1} \ar[d]_{\Gamma_\ast \otimes 1} & h \otimes \alpha_\ast \ar[d]^{1 \otimes \Gamma_\ast} \\
\beta_\ast \otimes g \ar[r]^-{\beta_1} & h \otimes \beta_\ast
}
\]
It is easy to see that this is equivalent to the commutativity of the diagram \eqref{eq_conjugate_action_2-groupoid_2-morphism}.
\eit
By comparing these two 2-categories, it is easy to conclude that the action 2-groupoid $\CG \dquotient \CG$ (with the conjugate action) is equivalent to the loop 2-groupoid $\mathrm L \mathrm B \CG$.
\end{expl}

\begin{rem} \label{rem_S1_action_LG}
There is a canonical $\mathrm B \Zb$-action on the 2-groupoid $\CG \dquotient \CG \simeq \mathrm L \mathrm B \CG = \fun(\mathrm B \Zb,\mathrm B \CG)$ induced from the regular $\mathrm B \Zb$-action on itself. Note that a $\mathrm B \Zb$-action on $\CG \dquotient \CG$ is equivalent to by a 2-natural autoequivalence $\phi \in \Omega \Aut(\CG \dquotient \CG) = \Aut(1_{\CG \dquotient \CG})$. This 2-natural autoequivalence $\phi$ is defined as follows:
\bit
\item For any $g \in \CG \dquotient \CG$, the 1-morphism $\phi_g \colon g \to g$ is $g \in \CG$ equipped with the morphism
\[
(g \otimes g) \otimes g^* \simeq g \otimes (g \otimes g^*) \xrightarrow{1 \otimes \ev} g \otimes e \simeq g .
\]
\item For any morphism $(k,c) \colon g \to h$ in $\CG \dquotient \CG$, the 2-isomorphism $\phi_{(k,c)}$ as depicted in the following diagram
\[
\xymatrix{
g \ar[r]^-{(k,c)} \ar[d]_{\phi_g} & h \ar[d]^{\phi_h} \\
g \ar[r]^-{(k,c)} & h \ultwocell<\omit>{\;\;\phi_{(k,c)}\quad}
}
\]
is given by the composite morphism $k \otimes g \xrightarrow{1 \otimes \coev} (k \otimes g) \otimes (k^* \otimes k) \simeq ((k \otimes g) \otimes k^*) \otimes k \xrightarrow{c \otimes 1} h \otimes k$.
\eit
\end{rem}

\begin{expl}
Let $\CG$ be a 2-group. Note that there is a canonical monoidal functor $\CG \to \pi_1(\CG)$. So there is a $\CG$-action on $\pi_2(\CG)$ (also by conjugation). Now we consider the action 2-groupoid $\pi_2(\CG) \dquotient \CG$.
\bit
\item The objects are elements in $\pi_2(\CG)$.
\item A 1-morphism $(g,1) \colon a \to b$ is of an object $g \in \CG$ such that $g \triangleright a = b$. Recall that $g \triangleright a = 1_g \otimes a \otimes 1_{g^*}$.
\item A 2-morphism $(g,1) \Rightarrow (h,1) \colon a \to b$ is a morphism $g \to h$ in $\CG$ satisfying no condition.
\eit
We are also interested in the \emph{sphere 2-groupoid} $\mathrm S \CG \coloneqq \fun(\mathrm B^2 \Zb,\mathrm B \CG)$ of $\CG$. We denote the only object in $\mathrm B^2 \Zb$ by $\ast$ and denote the only 1-morphism in $\mathrm B^2 \Zb$ by $\bullet$.
\bit
\item A 2-functor $F \colon \mathrm B^2 \Zb \to \mathrm B \CG$ is the same as a monoidal functor $F \colon \mathrm B \Zb \to \CG$, and also the same as a group homomorphism $\Zb \to \pi_2(\CG)$. So for every $a \in \pi_2(\CG)$ there is a 2-functor $S^a \colon \mathrm B^2 \Zb \to \mathrm B \CG$ that maps $1 \in \Zb$ to $a$, and these are the only 2-functors.
\item A 2-natural transformation $\alpha \colon S^a \Rightarrow S^b$ for $a,b \in \pi_2(\CG)$ consists of a 1-morphism $\alpha_\ast \colon \ast \to \ast$ in $\mathrm B \CG$ and a 2-morphism $\alpha_\bullet$ in $\mathrm B \CG$ as depicted in the following diagram:
\[
\xymatrix{
\ast \ar[r]^{S^a(\bullet) = e} \ar[d]_{\alpha_\ast} & \ast \ar[d]^{\alpha_\ast} \\
\ast \ar[r]_{S^b(\bullet) = e} & \ast \ultwocell<\omit>{\alpha_\bullet}
}
\]
By the naturality, $\alpha_\bullet$ is trivial and $\alpha_\ast$ satisfies the condition that
\[
1_{\alpha_\ast} \circ a = b \circ 1_{\alpha_\ast} .
\]
So $\alpha$ is determined by the object $\alpha_\ast \in \CG$ satisfying $\alpha_\ast \triangleright a = b$.
\item A modification $\Gamma \colon \alpha \Rrightarrow \beta$ consists of a 2-morphism $\Gamma_\ast \colon \alpha_\ast \Rightarrow \beta_\ast$ in $\mathrm B \CG$ satisfying a condition. So $\Gamma_\ast \colon \alpha_\ast \to \beta_\ast$ is a morphism in $\CG$. It is easy to see that the modification axiom for $\Gamma$ is trivial.
\eit
By comparing these two 2-categories, it is easy to see that the action 2-groupoid $\pi_2(\CG) \dquotient \CG = \Omega \CG \dquotient \CG$ (with the conjugate action) is equivalent to the sphere 2-groupoid $\mathrm S \mathrm B\CG$.
\end{expl}

\subsection{The duality between homotopy fixed points and homotopy quotients}

\begin{prop} \label{prop_duality_homotopy_fixed_point_quotient}
Let $\CG$ be a 2-group, $\CC$ be a category equipped with a $\CG$-action and $\CD$ be a category. Then there is a canonical equivalence of categories
\[
\fun(\CC \dquotient \CG,\CD) \simeq \fun(\CC,\CD)^\CG ,
\]
where the $\CG$-action on $\fun(\CC,\CD)$ is induced by that of $\CC$.
\end{prop}

\pf
First we construct the functor $\fun(\CC \dquotient \CG,\CD) \to \fun(\CC,\CD)^\CG$. Suppose $F \colon \CC \dquotient \CG \to \CD$ is a functor. The composite functor $\CC \to \CC \dquotient \CG \to \CD$ is denoted by $\overline F$. Note that for any $g \in \CG$ and $x \in \CC$, there is an isomorphism $(g,\ev_g \odot 1_x) \colon g^* \odot x \to x$ in $\CC \dquotient \CG$. Then we define
\[
u_{g,x} \coloneqq F(g,\ev_g \odot 1) \colon (g \odot \overline F)(x) = F(g^* \odot x) \to F(x) = \overline F(x) .
\]
Thus we obtain a natural isomorphism $u_g \colon g \odot \overline F \Rightarrow \overline F$ for every $g \in \CG$. We need to check that the diagrams in \eqref{diag_equivariantization_object} are commutative. It is not hard to compute that the composite 1-morphism
\[
(g \otimes h)^* \odot x \simeq h^* \odot (g^* \odot x) \xrightarrow{(h,\ev_h \odot 1)} g^* \odot x \xrightarrow{(g,\ev_g \odot 1)} x
\]
in $\CC \dquotient \CG$ is equal to $(g \otimes h,\ev_{g \otimes h} \odot 1_x) \colon (g \otimes h)^* \odot x \to x$. Then applying $F$ on this equality we obtain the first commutative diagram. For the second commutative diagram, note that there is a 2-morphism in $\CC \dquotient \CG$ as depicted in the following diagram:
\[
\xymatrix{
g^* \odot x \ar[rr]^{(g,\ev_g \odot 1)} \ar[dr]_{(e,a^* \odot 1)} \rrtwocell<\omit>{<3> a} & & x \\
 & h^* \odot x \ar[ur]_{(h,\ev_h \odot 1)}
}
\]
Then applying $F$ we obtain the second commutative diagram.

If $\phi \colon F \Rightarrow K \colon \CC \dquotient \CG \to \CD$ is a natural transformation, then by composing with the functor $\CC \to \CC \dquotient \CG$ we obtain a natural transformation $\overline \phi \colon \overline F \Rightarrow \overline K \colon \CC \to \CD$. Then we need to show that the diagram \eqref{diag_equivariantization_morphism} commutes. This is due to the naturality of $\phi$. Hence we obtain the functor $\fun(\CC \dquotient \CG,\CD) \to \fun(\CC,\CD)^\CG$ defined by $F \mapsto (\overline F,\{u_g\}_{g \in \CG})$.

Then we construct the functor $\fun(\CC,\CD)^\CG \to \fun(\CC \dquotient \CG,\CD)$. Suppose $(F,\{u_g\}_{g \in \CG}) \in \fun(\CC,\CD)^\CG$. Then we define a functor $\underline F \colon \CC \dquotient \CG \to \CD$ as follows:
\bit
\item For any $x \in \CC \dquotient \CG$, we define $\underline F(x) \coloneqq F(x)$.
\item For any morphism $(g,f) \colon x \to y$ in $\CC \dquotient \CG$, we define
\[
\underline F(g,f) \coloneqq \biggl( F(x) \simeq F(g^* \odot (g \odot x)) \xrightarrow{u_{g,g \odot x}} F(g \odot x) \xrightarrow{F(f)} F(y) \biggr) .
\]
\eit
If $\phi \colon (F,u) \Rightarrow (K,v)$ is a $\CG$-equivariant natural transformation, then we define $\underline \phi \colon \underline F \Rightarrow \underline K$ by $\underline \phi_x \coloneqq \phi_x$ for any $x \in \CC \dquotient \CG$. The naturality of $\underline \phi$ is due to that of $\phi$ and the commutative diagram \eqref{diag_equivariantization_morphism} for $\phi$.

Finally, it is not hard to verify that these two functors between $\fun(\CC \dquotient \CG,\CD)$ and $\fun(\CC,\CD)^\CG$ are mutually inverse to each other.
\epf

\begin{rem}
Furthermore, for a 2-group $\CG$ and a 2-category $\SC$ equipped with $\CG$-actions, we can also define the homotopy fixed points $\SC^\CG$ and homotopy quotients $\SC \dquotient \CG$, which are 2-categories. For example, for the trivial $\CG$-action on $2\vect$ we have $2\vect^\CG \simeq 2\rep(\CG)$. We should also have an equivalence of 2-categories $\fun(\SC \dquotient \CG,\SD) \simeq \fun(\SC,\SD)^\CG$.
\end{rem}

\begin{expl}
Recall that there is a $\mathrm B \Zb$-action on $\CG \dquotient \CG \simeq \fun(\mathrm B \Zb,\mathrm B \CG)$ defined in Remark \ref{rem_S1_action_LG}. Then we have an equivalence of 2-categories $(\CG \dquotient \CG)^{\mathrm B \Zb} \simeq \fun(\mathrm B \Zb \dquotient \mathrm B \Zb,\mathrm B \CG) \simeq \fun(\ast,\mathrm B \CG) \simeq \mathrm B \CG$.
\end{expl}

\subsection{Fourier 2-transform} \label{sec_Fourier_2-transform}

Let $\CG$ be a finite 2-group. There is a canonical equivalence
\[
\vect_\CG^\op \simeq \fun_\bk(\vect_\CG,\vect) \simeq \fun(\CG,\vect) ,
\]
where the first equivalence is the Yoneda embedding and the second equivalence is the universal property of Karoubi completion. However, there is no `canonical' equivalence $\vect_\CG \simeq \fun(\CG,\vect)$. This is similar to the fact that a vector space is not canonically isomorphic to its dual space. There are many different equivalences $\vect_\CG \simeq \fun(\CG,\vect)$. Now we are going to fix a choice.

First, there is a functor
\begin{align*}
\CG & \to \fun(\CG,\vect) \\
g & \mapsto \vect_\CG(I(-),g) ,
\end{align*}
where $I \colon \CG^\op \to \CG$ is the inverse morphism functor that sends each morphism to its inverse. By the universal property of Karoubi completion, it induces a $\bk$-linear functor $\Phi \colon \vect_\CG \to \fun(\CG,\vect)$. In other words, $\Phi(V) = \vect_\CG(I(-),V)$ for all $V \in \vect_\CG$. There is another obvious $\bk$-linear functor $\Phi' \colon \vect_\CG \to \fun(\CG,\vect)$ defined by $\Phi'(V) \coloneqq \vect_\CG(I(V),-)$ for every $V \in \vect_\CG$. Since $I^2$ is identity, $\Phi$ and $\Phi'$ are naturally isomorphic.

To show that $\Phi$ and $\Phi'$ are equivalence, note that $\Phi$ is equal to the composite functor
\[
\vect_\CG \simeq \vect_\CG^\op \simeq \fun_\bk(\vect_\CG,\vect) \simeq \fun(\CG,\vect) ,
\]
where the first equivalence is induced by $I$. We also give an explicit construction of the quasi-inverse of $\Phi$.

\begin{prop}
Viewed $\vect_\CG$ as a $\vect$-module. Then the functor
\begin{align*}
\Psi \colon \fun(\CG,\vect) & \to \vect_\CG \\
F & \mapsto \int^{g \in \CG} F(g) \odot I(g)
\end{align*}
is the quasi-inverse of $\Phi \colon \vect_\CG \to \fun(\CG,\vect)$.
\end{prop}

\pf
For every $g,x \in \CG$, we have
\begin{multline*}
\vect_\CG((\Psi \circ \Phi)(g),x) = \vect_\CG \biggl(\int^{h \in \CG} \vect_\CG(I(h),g) \odot I(h),x \biggr) \\
\simeq \int_{h \in \CG} \vect_\CG(\vect_\CG(I(h),g) \odot I(h),x) 
\simeq \int_{h \in \CG} \vect(\vect_\CG(I(h),g),\vect_\CG(I(h),x)) \\
\simeq \nat(\vect_\CG(I(-),g),\vect_\CG(I(-),x)) \simeq \vect_\CG(g,x) .
\end{multline*}
Then by the Yoneda lemma $(\Psi \circ \Phi)(g) \simeq g$. Conversely, for every $F \in \fun(\CG,\vect)$, we have
\begin{multline*}
(\Phi \circ \Psi)(F)(g) = \vect_\CG \biggl(I(g),\int^{h \in \CG} F(h) \odot I(h) \biggr) \simeq \int^{h \in \CG} \vect_\CG(I(g),F(h) \odot I(h)) \\
\simeq \int^{h \in \CG} F(h) \otimes \vect_\CG(I(g),I(h)) \simeq \int^{h \in \CG} F(h) \otimes \vect_\CG(h,g) \simeq F(g) ,
\end{multline*}
where the first isomorphism is due to the fact that $\vect_\CG$ is semisimple and the hom functor is exact in each variable.
\epf

Since $\vect_\CG$ is a multi-fusion category, its monoidal structure can be transferred to $\fun(\CG,\vect)$ so that $\Phi$ and $\Psi$ are monoidal equivalences. This monoidal structure on $\fun(\CG,\vect)$ is given by the Day convolution \cite{Day71}:
\[
(F \circledast G)(g) \coloneqq \int^{h,k \in \CG} \vect_\CG(h \otimes k,g) \otimes F(h) \otimes G(k) \simeq \int^{k \in \CG} F(g \otimes k^L) \otimes G(k) , \quad F,G \in \fun(\CG,\vect) .
\]
Indeed, for every $g,h,x \in \CG$ we have
\begin{multline*}
(\Phi(g) \circledast \Phi(h))(x) = \int^{k \in \CG} \vect_\CG(I(x \otimes k^L),g) \otimes \vect_\CG(I(k),h) \\
\simeq \int^{k \in \CG} \vect_\CG(I(x) \otimes k^*,g) \otimes \vect_\CG(I(h),k) \\
\simeq \vect_\CG(I(x) \otimes h^R,g) \simeq \vect_\CG(I(x),g \otimes h) = \Phi(g \otimes h)(x) ,
\end{multline*}
where we use the fact that the left and right dual of $h \in \CG$ coincide. It is not hard to see that this isomorphism is also the unique morphism rendering the following diagram commutative for all $y \in \CG$:
\be \label{diag_monoidal_structure_Fourier_2-transform}
\begin{array}{c}
\xymatrix{
\int^{k \in \CG} \vect_\CG(I(x \otimes k^L),g) \otimes \vect_\CG(I(k),h) \ar@{-->}[r] & \vect_\CG(I(x),g \otimes h) \\
\vect_\CG(I(x \otimes y^L),g) \otimes \vect_\CG(I(y),h) \ar[r]^-{\otimes} \ar[u]^{\tau_y} & \vect_\CG(I(x \otimes y^L \otimes y),g \otimes h) \ar[u]_{\simeq}
}
\end{array}
\ee

\medskip
Now we equip $\vect_\CG$ and $\fun(\CG,\vect)$ the conjugation $\CG$-action. Note that for every $g,h \in \CG$ we have
\begin{multline} \label{eq_G-module_functor_Phi}
\Phi((g \otimes h) \otimes g^*) = \vect_\CG(I(-),(g \otimes h) \otimes g^*) \simeq \vect_\CG((g^L \otimes I(-)) \otimes I(g),h) \\
\simeq \vect_\CG(I((g^* \otimes -) \otimes g),h) = \Phi(h)((g^* \otimes -) \otimes g) = g \odot \Phi(h) .
\end{multline}
Therefore, $\Phi$ and $\Psi$ are $\CG$-module equivalences. Moreover, the conjugation $\CG$-action is monoidal. So after taking the equivariantization we have monoidal equivalences between $\fun(\CG,\vect)^\CG$ and $\FZ_1(\vect_\CG) \simeq (\vect_\CG)^\CG$ (see Example \ref{expl_Drinfeld_center_equivariantization}), still denoted by $\Phi$ and $\Psi$.

Let us write down the equivalences $\Phi \colon \FZ_1(\vect_\CG) \to \fun(\CG,\vect)^\CG$ and $\Psi \colon \fun(\CG,\vect)^\CG \to \FZ_1(\vect_\CG)$ more explicitly. For an object $(V,u = \{u_g \colon (g \otimes V) \otimes g^* \to V\}_{g \in \CG}) \in (\vect_\CG)^\CG \simeq \FZ_1(\CG)$, its image under $\Phi$ is the functor $\Phi(V) = \vect_\CG(I(-),V)$ equipped with the $\CG$-equivariant structure
\[
g \odot \Phi(V) \simeq \Phi((g \otimes V) \otimes g^*) \xrightarrow{\Phi(u_g)} \Phi(V) ,
\]
where the first isomorphism is the $\CG$-module functor structure \eqref{eq_G-module_functor_Phi} of $\Phi$. Conversely, for an object $(F,\phi = \{\phi_g \colon g \odot F \Rightarrow F\}_{g \in \CG}) \in \fun(\CG,\vect)^\CG$, its image under $\Psi$ is the object $\Psi(F) \in \vect_\CG$ equipped with the $\CG$-equivariant structure
\be \label{eq_half_braiding_Psi}
(g \otimes \Psi(F)) \otimes g^* = \int^{h \in \CG} F(h) \odot ((g \otimes I(h)) \otimes g^*) \xrightarrow{\int^h \phi_g \odot 1} \int^{h \in \CG} F((g \otimes h) \otimes g^*) \odot I((g \otimes h) \otimes g^*) = \Psi(F) .
\ee
Furthermore, the braiding structure on $\FZ_1(\vect_\CG)$ is also transferred to $\fun(\CG,\vect)^\CG$. Let us describe this braiding structure. For $(F,\phi),(G,\psi) \in \fun(\CG,\vect)^\CG$, the braiding is a natural isomorphism
\[
(\beta_{F,G})_g \colon (F \circledast G)(g) \to (G \circledast F)(g) , \quad g \in \CG .
\]
This isomorphism is the unique morphism rendering the following diagram commutative for all $x \in \CG$:
\[
\xymatrix{
\int^{h \in \CG} F(g \otimes h^L) \otimes G(h) \ar@{-->}[rr]^-{(\beta_{F,G})_g} & & \int^{k \in \CG} G(g \otimes k^L) \otimes F(k) \\
F(g \otimes x^L) \otimes G(x) \ar[r]^{\simeq} \ar[u]^{\tau_x} & G(x) \otimes F(g \otimes x^L) \ar[r]^-{\psi_g \otimes 1} & G((g \otimes x) \otimes g^*) \otimes F(g \otimes x^L) \ar[u]_{\tau_{g \otimes x^L}}
}
\]
where the unlabeled isomorphism is given by switching two vector spaces.

\begin{defn}
The equivalences $\Phi \colon \FZ_1(\vect_\CG) \to \fun(\CG,\vect)^\CG$ and $\Psi \colon \fun(\CG,\vect)^\CG \to \FZ_1(\vect_\CG)$ defined above are called the \emph{Fourier 2-transform} of $\CG$.
\end{defn}

\begin{rem} \label{rem_Fourier_transform}
Let us explain why these two equivalences are called ``Fourier 2-transform'' by considering the de-categorified case. Suppose $G$ is a finite group. Then we have two isomorphisms
\begin{align*}
\Phi \colon \bk[G] & \to \fun(G) & \Psi \colon \fun(G) & \to \bk[G] \\
g & \mapsto \delta_g & f & \mapsto \sum_{g \in G} f(g) \cdot g
\end{align*}
where $\delta_g$ is the delta function defined by $\delta_g(h) = \delta_{g,h}$. Moreover, they are algebra homomorphisms if we equip $\fun(G)$ with the convolution product:
\[
(s \ast t)(g) \coloneqq \sum_{\substack{x,y \in G \\ xy=g}} s(x)t(y) = \sum_{y \in G} s(gy^{-1}) t(y) .
\]
If we equip $\bk[G]$ and $\fun(G)$ with the conjugation $G$-action, $\Phi$ and $\Psi$ are also $G$-equivariant maps. So after taking the $G$-fixed points, we obtain two isomorphisms of commutative algebras $\Phi \colon Z(\bk[G]) \simeq \bk[G]^G \to (\fun(G)^G,\ast)$ and $\Psi \colon (\fun(G)^G,\ast) \to Z(\bk[G])$.

When $G$ is abelian, $\bk[G]$ is a commutative algebra and are naturally isomorphic to $\fun(\hat G)$ with the point-wise multiplication, where $\hat G$ is the dual group of $G$. It is easy to see that the isomorphism $(\fun(G),\ast) \to (\fun(\hat G),\cdot)$ sends a function $f \in \fun(G)$ to the function $\hat f \in \fun(\hat G)$ defined by
\[
\hat f(\rho) = \sum_{g \in G} \rho(g) f(g) ,
\]
and its inverse $(\fun(\hat G),\cdot) \to (\fun(G),\ast)$ sends a function $\varphi \in \fun(\hat G)$ to the function $\check \varphi \in \fun(G)$ defined by
\[
\check \varphi(g) = \frac{1}{\lvert G \rvert} \sum_{\rho \in \hat G} \rho(g^{-1}) \varphi(\rho) .
\]
These are the Fourier transform on the finite abelian group $G$.

In general, $G$ may not be abelian. Then $Z(\bk[G])$ can be identified with the space of functions on the set $\Irr(\rep(G))$ of irreducible $G$-representations, and the equivalence between $Z(\bk[G])$ and $(\fun(G)^G,\ast)$ can be viewed as a non-abelian Fourier transform.
\end{rem}

Combining with Proposition \ref{prop_duality_homotopy_fixed_point_quotient} and Example \ref{expl_conjugation_2-groupoid_loop_2-groupoid}, we obtain the following result.

\begin{cor} \label{cor_fun_conjugation_Drinfeld_center}
Let $\CG$ be a finite 2-group, acting on itself by conjugation. Then we have $\fun(\CG \dquotient \CG,\vect) \simeq \FZ_1(\vect_\CG)$. Moreover, they are both equivalent to the space $\fun(\mathrm{L} \mathrm B \CG,\vect)$ of $\vect$-valued functors on the loop 2-groupoid. In particular, there is a natural braided monoidal structure on $\fun(\mathrm L \mathrm B \CG,\vect)$.
\end{cor}

\section{The 2-characters for finite 2-groups} \label{sec_2-character}

\subsection{Definition and basic properties}

Let $\CG$ be a 2-group and $\CV \in 2\rep(\CG)$. For every $g \in \CG$, we define
\be \label{eq_chi_V_g}
\chi_\CV(g) \coloneqq \bigl(\vect \xrightarrow{\coev} \CV^\op \boxtimes \CV \xrightarrow{1 \boxtimes (g \odot -)} \CV^\op \boxtimes \CV \xrightarrow{\ev} \vect\bigr) ,
\ee
where the evaluation and coevaluation 1-morphisms are given in Section \ref{sec_2-group}. Intuitively, $\chi_\CV(g)$ is the `trace' of the action by $g$. Since $\End_\bk(\vect) \simeq \vect$ as fusion categories, $\chi_\CV(g)$ can be viewed as a finite-dimensional vector space. It follows that $\chi_\CV$ defines a functor $\chi_\CV \colon \CG \to \vect$, called the \emph{2-character} of $\CV$. More explicitly, by direct computation we have
\[
\chi_\CV(g) \simeq \bigoplus_{v \in \Irr(\CV)} \CV(v,g \odot v) .
\]
Since $\CV$ is finite semisimple, by Corollary \ref{cor_semisimple_end_direct_sum} we also have
\[
\chi_\CV(g) \simeq \int_{v \in \CV} \CV(v,g \odot v) \simeq \nat(1_\CV,g \odot -) .
\]
In the following, we also take this as the definition of $\chi_\CV(g)$.

\begin{rem} \label{rem_categorical_trace_GK}
Let $\SC$ be a 2-category and $f \colon x \to x$ be a 1-endomorphism in $\SC$. In \cite{GK08}, Ganter and Kapranov defined the ``categorical trace'' of $f$ to be $\Hom_{\Hom_\SC(x,x)}(1_x,f)$. So when $\SC = 2\vect$, the categorical trace of the functor $g \odot - \colon \CV \to \CV$ coincides with the 2-character $\chi_\CV(g)$ defined as above. Ganter and Kapranov also studied the 2-character in the special case that $\CG$ is an ordinary group.
\end{rem}

\begin{expl}
For the trivial 2-representation $\vect$, the 2-character $\chi_\vect \colon \CG \to \vect$ is the constant functor on $\bk \in \vect$. More generally, let $\phi \colon \CG \to \vect$ be a monoidal functor and $\vect_\phi$ be the 2-representation defined by $\CG \xrightarrow{\phi} \vect \simeq \End_\bk(\vect)$. Then $\vect_\phi \in 2\rep(\CG)$ is an invertible object and its associated 2-character is $\phi$.
\end{expl}

\begin{prop} \label{prop_2-character_direct_sum_tensor_product}
Let $\CV,\CW \in 2\rep(\CG)$. Then we have the following results for $g \in \CG$:
\bnu[(1)]
\item $\chi_{\CV \oplus \CW}(g) \simeq \chi_\CV(g) \oplus \chi_\CW(g)$.
\item $\chi_{\CV \boxtimes \CW}(g) \simeq \chi_\CV(g) \otimes \chi_\CW(g)$.
\item $\chi_{\CV^\op}(g) \simeq \chi_{\CV}(g^*)$.
\enu
\end{prop}

\pf
By direct computation we have
\[
\chi_{\CV \oplus \CW}(g) = \bigoplus_{v \in \Irr(\CV)} \CV(v,g \odot v) \oplus \bigoplus_{w \in \Irr(\CW)} \CW(w,g \odot w) \simeq \chi_\CV(g) \oplus \chi_\CW(g) ,
\]
and
\begin{align*}
\chi_{\CV \boxtimes \CW}(g) & = \bigoplus_{v \in \Irr(\CV)} \bigoplus_{w \in \Irr(\CW)} \CV \boxtimes \CW(v \boxtimes w,g \odot (v \boxtimes w)) \\
& = \bigoplus_{v \in \Irr(\CV)} \bigoplus_{w \in \Irr(\CW)} \CV \boxtimes \CW(v \boxtimes w,(g \odot v) \boxtimes (g \odot w)) \\
& \simeq \bigoplus_{v \in \Irr(\CV)} \bigoplus_{w \in \Irr(\CW)} \CV(v,g \odot v) \otimes \CW(w,g \odot w) \simeq \chi_\CV(g) \otimes \chi_\CW(g) .
\end{align*}
and
\[
\chi_{\CV^\op}(g) = \bigoplus_{v \in \Irr(\CV^\op)} \CV^\op(v,g \odot v) = \bigoplus_{v \in \Irr(\CV)} \CV(g \odot v,v) \simeq \bigoplus_{v \in \Irr(\CV)} \CV(v,g^* \odot v) = \chi_{\CV}(g^*).
\]
This is indeed the natural isomorphism of two functors $\chi_{\CV^\op} \simeq \chi_\CV \circ (-)^*$.
\epf

\begin{rem} \label{rem_G_TQFT}
The 2-character $\chi_\CV$ and its properties can be understood from a TQFT point of view. However, all of our discussion of TQFT in several remarks is intended only as motivation. We do not aim to construct the relevant TQFT nor seek to give rigorous proofs of every statement.

A key observation comes from the cobordism hypothesis, which was proposed by Baez and Dolan \cite{BD95}. Later, Lurie gave an outline of a proof in \cite{Lur08}, and Grady and Pavlov gave a proof in \cite{GP21}. The cobordism hypothesis says that a 2D extended framed TQFT, that is, a symmetric monoidal 2-functor from the framed 2D cobordism 2-category to another 2-category, is determined by its image on a single point, which must be a fully dualizable object. Here we consider a $\CG$-equivariant version (see also \cite[Theorem 2.4.26]{Lur08}).

A \emph{2D extended $\CG$-TQFT} is a symmetric monoidal 2-functor $Z \colon \mathrm{Cob}_2^\CG \to 2\vect$, where $\mathrm{Cob}_2^\CG$ is the 2-category of 2D cobordisms equipped with $\CG$-structures, i.e., continuous maps to the classifying space $\lvert \mathrm B \CG \rvert$ of $\CG$. The classifying space $\lvert \mathrm B \CG \rvert$ can be defined as the geometric realization of the Duskin nerve of the 2-groupoid $\mathrm B \CG$. The cobordism hypothesis suggests that a 2D extended $\CG$-TQFT is determined by its valued on a single point, which is a finite semisimple 2-representation $\CV \in 2\rep(\CG) \simeq 2\vect^\CG$ (recall that it is always fully dualizable). Such a TQFT is denoted by $Z_\CV \colon \mathrm{Cob}_2^\CG \to 2\vect$. It is also a topological sigma model with the target space $\lvert \mathrm B \CG \rvert$.

In the 2-category $\mathrm{Cob}_2^\CG$, a circle $S^1$ equipped with a $\CG$-structure is a 1-morphism from the empty set $\emptyset$ to itself. What does a $\CG$-structure on $S^1$ look like? First we fix a point $\ast \in S^1$ and consider \emph{pointed} continuous maps $f \colon S^1 \to \lvert \mathrm B \CG \rvert$, and the complete answer is given in Remark \ref{rem_TQFT_conjugation_invariance}. Up to homotopy, $f$ should maps the 1-cell $S^1 \setminus \{\ast\}$ to a 1-cell in $\lvert \mathrm B \CG \rvert$, which corresponds to an element $g \in \CG$. This element $g$ can be viewed as the holonomy around this circle (depending on the fixed reference point $\ast$). So for every $g \in \CG$ there is a circle in $\mathrm{Cob}_2^\CG$ with holonomy $g$. We denote it by $S^1_g$ and depict it in the following figure:
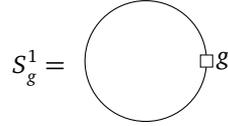
\[
S^1_g =
\begin{array}{c}
\begin{tikzpicture}[scale=0.8]
\draw (0,0) circle (1) ;
\draw[very thin,fill=white] (0.9,-0.1) rectangle (1.1,0.1) node[midway,right] {$g$} ;
\end{tikzpicture}
\end{array}
\]
Here the holonomy is depicted by the defect labeled by $g$, which is the Poincar\'{e} dual of the 1-cell $S^1 \setminus \{\ast\}$.

The circle $S^1_g$ can be decomposed as the composition of three 1-morphisms in $\mathrm{Cob}_2^\CG$ (read from bottom to top):
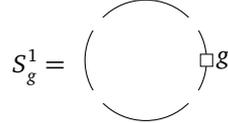
\[
S^1_g =
\begin{array}{c}
\begin{tikzpicture}[scale=0.8]
\draw (0,0) circle (1) ;
\draw[very thin,fill=white] (0.9,-0.1) rectangle (1.1,0.1) node[midway,right] {$g$} ;
\fill[white] (-1,0.5) rectangle (1,0.7) ;
\fill[white] (-1,0.-0.5) rectangle (1,-0.7) ;
\end{tikzpicture}
\end{array}
\]
Since $Z_\CV$ maps a single point to $\CV$, it maps the cup and cap to the coevaluation and evaluation 1-morphism associated to $\CV$, respectively, and the defect $g$ should be mapped to the $\CG$-action on $\CV$ by $g$. Thus $Z_\CV(S^1_g)$ is the composition of three corresponding 1-morphisms in $2\vect$, which is $\chi_\CV(g)$ defined in \eqref{eq_chi_V_g}. So briefly speaking, the 2-character $\chi_\CV$ is the value of the TQFT $Z_\CV$ on the circle. This is also the usual way of defining characters by graph calculus.
\end{rem}

\subsection{The conjugation invariance of 2-characters}

Let $\CG$ be a 2-group and $\CV \in 2\rep(\CG)$. For every $g,x \in G$, there is an isomorphism of vector spaces
\[
\psi_{g,x} \colon \chi_\CV(x) \to \chi_\CV((g \otimes x) \otimes g^*)
\]
defined as follows:
\begin{multline*}
\chi_\CV(x) = \nat(1_\CV,x \odot -) \xrightarrow{(g \odot -) \circ - \circ (g^* \odot -)} \nat(g \odot (g^* \odot -),g \odot (x \odot (g^* \odot -))) \\
\simeq \nat(1_\CV,((g \otimes x) \otimes g^*) \odot -) = \chi_\CV((g \otimes x) \otimes g^*) .
\end{multline*}
We call these isomorphisms the \emph{conjugation invariance} of the 2-character $\chi_\CV$. This isomorphism $\psi_{g,x}$ can also be depicted as the following pasting diagram (some coherence 2-morphisms are omitted):
\[
\begin{array}{c}
\xymatrix{
\CV \rtwocell^{1_\CV}_{x \odot -}{\alpha} & \CV
}
\end{array}
\mapsto
\begin{array}{c}
\xymatrix@C=3em{
\CV \ar[r]^-{g^* \odot -} & \CV \rtwocell^{1_\CV}_{x \odot -}{\alpha} & \CV \ar[r]^{g \odot -} & \CV
}
\end{array}
\]
Using the definition of 2-characters by ends (or equivalently, using the components of natural transformations), the isomorphism $\psi_{g,x}$ is the unique morphism rendering the following diagram commutative for all $v \in \CV$:
\be \label{diag_character_conjugation_invariance}
\begin{array}{c}
\xymatrix{
\chi_\CV(x) \ar@{-->}[rr]^-{\psi_{g,x}} \ar[d]_{\tau_{g^* \odot v}} & & \chi_\CV((g \otimes x) \otimes g^*) \ar[d]^{\tau_v} \\
\CV(g^* \odot v,x \odot (g^* \odot v)) \ar[r]^-{\simeq} & \CV(v,g \odot (x \odot (g^* \odot v)) \ar[r]^-{\simeq} & \CV(v,((g \otimes x) \otimes g^*) \odot v)
}
\end{array}
\ee
It is clear that $\psi_{g,x}$ is natural in $x$, i.e., the following diagram commutes for every morphism $a \colon x \to y$ in $\CG$:
\[
\xymatrix@C=5em{
\chi_\CV(x) \ar[r]^-{\psi_{g,x}} \ar[d]_{\chi_\CV(a)} & \chi_\CV((g \otimes x) \otimes g^*) \ar[d]^{\chi_{\CV((1 \otimes a) \otimes 1)}} \\
\chi_\CV(y) \ar[r]^-{\psi_{g,y}} & \chi_\CV((g \otimes y) \otimes g^*)
}
\]
Thus $\psi_g \colon \chi_\CV \Rightarrow g^* \odot \chi_\CV \colon \CG \to \vect$ is a natural isomorphism, where the $\CG$-action on $\fun(\CG,\vect)$ is induced by the conjugation $\CG$-action on $\CG$.

Moreover, using the universal property of ends, it is not hard to verify that the following diagrams commute for all $g,h,x \in \CG$ and morphism $a \colon g \to h$ in $\CG$:
\be \label{diag_conjugation_invariance_associativity}
\begin{array}{c}
\xymatrix@C=5em{
\chi_\CV(x) \ar[r]^-{\psi_{g \otimes h,x}} \ar[d]_{\psi_{h,x}} & \chi_\CV(((g \otimes h) \otimes x) \otimes (g \otimes h)^*) \ar[d]^{\simeq} \\
\chi_\CV((h \otimes x) \otimes h^*) \ar[r]^-{\psi_{g,(h \otimes x) \otimes h^*}} & \chi_\CV((g \otimes ((h \otimes x) \otimes h^*)) \otimes g^*)
}
\end{array}
\quad
\begin{array}{c}
\xymatrix{
\chi_\CV(x) \ar[r]^-{\psi_{g,x}} \ar[dr]_{\psi_{h,x}} & \chi_\CV((g \otimes x) \otimes g^*) \ar[d]^{\chi_\CV((a \otimes 1) \otimes a^*)} \\
 & \chi_\CV((h \otimes x) \otimes h^*)
}
\end{array}
\ee
Hence we see that $\chi_\CV$ together with the natural isomorphisms $\{\psi_g\}_{g \in \CG}$ is an object in the equivariantization $\fun(\CG,\vect)^\CG$. Then by Proposition \ref{prop_duality_homotopy_fixed_point_quotient} we obtain the following result.

\begin{thm} \label{thm_2-character_trace_functor}
Let $\CG$ be a 2-group and $\CV \in 2\rep(\CG)$. Then the 2-character $\chi_\CV \colon \CG \to \vect$ together with the conjugation invariance natural isomorphisms $\{\psi_g\}_{g \in \CG}$ is an object in $\fun(\CG,\vect)^\CG$.
\end{thm}

The category $\fun(\CG,\vect)^\CG \simeq \fun(\CG \dquotient \CG,\vect)$ is called the space of \emph{class functors} on $\CG$. It is the categorification of the space of class functions on a group $G$. Since $\fun(\CG,\vect)^\CG$ is equivalent to $\FZ_1(\vect_\CG)$ via the Frourier 2-transform defined in Section \ref{sec_Fourier_2-transform}, the 2-characters can also be viewed as objects in the Drinfeld center $\FZ_1(\vect_\CG)$.

\begin{rem}
Theorem \ref{thm_2-character_Lagrangian_algebra} was also proved by Ganter and Usher in \cite[Theorem 4.15]{GU16}. By viewing the 2-character $\chi_\CV(g)$ as a ``categorical trace'' in the sense of \cite{GK08} (see Remark \ref{rem_categorical_trace_GK}), the conjugation invariance of $\chi_\CV$ is a special case of the conjugation invariance of the categorical trace as developed in \cite[Proposition 3.8]{GK08}. The first diagram of \eqref{diag_conjugation_invariance_associativity} is the same as \cite[Proposition 3.8 (b)]{GK08}, in which the associators are omitted. When $\CG = G$ is an ordinary finite group, Bartlett \cite{Bar09,Bar09a} also proved that the 2-characters $\chi_\CV$ are objects in $\FZ_1(\vect_G)$. In \cite[Corollary 9.11]{Bar09a}, he concluded that the 2-characters, as objects in $\FZ_1(\vect_G)$, are not enough to distinguish the 2-representations. In the next subsection, we consider more structures on the 2-characters so that the 2-representations can be distinguished by their 2-characters (see Corollary \ref{cor_2-character_correspond_2-representation}).
\end{rem}

\begin{rem} \label{rem_conjugation_invariance_commutator_functor_shadow}
The conjugation invariance of a 2-character $\chi_\CV$ is equivalent to a family of isomorphisms
\[
\theta_{g,h} \colon \chi_\CV(g \otimes h) \to \chi_\CV(h \otimes g) , \quad g,h \in \CG
\]
satisfying certain coherence conditions. More precisely, if we take
\[
\theta_{g,h} \coloneqq \bigl( \chi_\CV(g \otimes h) \simeq \chi_\CV((g \otimes (h \otimes g)) \otimes g^*) \xrightarrow{\psi_{g,h \otimes g}^{-1}} \chi_\CV(h \otimes g) \bigr) ,
\]
then the first diagram of \eqref{diag_conjugation_invariance_associativity} is equivalent to say that $(\chi_\CV,\theta)$ is a ``commutator functor'' in the sense of \cite{BFO09,Ost14} or a ``categorified trace'' in the sense of \cite{HPT16}. By \cite[Lemma]{HPT16} (see also \cite[Remark 2.2]{Ost14}), $\theta_{g,\one} = \psi_{g,g}^{-1}$ is identity. This is related to the $S^1$-invariance of $\chi_\CV$ (see Theorem \ref{thm_2-character_S1_invariance} below). Moreover, it is obvious that $\theta_{\one,g} = \psi_{\one,g}^{-1}$ is also identity. This condition implies that $\chi_\CV$ is a ``shadow'' in the sense of \cite{Pon18,PS13} (see \cite[Remark 3.2]{HPT16} for more explanation of the difference between ``categorified trace'' and ``shadow''). More precisely, if we take
\[
\theta'_{g,h} \coloneqq \bigl( \chi_\CV(g \otimes h) \xrightarrow{\psi_{g^*,g \otimes h}} \chi_\CV((g^* \otimes (g \otimes h)) \otimes g) \simeq \chi_\CV(h \otimes g) \bigr) ,
\]
then $(\chi_\CV,\theta')$ is a ``shadow'' on the 2-category $\mathrm B \CG$ as defined in of \cite{Pon18,PS13}.

Viewing $\chi_\CV$ as a functor from $\vect_\CG$ to $\vect$, by \cite[Proposition 2.5]{Ost14} (see also \cite[Section 6]{BFO09}) its left adjoint $\chi_\CV^L \colon \vect \to \vect_\CG$ factors through $\FZ_1(\vect_\CG)$. Then $\chi_\CV^L(\bk) \in \FZ_1(\vect_\CG)$ is exactly the 2-character $\chi_\CV$ as an object in $\FZ_1(\vect_\CG)$. More precisely, $\Phi(\chi_\CV^L(\bk)) \simeq \chi_\CV \in \fun(\CG,\vect)^\CG$, where $\Phi \colon \FZ_1(\vect_\CG) \to \fun(\CG,\vect)^\CG$ is the Fourier 2-transform defined in Section \ref{sec_Fourier_2-transform}.
\end{rem}

\begin{defn}
Let $\CX$ be a finite semisimple category. Its \emph{dimension} is the $\bk$-algebra $\dim(\CX) \coloneqq \nat(1_\CX,1_\CX)$.
\end{defn}

Clearly for any $\CV \in 2\rep(\CG)$ we have $\chi_\CV(e) = \dim(\CV)$. For every $g \in \CG$, there is a $\bk$-algebra automorphism on $\dim(\CV)$:
\be \label{eq_pi_1_action_dim}
\dim(\CV) = \chi_\CV(e) \xrightarrow{\psi_{g,e}} \chi_\CV((g \otimes e) \otimes g^*) \simeq \chi_\CV(e) = \dim(\CV) .
\ee
It is not hard to verify that this gives a $\pi_1(\CG)$-representation on $\dim(\CV)$. Similarly, for any $F \in \fun(\CG,\vect)^\CG$, there is a natural $\pi_1(\CG)$-action on the vector space $F(e)$ (see the discussion before Lemma \ref{lem_invariant_pairing_hom}).

\begin{rem} \label{rem_TQFT_conjugation_invariance}
The conjugation invariance of the 2-character $\chi_\CV$ can also be understood from the TQFT point of view discussed in Remark \ref{rem_G_TQFT}. Recall we have found a 1-morphism $S^1_g$ in $\mathrm{Cob}_2^\CG$ for every $g \in \CG$, which is a circle with holonomy $g$ (and a fixed reference point $\ast$). Now we show that for every $g,x \in \CG$, the 1-morphisms $S^1_{gxg^*}$ and $S^1_x$ are isomorphic in $\mathrm{Cob}_2^\CG$. The following figure depicts a 2-isomorphism $S_x^1 \to S_{gxg^*}^1$ in $\mathrm{Cob}_2^\CG$ (we also draw the direction of the compositions for morphisms and the tensor product for objects):
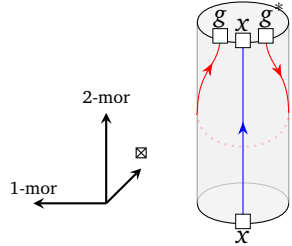
\[
\begin{tikzpicture}[scale=1.2]
\useasboundingbox (0,-0.3) rectangle (1,2.3) ;
%% cylinder
\draw (0,0)--(0,2) ;
\draw (1,0)--(1,2) ;
\fill[gray!10,opacity=0.8] (0,0)--(0,2) .. controls (0,2.3) and (1,2.3) .. (1,2)--(1,0) .. controls (1,0.3) and (0,0.3) .. cycle ;
\draw (0,2) .. controls (0,2.3) and (1,2.3) .. (1,2) ;
\draw (0,0) .. controls (0,0.3) and (1,0.3) .. (1,0) ;
% defects
\draw[thick,red,dotted] (0,1) .. controls (0,0.5) and (1,0.5) .. (1,1) ;
\fill[gray!10,opacity=0.8] (0,0)--(0,2) .. controls (0,1.7) and (1,1.7) .. (1,2)--(1,0) .. controls (1,-0.3) and (0,-0.3) .. cycle ;
\draw (0,2) .. controls (0,1.7) and (1,1.7) .. (1,2) ;
\draw (0,0) .. controls (0,-0.3) and (1,-0.3) .. (1,0) ;
%% defects
\draw[->-,thin,blue] (0.5,-0.2)--(0.5,1.8) ;
\draw[->-=0.65,thin,red] (0,1) .. controls (0,1.5) and (0.25,1.5) .. (0.25,1.85) ;
\draw[->-,thin,red] (0.75,1.85) .. controls (0.75,1.5) and (1,1.5) .. (1,1) ;
\draw[very thin,fill = white] (0.42,-0.28) rectangle (0.58,-0.12) node[midway,below] {$x$} ;
\draw[very thin,fill = white] (0.42,1.72) rectangle (0.58,1.88) node[midway,above] {$x$} ;
\draw[very thin,fill = white] (0.17,1.77) rectangle (0.33,1.93) node[midway,above] {$g$} ;
\draw[very thin,fill = white] (0.67,1.77) rectangle (0.83,1.93) node[midway,above] {$\;\;g^*$} ;
%% frame
\draw[-stealth,thick] (-1,0)--(-1.8,0) node[above] {\scriptsize 1-mor} ;
\draw[-stealth,thick] (-1,0)--(-1,1) node[above] {\scriptsize 2-mor} ;
\draw[-stealth,thick] (-1,0,0)--(-1,0,-1) node[above] {\scriptsize $\boxtimes$} ;
\end{tikzpicture}
\]
Here the $\CG$-structure on the cylinder is depicted by the defects labeled by $g,x$. Hence, a $\CG$-structure on $S^1$, i.e., a continuous map $S^1 \to \lvert \mathrm B \CG \rvert$, up to homotopy, is determined by the conjugacy class of the holonomy $g \in \pi_1(\CG)$. By decomposing this cylinder into small pieces, the image of this cylinder under $Z_\CV \colon \mathrm{Cob}_2^\CG \to 2\vect$ should be
\begin{multline*}
\chi_\CV(x) = \int_{v \in \CV} \CV(v,x \odot v) \simeq \int_{v \in \CV} \CV((g^* \otimes g) \odot v,x \odot v) \simeq \int_{v \in \CV} \CV(g^* \odot (g \odot v),x \odot v) \\
\simeq \int_{w \in \CV} \CV(g^* \odot w,x \odot (g^* \odot w)) \simeq \int_{w \in \CV} \CV(w,((g \otimes x) \otimes g^*) \odot w) = \chi_\CV((g \otimes x) \otimes g^*) .
\end{multline*}
Comparing with the diagram \eqref{diag_character_conjugation_invariance}, this isomorphism is equal to the conjugation invariance $\psi_{g,x}$.

It follows that the isomorphism $\chi_\CV(e) \to \chi_\CV(e)$ defined by \eqref{eq_pi_1_action_dim} is the image of the following 2-morphism in $\mathrm{Cob}_2^\CG$ under $Z_\CV$:
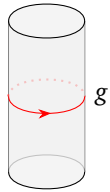
\[
\begin{tikzpicture}
%% cylinder
\draw (0,0)--(0,2) ;
\draw (1,0)--(1,2) ;
\fill[gray!10,opacity=0.8] (0,0)--(0,2) .. controls (0,2.3) and (1,2.3) .. (1,2)--(1,0) .. controls (1,0.3) and (0,0.3) .. cycle ;
\draw (0,2) .. controls (0,2.3) and (1,2.3) .. (1,2) ;
\draw (0,0) .. controls (0,0.3) and (1,0.3) .. (1,0) ;
% defects
\draw[thick,red,dotted] (0,1) .. controls (0,1.3) and (1,1.3) .. (1,1) ;
\fill[gray!10,opacity=0.8] (0,0)--(0,2) .. controls (0,1.7) and (1,1.7) .. (1,2)--(1,0) .. controls (1,-0.3) and (0,-0.3) .. cycle ;
\draw (0,2) .. controls (0,1.7) and (1,1.7) .. (1,2) ;
\draw (0,0) .. controls (0,-0.3) and (1,-0.3) .. (1,0) ;
% defects
\draw[->-,thin,red] (0,1) .. controls (0,0.7) and (1,0.7) .. (1,1) node[right,black] {$g$} ;
\end{tikzpicture}
\]
From this figure we can easily see that \eqref{eq_pi_1_action_dim} defines a $\pi_1(\CG)$-action on the $\bk$-algebra $\chi_\CV(e) = \dim(\CV)$.
\end{rem}

%Now we assume that $\CG = \CG(G,A,\alpha)$ is skeletal.
%
%Let $\CV \in 2\rep(\CG)$. By \cite{GK08}, the 2-character $\chi_\CV$ has the following conjugation invariance: for every $g,x \in G$, there is an isomorphism of vector spaces
%\[
%\psi_g \colon \chi_\CV(x) \to \chi_\CV(gxg^{-1})
%\]
%defined by
%\begin{multline*}
%\chi_\CV(x) = \nat(1_\CV,x \odot -) \xrightarrow{(g \odot -) \circ - \circ (g^{-1} \odot -)} \nat(g \odot (g^{-1} \odot -),g \odot (x \odot (g^{-1} \odot -))) \\
%\simeq \nat(1_\CV,(gxg^{-1}) \odot -) = \chi_\CV(gxg^{-1}) .
%\end{multline*}
%
%For every $g \in G$, the vector space $\chi_\CV(g)$ is equipped with a $\End_\CG(g) = A$-action. It is not hard to see that $\psi_g$ has the following $A$-equivariance condition for every $a \in A = \End_\CG(x)$:
%\[
%\xymatrix{
%\chi_\CV(x) \ar[r]^-{\psi_g} \ar[d]_{\chi_\CV(a)} & \chi_\CV(gxg^{-1}) \ar[d]^{\chi_\CV(g \triangleright a)} \\
%\chi_\CV(x) \ar[r]^-{\psi_g} & \chi_\CV(gxg^{-1})
%}
%\]
%
%$\psi_g \circ \psi_h = ? \circ \psi_{gh}$

\subsection{The \texorpdfstring{$S^1$}{S1}-invariance of 2-characters}

Recall that there is a $\mathrm B \Zb$-action on $\CG \dquotient \CG$ defined in Remark \ref{rem_S1_action_LG}. Then it induces a $\mathrm B \Zb$-action on the space of class functors $\fun(\CG,\vect)^\CG \simeq \fun(\CG \dquotient \CG,\vect)$, given by a natural automorphism $\phi \colon 1_{\fun(\CG,\vect)^\CG} \Rightarrow 1_{\fun(\CG,\vect)^\CG}$ defined by
\[
(\phi_{(F,\psi)})_g \coloneqq \bigl(F(g) \xrightarrow{\psi_{g,g}} F((g \otimes g) \otimes g^*) \simeq F(g) \bigr) , \quad (F,\psi) \in \fun(\CG,\vect)^\CG , \, g \in \CG .
\]
For a class functor, being a $\mathrm B \Zb$-fixed point is only a property and does not require additional structures, because $\mathrm B \Zb$ has no nontrivial object. In other words, the forgetful functor $(\fun(\CG,\vect)^\CG)^{\mathrm B \Zb} \to \fun(\CG,\vect)^\CG$ is fully faithful.

\begin{thm} \label{thm_2-character_S1_invariance}
For any $\CV \in 2\rep(\CG)$, the 2-character $\chi_\CV \in \fun(\CG,\vect)^\CG$ is a $\mathrm B \Zb$-fixed point.
\end{thm}

\pf
We need to show that the conjugation invariance
\[
\psi_{g,g} \colon \chi_\CV(g) \to \chi_\CV((g \otimes g) \otimes g^*) \simeq \chi_\CV(g)
\]
is identity for every $g \in \CG$. As discussed in Remark \ref{rem_conjugation_invariance_commutator_functor_shadow}, this is a corollary of \cite[Lemma 3.3]{HPT16} (see also \cite[Remark 2.2]{Ost14}). Here we give a direct proof. By the diagram \eqref{diag_character_conjugation_invariance}, it is equivalent to show that the following diagram commutes for any $v \in \CV$ and any natural transformation $\alpha \in \nat(1_\CV,g \odot -) \simeq \chi_\CV(g)$:
\[
\xymatrix@C=4em{
(g \otimes g^*) \odot v \ar[r]^-{\simeq} \ar[d]_{\simeq} & g \odot (g^* \odot v) \ar[r]^-{1 \odot \alpha_{g^* \odot v}} & g \odot (g \odot (g^* \odot v)) \ar[d]^{\simeq} \\
v \ar[r]^-{\alpha_v} & g \odot v & ((g \otimes g) \otimes g^*) \odot v \ar[l]_-{\simeq}
}
\]
This is the outer diagram of the following diagram, which is commutative by the naturality of $\alpha$:
\[
\xymatrix@C=4em{
g^* \odot v \ar[r]^-{\alpha_{g^* \odot v}} \ar[d]_{\alpha_{g^* \odot v}} & g \odot (g^* \odot v) \ar[r]^-{\simeq} \ar[d]^{\alpha_{g \odot (g^* \odot v)}} & v \ar[d]^{\alpha_v} \\
g \odot (g^* \odot v) \ar[r]^-{1 \odot \alpha_{g^* \odot v}} & g \odot (g \odot (g^* \odot v)) \ar[r]^-{\simeq} & g \odot v
}
\]
Here to show the commutativity of the right rectangle, recall that the bottom right isomorphism is induced from $g \otimes (g \otimes g^*) \xrightarrow{1 \otimes \ev} g$.
\epf

\begin{rem}
The 2-group $\mathrm B \Zb$ can be viewed as an algebraic model of the topological group $S^1$. From the TQFT point of view discussed in Remark \ref{rem_G_TQFT} and \ref{rem_TQFT_conjugation_invariance}, the $\mathrm B \Zb$-invariance means that the image of the following cylinder
\[
\begin{array}{c}
\begin{tikzpicture}[scale=1.2]
\useasboundingbox (0,-0.3) rectangle (1,2.3) ;
%% cylinder
\draw (0,0)--(0,2) ;
\draw (1,0)--(1,2) ;
\fill[gray!10,opacity=0.8] (0,0)--(0,2) .. controls (0,2.3) and (1,2.3) .. (1,2)--(1,0) .. controls (1,0.3) and (0,0.3) .. cycle ;
\draw (0,2) .. controls (0,2.3) and (1,2.3) .. (1,2) ;
\draw (0,0) .. controls (0,0.3) and (1,0.3) .. (1,0) ;
% defects
\draw[thick,red,dotted] (0,1) .. controls (0,0.5) and (1,0.5) .. (1,1) ;
\fill[gray!10,opacity=0.8] (0,0)--(0,2) .. controls (0,1.7) and (1,1.7) .. (1,2)--(1,0) .. controls (1,-0.3) and (0,-0.3) .. cycle ;
\draw (0,2) .. controls (0,1.7) and (1,1.7) .. (1,2) ;
\draw (0,0) .. controls (0,-0.3) and (1,-0.3) .. (1,0) ;
%% defects
\draw[->-,thin,red] (0.5,-0.2)--(0.5,1.8) ;
\draw[->-=0.65,thin,red] (0,1) .. controls (0,1.5) and (0.25,1.5) .. (0.25,1.85) ;
\draw[->-,thin,red] (0.75,1.85) .. controls (0.75,1.5) and (1,1.5) .. (1,1) ;
\draw[very thin,fill = white] (0.42,-0.28) rectangle (0.58,-0.12) node[midway,below] {$g$} ;
\draw[very thin,fill = white] (0.42,1.72) rectangle (0.58,1.88) node[midway,above] {$g$} ;
\draw[very thin,fill = white] (0.17,1.77) rectangle (0.33,1.93) node[midway,above] {$g$} ;
\draw[very thin,fill = white] (0.67,1.77) rectangle (0.83,1.93) node[midway,above] {$\;\;g^*$} ;
\end{tikzpicture}
\end{array}
=
\begin{array}{c}
\begin{tikzpicture}[scale=1.2]
\useasboundingbox (0,-0.3) rectangle (1,2.3) ;
%% cylinder
\draw (0,0)--(0,2) ;
\draw (1,0)--(1,2) ;
\fill[gray!10,opacity=0.8] (0,0)--(0,2) .. controls (0,2.3) and (1,2.3) .. (1,2)--(1,0) .. controls (1,0.3) and (0,0.3) .. cycle ;
\draw (0,2) .. controls (0,2.3) and (1,2.3) .. (1,2) ;
\draw (0,0) .. controls (0,0.3) and (1,0.3) .. (1,0) ;
% defects
\draw[thick,red,dotted] (1,0.4) .. controls (1,0.8) and (0,0.8) .. (0,1.2) ;
\fill[gray!10,opacity=0.8] (0,0)--(0,2) .. controls (0,1.7) and (1,1.7) .. (1,2)--(1,0) .. controls (1,-0.3) and (0,-0.3) .. cycle ;
\draw (0,2) .. controls (0,1.7) and (1,1.7) .. (1,2) ;
\draw (0,0) .. controls (0,-0.3) and (1,-0.3) .. (1,0) ;
%% defects
\draw[->-,thin,red] (0.5,-0.2) .. controls (0.5,0.2) and (1,0.1) .. (1,0.4) ;
\draw[->-,thin,red] (0,1.2) .. controls (0,1.5) and (0.5,1.4) .. (0.5,1.8) ;
%\draw[->-=0.65,thin,red] (0,1) .. controls (0,1.5) and (0.25,1.5) .. (0.25,1.85) ;
%\draw[->-,thin,red] (0.75,1.85) .. controls (0.75,1.5) and (1,1.5) .. (1,1) ;
\draw[very thin,fill = white] (0.42,-0.28) rectangle (0.58,-0.12) node[midway,below] {$g$} ;
\draw[very thin,fill = white] (0.42,1.72) rectangle (0.58,1.88) node[midway,above] {$g$} ;
\end{tikzpicture}
\end{array}
\]
is identity:
\[
Z_\CV \biggl(
\begin{array}{c}
\begin{tikzpicture}[scale=1.2]
\useasboundingbox (0,-0.3) rectangle (1,2.3) ;
%% cylinder
\draw (0,0)--(0,2) ;
\draw (1,0)--(1,2) ;
\fill[gray!10,opacity=0.8] (0,0)--(0,2) .. controls (0,2.3) and (1,2.3) .. (1,2)--(1,0) .. controls (1,0.3) and (0,0.3) .. cycle ;
\draw (0,2) .. controls (0,2.3) and (1,2.3) .. (1,2) ;
\draw (0,0) .. controls (0,0.3) and (1,0.3) .. (1,0) ;
% defects
\draw[thick,red,dotted] (1,0.4) .. controls (1,0.8) and (0,0.8) .. (0,1.2) ;
\fill[gray!10,opacity=0.8] (0,0)--(0,2) .. controls (0,1.7) and (1,1.7) .. (1,2)--(1,0) .. controls (1,-0.3) and (0,-0.3) .. cycle ;
\draw (0,2) .. controls (0,1.7) and (1,1.7) .. (1,2) ;
\draw (0,0) .. controls (0,-0.3) and (1,-0.3) .. (1,0) ;
%% defects
\draw[->-,thin,red] (0.5,-0.2) .. controls (0.5,0.2) and (1,0.1) .. (1,0.4) ;
\draw[->-,thin,red] (0,1.2) .. controls (0,1.5) and (0.5,1.4) .. (0.5,1.8) ;
%\draw[->-=0.65,thin,red] (0,1) .. controls (0,1.5) and (0.25,1.5) .. (0.25,1.85) ;
%\draw[->-,thin,red] (0.75,1.85) .. controls (0.75,1.5) and (1,1.5) .. (1,1) ;
\draw[very thin,fill = white] (0.42,-0.28) rectangle (0.58,-0.12) node[midway,below] {$g$} ;
\draw[very thin,fill = white] (0.42,1.72) rectangle (0.58,1.88) node[midway,above] {$g$} ;
\end{tikzpicture}
\end{array}
\biggr) = Z_\CV \biggl(
\begin{array}{c}
\begin{tikzpicture}[scale=1.2]
\useasboundingbox (0,-0.3) rectangle (1,2.3) ;
%% cylinder
\draw (0,0)--(0,2) ;
\draw (1,0)--(1,2) ;
\fill[gray!10,opacity=0.8] (0,0)--(0,2) .. controls (0,2.3) and (1,2.3) .. (1,2)--(1,0) .. controls (1,0.3) and (0,0.3) .. cycle ;
\draw (0,2) .. controls (0,2.3) and (1,2.3) .. (1,2) ;
\draw (0,0) .. controls (0,0.3) and (1,0.3) .. (1,0) ;
\fill[gray!10,opacity=0.8] (0,0)--(0,2) .. controls (0,1.7) and (1,1.7) .. (1,2)--(1,0) .. controls (1,-0.3) and (0,-0.3) .. cycle ;
\draw (0,2) .. controls (0,1.7) and (1,1.7) .. (1,2) ;
\draw (0,0) .. controls (0,-0.3) and (1,-0.3) .. (1,0) ;
%% defects
\draw[->-,thin,red] (0.5,-0.2)--(0.5,1.8) ;
\draw[very thin,fill = white] (0.42,-0.28) rectangle (0.58,-0.12) node[midway,below] {$g$} ;
\draw[very thin,fill = white] (0.42,1.72) rectangle (0.58,1.88) node[midway,above] {$g$} ;
\end{tikzpicture}
\end{array}
\biggr) = 1_{\chi_\CV(g)} .
\]
In other words, the 2-character $\chi_\CV$, as the image of the circle $S^1$ under $Z_\CV$, is invariant under the rotation of $S^1$ (i.e., the Dehn twist).
\end{rem}

\begin{rem}
The $\mathrm{B}\mathbb{Z}$-invariance of the 2-character $\chi_\mathcal{V}$, defined as the value of the $G$-TQFT $Z_\mathcal{V}$ on the circle $S^1$, fits into a broader pattern connecting categorical traces, Hochschild homology, and loop spaces, which we briefly explain.

For a ring spectrum $R$, there is a \emph{Dennis trace} map $$K(R) \to \mathrm{THH}(R)$$ from its \emph{algebraic $K$-theory} to \emph{topological Hochschild homology} \cite{BHM93,DGM12}, which can be regarded as a categorified Chern character. When $R = \bk[G]$ is the ordinary finite group algebra, we have $K_0(R)$ being the representation ring of $G$ and its image under the Dennis trace is precisely the associated characters valued in $\mathrm{THH}_0(R)$ being the center of $R$. 

The target $\mathrm{THH}(R)$ is the geometric realization of the cyclic bar construction on $R$ and carries a canonical $S^1$-action arising from rotating the circle. The homotopy fixed points $\mathrm{THH}(R)^{hS^1}$ is denoted by the \emph{topological negative cyclic homology} $\mathrm{TC}^-(R)$; after adding more refined cyclotomic fixed points yield  and topological cyclic homology $\mathrm{TC}(R)$, respectively. Dennis trace map factors through homotopy $S^1$ fixed points, providing the \emph{cyclotomic trace}, the primary computational tool for chromatic homotopy theory \cite{BHM93, NS18}: \[K(R) \xrightarrow{\text{cyclotomic character}} \mathrm{TC}(R) \xrightarrow{\text{cyclotomic information}} \mathrm{TC}^-(R) \xrightarrow{S^1\text{ fixed points}} \mathrm{THH}(R).\]

In derived algebraic geometry, Ben-Zvi, Francis, and Nadler \cite{BZFN10} proved that for a perfect stack $X$, the Drinfeld center $\FZ_1(\mathrm{QC}(X))$, the categorical trace (Hochschild homology category) of $\mathrm{QC}(X)$, and the category $\mathrm{QC}(\mathrm{L}X)$ of quasi-coherent sheaves on the free loop space $\mathrm{L}X$ are all equivalent. Our Corollary \ref{cor_fun_conjugation_Drinfeld_center} — the equivalence $\fun(G \dquotient G, \vect) \simeq \FZ_1(\vect_G)$ — is a finite, discrete instance of this general phenomenon, where $X = \mathrm{B}G$ and $\mathrm{L}X = \mathrm{LB}G = G \dquotient G$ (with the conjugation action). In this context, the 2-character $\chi_\CV$ is the categorical analogue of the Chern character: the Dennis trace $K(R) \to \mathrm{THH}(R)$ is categorified to the map from $2\rep(G)$ (playing the role of $K$-theory) to $\fun(G \dquotient G, \vect) \simeq \FZ_1(\vect_G)$ (playing the role of Hochschild homology). In the language of Hoyois, Scherotzke, and Sibilla \cite{HSS17}, the categorified Chern character is a symmetric monoidal functor to $S^1$-equivariant complexes on the derived free loop stack.

The $S^1$-action on $\mathrm{THH}(R)$, which ultimately comes from the rotational symmetry of the cyclic bar construction, has a categorical analogue in our setting. The conjugation invariance data $(\psi_g)_{g \in G}$ of the 2-character (Theorem \ref{thm_2-character_trace_functor}) can be understood as the finite, discrete shadow of this $S^1$-equivariance: the mapping space $\map(S^1,\lvert \mathrm{B}G \rvert) \simeq \lvert \mathrm{LB}G \rvert \simeq \lvert G \dquotient G \rvert$ inherits an $S^1$-action from reparametrization of the circle. Theorem \ref{thm_2-character_S1_invariance} further shows that 2-character $\chi_\CV$, as a functor on $G \dquotient G$, is invariant with respect to this action, as a reminiscent of the $S^1$-invariance of the Dennis trace.

See also Ben-Zvi and Nadler \cite{BZN21, BZN13} for the theory of nonlinear and secondary traces in derived algebraic geometry, and Carmeli, Cnossen, Ramzi, and Yanovski \cite{CCRY25} for a recent $\infty$-categorical character theory via categorified traces that provides the general framework encompassing both the Dennis trace and our construction in principle.
\end{rem}

\subsection{The irreducible 2-characters are Lagrangian algebras}

\begin{defn}
Let $\CC$ be a braided multi-fusion category. Then it is the direct sum of braided fusion categories $\CC = \bigoplus_{i=1}^n \CC_i$. We say $\CC$ is \emph{nondegenerate} if each $\CC_i$ is. A \emph{Lagrangian algebra} in $\CC$ is a Lagrangian algebra in a direct summand $\CC_i$.
\end{defn}

Let $\CG$ be a finite 2-group and $\CV$ be an indecomposable (i.e., irreducible) finite semisimple 2-representation of $\CG$. By Remark \ref{rem_full_center_Lagrangian}, the full center $Z(\CV)$ of $\CV$ is a Lagrangian algebra in the braided multi-fusion category $\FZ_1(\vect_\CG)$. Recall Example \ref{expl_full_center_end} that $Z(\CV)$ is the end $\int_{v \in \CV} [v,v] \in \vect_\CG$ equipped with a canonical half-braiding.

\begin{thm} \label{thm_2-character_Lagrangian_algebra}
Let $\Phi \colon \FZ_1(\vect_\CG) \to \fun(\CG,\vect)^\CG$ be the Fourier 2-transform defined in Section \ref{sec_Fourier_2-transform}. There is a canonical isomorphism $\Phi(Z(\CV)) \simeq \chi_{\CV^\op}$ in $\fun(\CG,\vect)^\CG$. Hence the irreducible 2-character $\chi_\CV$ admits a structure of a Lagrangian algebra in $\fun(\CG,\vect)^\CG$, inherited from the full center $Z(\CV^\op)$.
\end{thm}

\pf
First we construct an isomorphism $\sigma \colon \Phi(Z(\CV)) \simeq \chi_{\CV^\op}$ in $\fun(\CG,\vect)$. For every $x \in \CG$, we define $\sigma_x$ to be the following composite isomorphism:
\begin{multline*}
\Phi(Z(\CV))(x) = \vect_\CG(I(x),Z(\CV)) \simeq \vect_\CG\biggl( I(x),\int_{v \in \CV} [v,v] \biggr) \\
\simeq \int_{v \in \CV} \vect_\CG(I(x),[v,v]) \simeq \int_{v \in \CV} \CV(I(x) \odot v,v) = \int_{v \in \CV^\op} \CV^\op(v,x \odot v) = \chi_{\CV^\op}(x) .
\end{multline*}
The naturality of $\sigma$ is clear. Then we need to show that $\sigma$ is a morphism in $\fun(\CG,\vect)^\CG$, i.e., the following diagram commutes for every $g \in \CG$:
\[
\xymatrix@C=4em{
\Phi(Z(\CV)) \ar[rr]^-{\sigma} \ar[d] & & \chi_{\CV^\op} \ar[d] \\
\Phi(g^* \odot Z(\CV)) \ar[r]^-{\simeq} & g^* \odot \Phi(Z(\CV)) \ar[r]^-{1 \odot \sigma} & g^* \odot \chi_{\CV^\op}
}
\]
where the vertical arrows are induced by the $\CG$-equivariant structures on $Z(\CV)$ and $\chi_{\CV^\op}$, respectively (see the diagrams \eqref{diag_full_center_half_braiding} and \eqref{diag_character_conjugation_invariance}). In other words, we need to show the commutativity of the following diagram for all $g,x \in \CG$:
\[
\xymatrix@C=5em{
\Phi(Z(\CV))(x) \ar[rr]^-{\sigma_x} \ar[d] & & \chi_{\CV^\op}(x) \ar[d]^{\psi_{g,x}} \\
\Phi((g^* \otimes Z(\CV)) \otimes g)(x) \ar[r]^-{\simeq} & \Phi(Z(\CV))((g \otimes x) \otimes g^*) \ar[r]^-{\sigma_{(g \otimes x) \otimes g^*}} & \chi_{\CV^\op}((g \otimes x) \otimes g^*)
}
\]
By the universal property of ends and expanding the vertical arrows, it suffices to show that for any $g,x \in \CG$ and $v \in \CV$ the following diagram commutes:
\[
\xymatrix{
\vect_\CG(x,[g^* \odot v,g^* \odot v]) \ar[rr] \ar[d] & & \CV(x \odot (g^* \odot v),g^* \odot v) \ar[d] \\
\vect_\CG(x,(g^* \otimes [v,v]) \otimes g) \ar[r] & \vect_\CG((g \otimes x) \otimes g^*,[v,v]) \ar[r] & \CV(((g \otimes x) \otimes g^*) \odot v,v)
}
\]
The commutativity follows from the definition of the canonical isomorphism $[g^* \odot v,g^* \odot v] \simeq g^* \otimes [v,v] \otimes g$. Hence we prove the isomorphism $\Phi(Z(\CV)) \simeq \chi_{\CV^\op}$.
\epf

By \cite[Proposition 4.8]{DMNO13}, for a fusion category $\CC$, taking the full center defines a one-to-one correspondence between equivalence classes of indecomposable finite semisimple $\CC$-modules and isomorphism classes of Lagrangian algebras in $\FZ_1(\CC)$. The same statement also holds when $\CC$ is multi-fusion. Hence we obtain the following result.

\begin{cor} \label{cor_2-character_correspond_2-representation}
For a finite 2-group $\CG$, its simple 2-representations are determined by their 2-characters. More precisely, two simple 2-representations $\CV,\CW \in 2\rep(\CG)$ are equivalent if and only if their 2-characters $\chi_\CV,\chi_\CW \in \fun(\CG,\vect)^\CG$ are isomorphic as Lagrangian algebras.
\end{cor}

\begin{rem}
Moreover, by \cite{DKR11a}, this one-to-one correspondence can be lifted to an equivalence between the groupoid of finite semisimple 2-representations and equivalences between 2-representations, and the groupoid of Lagrangian algebras (i.e., 2-characters) in $\fun(\CG,\vect)^\CG$ and algebra isomorphisms (i.e., $\CG$-equivariant monoidal natural isomorphisms).
\end{rem}

\begin{rem}
For a finite group $G$, the fusion rule of $\rep(G)$ can be calculated from the pointwise product of irreducible characters as class functions. Similarly, for a finite 2-group $\CG$, the fusion rule of $2\rep(\CG)$ can be calculated from the pointwise product of Lagrangian algebras in $\FZ_1(2\vect_\CG)$ as class functors.
\end{rem}

\begin{rem}
It may seems un-natural that $\Phi(Z(\CV)) \simeq \chi_{\CV^\op}$, not $\chi_\CV$. Indeed, in the de-categorified case (see Remark \ref{rem_Fourier_transform}), for a finite group $G$ and an irreducible $G$-representation $V$, there is a central idempotent
\be \label{eq_central_idempotent}
p_V \coloneqq \frac{\dim(V)}{\lvert G \rvert} \sum_{g \in G} \chi_V(g^{-1}) \cdot g \in Z(\bk[G]) .
\ee
It is known that irreducible $G$-representations are one-to-one correspond to primitive (minimal) central idempotent in $\bk[G]$. More precisely, the central idempotent corresponding to $V$ is $p_V$, because $p_V$ acts on $V$ as identity and acts on other irreducible $G$-representations as zero. Then it is easy to see that $\Phi(p_V) = \chi_{V^*}$, not $\chi_V$.
\end{rem}

It is known that the algebras in the functor category $\fun(\CG,\vect)$ with the tensor product being the Day convolution are lax monoidal functors \cite{Day70}. Let us compute the multiplication of the algebra $\chi_{\CV^\op}$, i.e., its lax monoidal structure. The multiplication $\Phi(Z(\CV)) \circledast \Phi(Z(\CV)) \Rightarrow \Phi(Z(\CV))$ is given by
\[
(\Phi(Z(\CV)) \circledast \Phi(Z(\CV)))(x) \simeq \Phi(Z(\CV) \otimes Z(\CV))(x) \xrightarrow{\Phi(\mu)(x)} \Phi(Z(\CV))(x) ,
\]
where the first equivalence is the monoidal structure of $\Phi$ induced from \eqref{diag_monoidal_structure_Fourier_2-transform}, and the multiplication $\mu \colon Z(\CV) \otimes Z(\CV) \to Z(\CV)$ is given by \eqref{diag_full_center_multiplication}. Thus for $h,k \in \CG$, the lax monoidal structure
\[
\Phi(Z(\CV))(h) \otimes \Phi(Z(\CV))(k) \to \Phi(Z(\CV))(h \otimes k)
\]
is given by
\begin{multline*}
\Phi(Z(\CV))(h) \otimes \Phi(Z(\CV))(k) = \vect_\CG(I(h),Z(\CV)) \otimes \vect_\CG(I(k),Z(\CV)) \\
\xrightarrow{\otimes} \vect_\CG(I(h \otimes k),Z(\CV) \otimes Z(\CV)) \xrightarrow{\vect_\CG(1,\mu)} \vect_\CG(I(h \otimes k),Z(\CV)) = \Phi(Z(\CV))(h \otimes k) .
\end{multline*}
Using the definition of $Z(\CV)$ by ends, this morphism is the unique morphism rendering the following diagram commutative for all $u \in \CV$:
\be \label{diag_lax_monoidal_chi_V}
\begin{array}{c}
\xymatrix@C=2em{
\chi_{\CV^\op}(h) \otimes \chi_{\CV^\op}(k) \ar@{-->}[rr] \ar@{=}[d] & & \chi_{\CV^\op}(h \otimes k) \ar@{=}[d] \\
\int_{v \in \CV} \vect_\CG(h,[v,v]) \otimes \vect_\CG(k,[w,w]) \ar[d]_{\tau_u \otimes \tau_u} & & \int_{v \in \CV} \vect_\CG(h \otimes k,[v,v]) \ar[d]^{\tau_u} \\
\vect_\CG(h,[u,u]) \otimes \vect_\CG(k,[u,u]) \ar[r]^-{\otimes} & \vect_\CG(h \otimes k,[u,u] \otimes [u,u]) \ar[r] & \vect_\CG(h \otimes k,[u,u])
}
\end{array}
\ee

\begin{rem}
We also obtain another construction of the full center $Z(\CV)$. Its underlying object is
\[
Z(\CV) \simeq \Psi(\chi_{\CV^\op}) = \int^{g \in \CG} \chi_\CV(g^*) \odot I(g) \in \vect_\CG ,
\]
and the half-braiding is induced from the conjugation invariance (see \eqref{eq_half_braiding_Psi}), and the multiplication is induced from the lax monoidal structure of $\chi_{\CV^\op}$. This categorifies the construction of the central idempotent \eqref{eq_central_idempotent}.
\end{rem}

\begin{rem}
The commutative (Frobenius) algebra structure of the 2-character $\chi_\CV$ can also be understood from the TQFT point of view discussed in Remark \ref{rem_G_TQFT}. For every $h,k \in \CG$, consider the pants depicted in the following figure, which is a 2-morphism in $\mathrm{Cob}_2^\CG$:
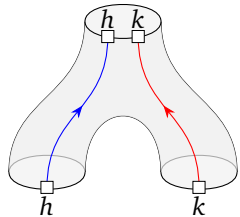
\[
\begin{tikzpicture}
\draw (-1,0) .. controls (-1,1) and (0,1) .. (0,2) ;
\draw (2,0) .. controls (2,1) and (1,1) .. (1,2) ;
\draw (0,0) .. controls (0,1) and (1,1) .. (1,0) ;
\fill[gray!10,opacity=0.7] (-1,0) .. controls (-1,1) and (0,1) .. (0,2) .. controls (0,2.3) and (1,2.3) .. (1,2) .. controls (1,1) and (2,1) .. (2,0) .. controls (2,0.3) and (1,0.3) .. (1,0) .. controls (1,1) and (0,1) .. (0,0) .. controls (0,0.3) and (-1,0.3) .. (-1,0) ;
\draw (0,2) .. controls (0,2.3) and (1,2.3) .. (1,2) ;
\draw (-1,0) .. controls (-1,0.3) and (0,0.3) .. (0,0) ;
\draw (1,0) .. controls (1,0.3) and (2,0.3) .. (2,0) ;
\fill[gray!10,opacity=0.7] (-1,0) .. controls (-1,1) and (0,1) .. (0,2) .. controls (0,1.7) and (1,1.7) .. (1,2) .. controls (1,1) and (2,1) .. (2,0) .. controls (2,-0.3) and (1,-0.3) .. (1,0) .. controls (1,1) and (0,1) .. (0,0) .. controls (0,-0.3) and (-1,-0.3) .. (-1,0) ;
\draw (0,2) .. controls (0,1.7) and (1,1.7) .. (1,2) ;
\draw (-1,0) .. controls (-1,-0.3) and (0,-0.3) .. (0,0) ;
\draw (1,0) .. controls (1,-0.3) and (2,-0.3) .. (2,0) ;
\draw[->-,thin,blue] (-0.5,-0.2) .. controls (-0.5,0.8) and (0.3,0.8) .. (0.3,1.8) ;
\draw[->-,thin,red] (1.5,-0.2) .. controls (1.5,0.8) and (0.7,0.8) .. (0.7,1.8) ;
\draw[very thin,fill=white] (-0.58,-0.28) rectangle (-0.42,-0.12) node[midway,below] {$h$} ;
\draw[very thin,fill=white] (0.22,1.72) rectangle (0.38,1.88) node[midway,above] {$h$} ;
\draw[very thin,fill=white] (1.42,-0.28) rectangle (1.58,-0.12) node[midway,below] {$k$} ;
\draw[very thin,fill=white] (0.62,1.72) rectangle (0.78,1.88) node[midway,above] {$k$} ;
\end{tikzpicture}
\]
By decomposing this pants into small pieces, the most nontrivial one is the following cobordism:
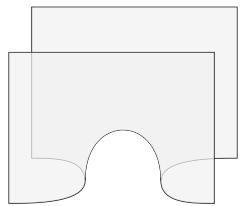
\[
\begin{tikzpicture}
\draw[very thin] (1,0.3) .. controls (1,1.2) and (2,1.2) .. (2,0.3) ;
\draw (2,0.3) .. controls (2,0.6) and (2.7,0.6) .. (3,0.6)--(3,2.6)--(0.3,2.6)--(0.3,0.6) .. controls (0.6,0.6) and (1,0.6) .. (1,0.3) ;
\fill[gray!10,opacity=0.7] (0.3,0.6) .. controls (0.6,0.6) and (1,0.6) .. (1,0.3) .. controls (1,1.2) and (2,1.2) .. (2,0.3) .. controls (2,0.6) and (2.7,0.6) .. (3,0.6)--(3,2.6)--(0.3,2.6)--cycle ;
\draw (2,0.3) .. controls (2,0) and (2.4,0) .. (2.7,0)--(2.7,2)--(0,2)--(0,0) .. controls (0.3,0) and (1,0) .. (1,0.3) ;
\fill[gray!10,opacity=0.7] (0,0) .. controls (0.3,0) and (1,0) .. (1,0.3) .. controls (1,1.2) and (2,1.2) .. (2,0.3) .. controls (2,0) and (2.4,0) .. (2.7,0)--(2.7,2)--(0,2)--cycle ;
\end{tikzpicture}
\]
Its image under $Z_\CV$ is the counit $\varepsilon$ for the adjunction between the coevaluation 1-morphism $\coev \colon \vect \to \CV^\op \boxtimes \CV$ and its right adjoint $\Hom_\CV \colon \CV^\op \boxtimes \CV \to \vect$. So for every $v,w \in \CV$, the morphism
\[
\varepsilon_{v \boxtimes w} \colon \CV(v,w) \odot \int_{u \in \CV} u \boxtimes u \to v \boxtimes w
\]
is equal to both two composite morphisms in the following commutative diagram:
\be \label{diag_pants_two_way}
\begin{array}{c}
\xymatrix{
\CV(v,w) \odot \int_{u \in \CV} u \boxtimes u \ar[r]^-{1 \odot \tau_v} \ar[d]_{1 \odot \tau_w} & v \boxtimes (\CV(v,w) \odot v) \ar[d] \\
(\CV(v,w) \odot w) \boxtimes w \ar[r] & v \boxtimes w
}
\end{array}
\ee
So the image of the above pants under $Z_\CV$ is the unique $\bk$-linear map $\chi_\CV(h) \otimes \chi_\CV(k) \to \chi_\CV(h \otimes k)$ rendering the following diagram commutative for all $u \in \CV$:
\[
\xymatrix{
\chi_\CV(h) \otimes \chi_\CV(k) \ar@{-->}[rr] \ar@{=}[d] & & \chi_\CV(h \otimes k) \ar@{=}[d] \\
\int_{v \in \CV} \CV(v,h \odot v) \otimes \int_{w \in \CW} \CV(w,k \odot w) \ar[d]_{\tau_{k \odot u} \otimes \tau_u} & & \int_{v \in \CV} \CV(v,(h \otimes k) \odot v) \ar[d]^{\tau_u} \\
\CV(k \odot u,h \odot (k \odot u)) \otimes \CV(u,k \odot u) \ar[r]^-{\circ} & \CV(u,h \odot (k \odot u)) \ar[r]^-{\simeq} & \CV(u,(h \otimes k) \odot u)
}
\]
By the commutative diagram \eqref{diag_pants_two_way}, this map is also equal to the one rendering the following diagram commutative for all $u \in \CV$:
\be \label{diag_multiplication_2-character}
\xymatrix@C=5em{
\int_{v \in \CV} \CV(v,h \odot v) \otimes \int_{w \in \CW} \CV(w,k \odot w) \ar@{-->}[r] \ar[d]_{\tau_u \otimes \tau_u} & \int_{v \in \CV} \CV(v,(h \otimes k) \odot v) \ar[d]^{\tau_u} \\
\CV(u,h \odot u) \otimes \CV(u,k \odot u) \ar[d]_{1 \otimes (h \odot -)} & \CV(u,(h \otimes k) \odot u) \ar[d]^{\simeq} \\
\CV(u,h \odot u) \otimes \CV(h \odot u,h \odot (k \odot u)) \ar[r]^-{\circ} \ar[r] & \CV(u,h \odot (k \odot u))
}
\ee
It is not hard to see that this is the same as \eqref{diag_lax_monoidal_chi_V} if $\CV$ is replaced by $\CV^\op$. Hence, the Lagrangian algebra structure on the 2-character $\chi_{\CV^\op}$ obtained in Theorem \ref{thm_2-character_Lagrangian_algebra} coincides with the intuition from the TQFT point of view.

As a special case, when $\CG = G$ is an ordinary group, the (non-extended) TQFT $Z_\CV$ should be a homotopy quantum field theory (HQFT) in the sense of Turaev \cite{Tur99,Tur10}. Turaev also showed that a 2D HQFT is determined by its values on the circles with different holonomies, which form a crossed Frobenius $G$-algebra, that is, a commutative Frobenius algebra in $\FZ_1(\vect_G)$ (see \cite{PT08,Por07} for a generalization to HQFTs with the target space modeled by crossed modules). A 3-dimensional generalization can be found in \cite{SV23}. By \cite{Dav10a}, the Lagrangian algebras in $\FZ_1(\vect_G)$ are classified by pairs $(H,[\omega])$ up to conjugation, where $H \subseteq G$ is a subgroup and $[\omega] \in H^2(H;\bk^\times)$ is a 2-cohomology class. The HQFTs corresponding to these Lagrangian algebras are called ``cohomological HQFTs'' by Turaev.
\end{rem}

\begin{rem}
In a 2D open-closed TQFT or CFT, the \emph{open-closed duality} \cite{MS06,RFFS07,KR09} states that the closed string algebra is the full center of the open string algebra. Here the open string algebra and the closed string algebra are the values of the quantum field theory on the open string (interval) $[0,1]$ and the closed string (circle) $S^1$, respectively. The boundary-bulk relation in topological orders \cite{KWZ15,KWZ17} can also be viewed as a variant of the open-closed duality.

Theorem \ref{thm_2-character_Lagrangian_algebra} also reveals the open-closed duality. In this case, the 2-representation $\CV$ is the \emph{category of boundary conditions}. Given a boundary condition $v \in \CV$, the open string algebra is the internal hom $[v,v] \in \vect_\CG$, and the closed string algebra is the full center of $[v,v]$. When $\CV \in 2\rep(\CG)$ is simple, the Morita class of $[v,v]$ is independent of the choice of nonzero $v \in \CV$ \cite{Ost03}, and the full center of $[v,v]$ is the full center $Z(\CV)$ \cite{Dav10}.
\end{rem}

\subsection{Examples}

\begin{expl}
Let $G$ be a finite group, viewed as a finite 2-group with only identity morphisms. Then $\vect_G$ is the usual category of finite-dimensional $G$-graded vector spaces. By \cite{Ost03}, the irreducible 2-representations of $G$, or equivalently, the indecomposable finite semisimple $\vect_G$-modules, are classified by pairs $(H,[\beta])$ up to conjugation, where $H$ is a subgroup of $G$ and $[\beta] \in H^2(H;\bk^\times)$.

Here we focus on the case $H = G$. For every $\beta \in Z^2(G;\bk^\times)$, there is a $\vect_G$-module structure on $\vect$ defined by the forgetful functor $\vect_G \to \vect$ equipped with the monoidal structure twisted by $\beta$. We denote this monoidal functor by $\forget^\beta \colon \vect_G \to \vect$ and the corresponding $\vect_G$-module by $\vect^\beta$. Two modules $\vect^\beta$ and $\vect^{\beta'}$ are equivalent if and only if $[\beta] = [\beta'] \in H^2(G;\bk^\times)$. It is not hard to see that $(\vect^\beta)^\op \simeq \vect^{\beta^{-1}}$ as 2-representations of $G$.

The 2-character $\chi_{\vect^\beta} \in \fun(G,\vect)$ is the constant functor on $\bk$:
\[
\chi_{\vect^\beta}(g) = \bk , \quad \forall g \in G .
\]
Then $\Psi(\chi_{\vect^\beta}) \in \vect_G$ is isomorphic to $\bk[G]$ as $G$-graded vector spaces. The multiplication, by \eqref{diag_multiplication_2-character}, is given by
\[
\chi_{\vect^\beta}(h) \otimes \chi_{\vect^\beta}(k) = \bk \otimes \bk \simeq \bk \xrightarrow{\beta(h,k)^{-1}} \bk = \chi_{\vect^\beta}(hk) .
\]
Thus $\Psi(\chi_{\vect^\beta})$ is the twisted group algebra $\bk[G,\beta^{-1}]$, which is a Lagrangian algebra in $\FZ_1(\vect_G)$ \cite{Dav10a}. This is exactly the full center of $\vect^{\beta^{-1}}$ in $\FZ_1(\vect_G)$.
\end{expl}

\begin{expl}
Let $A$ be a finite abelian group. Then $\mathrm B A$ is a finite 2-group and $\vect_{\mathrm B A} \simeq \vect^{\oplus \hat A}$ as multi-fusion categories \cite[Example 3.25]{HZ23}. The irreducible 2-representations of $\mathrm B A$ are the category $\vect$ equipped with the natural isomorphisms $\{\alpha^{\rho,a} \colon 1_\vect \Rightarrow 1_\vect\}_{a \in A}$ defined by
\[
\alpha^{\rho,a}_V \coloneqq \rho(a) \cdot 1_V ,
\]
where $\rho \in \hat A$. We denote this 2-representation by $\vect^\rho$. Clearly $(\vect^\rho)^\op \simeq \vect^\rho$.

The 2-character $\chi_{\vect^\rho} \in \fun(\mathrm B A,\vect) \simeq \rep(A)$ is given by the character $\rho$:
\[
\chi_{\vect^\rho}(\bullet) = \bk , \quad \chi_{\vect^\rho}(a) = \rho(a) \cdot 1_\bk , \quad \forall a \in A .
\]
It is the simple object of $\rho$-th component of $\vect_{\mathrm B A} \simeq \vect^{\oplus \hat A} \simeq \FZ_1(\vect_{\mathrm B A})$, which admits an obvious Lagrangian algebra structure.
\end{expl}

\begin{expl}
Let $\CG$ be the finite 2-group defined as follows:
\bit
\item Its first homotopy group is $\Zb_2$, and the second homotopy group is $\Zb_3$.
\item The action of $\Zb_2$ on $\Zb_3$ is the nontrivial one.
\item The 3-cocycle representing the associator is trivial.
\eit
We denote $\pi_1(\CG) = \{e,x\}$ with $x^2 = e$ and $\pi_2(\CG) = \{1,y,y^2\}$ with $y^3 = 1$. The 2-representations of $\CG$ were explicitly computed in \cite[Example 3.35]{HZ23}. As a multi-fusion category, $\vect_\CG$ is equivalent to $\vect_{\Zb_2} \oplus \End_\bk(\vect \oplus \vect)$. There are 3 irreducible 2-representations of $\CG$:
\bit
\item The first one is $\vect$ equipped with the trivial $\CG$-action, denoted by $\one$. This is the tensor unit of $2\rep(\CG)$.
\item The second one is $\vect_{\Zb_2}$ equipped with the following $\CG$-action. The object $x \in \CG$ acts by permuting two simple objects, and the morphism $y \in \pi_2(\CG)$ acts as identity. This 2-representation is denoted by $\one_c$.
\item The last one is $\vect \oplus \vect$ equipped with the following $\CG$-action. The object $x \in \CG$ acts by permuting two simple objects, and the morphism $y \in \pi_2(\CG)$ acts as the natural transformation whose components on two simple objects are $\omega$ and $\omega^2$. Here $\omega \in \bk^\times$ is a 3rd root of unity. This 2-representation is denoted by $S$.
\eit
The hom categories between these simple objects are depicted in the following graph:
\[
\xymatrix{
\one \ar@(ul,ur)[]^{\rep(\Zb_2)}  \ar@/^/[rr]^{\vect} & & \one_c \ar@(ul,ur)[]^{\vect_{\Zb_2}} \ar@/^/[ll]^{\vect} & S \ar@(ul,ur)[]^{\vect}
} .
\]

Now we compute the 2-characters of these irreducible 2-representations. First, as functors from $\CG$ to $\vect$,
\bit
\item $\chi_\one$ is the constant functor: $\chi_\one(e) = \chi_\one(x) = \bk$ equipped with the trivial $\pi_2(\CG) = \Zb_3$-action;
\item $\chi_{\one_c}(e) = \bk^2$ equipped with the trivial $\Zb_3$-action, and $\chi_{\one_c}(x) = 0$;
\item $\chi_S(e) = \bk^2$ equipped with the action of $y$ given by $\omega \oplus \omega^2$, and $\chi_S(x) = 0$.
\eit
By Proposition \ref{prop_2-character_direct_sum_tensor_product}, we can compute the fusion rules from these functors:
\[
\one_c \boxtimes \one_c \simeq \one_c \oplus \one_c , \quad S \boxtimes S \simeq \one_c \oplus S , \quad \one_c \boxtimes S \simeq S \boxtimes \one_c \simeq S \oplus S .
\]
These fusion rules were also obtained in \cite[Example 3.35]{HZ23}. As Lagrangian algebras in $\FZ_1(\vect_\CG) \simeq \FZ_1(\vect_{\Zb_2}) \oplus \vect$,
\bit
\item $\chi_\one$ is the group algebra $\bk[\Zb_2] \in \FZ_1(\vect_{\Zb_2})$;
\item $\chi_{\one_c}$ is the function algebra $\fun(\Zb_2) \in \FZ_1(\vect_{\Zb_2})$;
\item $\chi_S$ is the tensor unit in $\FZ_1(\End_\bk(\vect \oplus \vect)) \simeq \vect$.
\eit
\end{expl}

\begin{expl}
Let $\CG$ be the finite 2-group defined as follows:
\bit
\item Its first homotopy group is $\Zb_2$, and the second homotopy group is $\Zb_2$.
\item The action of $\Zb_2$ on $\Zb_2$ is trivial.
\item The 3-cocycle representing the associator is the nontrivial element $\alpha \in H^3(\Zb_2;\Zb_2) \simeq \Zb_2$.
\eit
We denote $\pi_1(\CG) = \{e,x\}$ with $x^2 = e$ and $\pi_2(\CG) = \{1,y\}$ with $y^2 = 1$. The 2-representations of $\CG$ were explicitly computed in \cite[Example 3.36]{HZ23}. As a multi-fusion category, $\vect_\CG$ is equivalent to $\vect_{\Zb_2} \oplus \vect_{\Zb_2}^\alpha$. There are 3 irreducible 2-representations of $\CG$:
\bit
\item The first one is $\vect$ equipped with the trivial $\CG$-action, denoted by $\one$. This is the tensor unit of $2\rep(\CG)$.
\item The second one is $\vect_{\Zb_2}$ equipped with the following $\CG$-action. The object $x \in \CG$ acts by permuting two simple objects, and the morphism $y \in \pi_2(\CG)$ acts as identity. This 2-representation is denoted by $\one_c$.
\item The last one is $\vect_{\Zb_2}^\alpha$ equipped with the following $\CG$-action. The object $x \in \CG$ acts by permuting two simple objects, and the morphism $y \in \pi_2(\CG)$ acts as $-1$. This 2-representation is denoted by $T$.
\eit
The hom categories between these simple objects are depicted in the following graph:
\[
\xymatrix{
\one \ar@(ul,ur)[]^{\rep(\Zb_2)}  \ar@/^/[rr]^{\vect} & & \one_c \ar@(ul,ur)[]^{\vect_{\Zb_2}} \ar@/^/[ll]^{\vect} & T \ar@(ul,ur)[]^{\vect_{\Zb_2}^\alpha}
} .
\]

Now we compute the 2-characters of these irreducible 2-representations. First, as functors from $\CG$ to $\vect$,
\bit
\item $\chi_\one$ is the constant functor: $\chi_\one(e) = \chi_\one(x) = \bk$ equipped with the trivial $\pi_2(\CG) = \Zb_2$-action;
\item $\chi_{\one_c}(e) = \bk^2$ equipped with the trivial $\Zb_2$-action, and $\chi_{\one_c}(x) = 0$;
\item $\chi_T(e) = \bk^2$ equipped with the action of $y$ given by $-1$, and $\chi_T(x) = 0$.
\eit
By Proposition \ref{prop_2-character_direct_sum_tensor_product}, we can compute the fusion rules from these functors:
\[
\one_c \boxtimes \one_c \simeq \one_c \oplus \one_c , \quad T \boxtimes T \simeq \one_c \oplus \one_c , \quad \one_c \boxtimes T \simeq T \boxtimes \one_c \simeq T \oplus T .
\]
These fusion rules were also obtained in \cite[Example 3.36]{HZ23}. As Lagrangian algebras in $\FZ_1(\vect_\CG) \simeq \FZ_1(\vect_{\Zb_2}) \oplus \FZ_1(\vect_{\Zb_2}^\alpha)$,
\bit
\item $\chi_\one$ is the group algebra $\bk[\Zb_2] \in \FZ_1(\vect_{\Zb_2})$;
\item $\chi_{\one_c}$ is the function algebra $\fun(\Zb_2) \in \FZ_1(\vect_{\Zb_2})$;
\item $\chi_T$ is also the function algebra $\fun(\Zb_2)$, which is the unique Lagrangian algebra in $\FZ_1(\vect_{\Zb_2}^\alpha)$.
\eit
\end{expl}

%% modular invariance and Lagrangian algebra?

\subsection{Joint 2-characters} \label{sec_numerical_2-character}

Let $\CG$ be a finite 2-group and $\CV \in 2\rep(\CG)$. Suppose $g,h \in \CG$ and $a \colon g \otimes h \to h \otimes g$ is a morphism in $\CG$. We define the \emph{joint 2-character} of $a$ to be
\[
\chi_\CV^{(2)}(a) \coloneqq \tr \bigl( \chi_\CV(h) \xrightarrow{\psi_{g,h}} \chi_\CV((g \otimes h) \otimes g^*) \xrightarrow{\chi_\CV(a \otimes 1)} \chi_\CV((h \otimes g) \otimes g^*) \simeq \chi_\CV(h) \bigr) \in \bk .
\]
From the TQFT point of view (see Remark \ref{rem_G_TQFT}), $\chi_\CV^{(2)}(a)$ is the partition function of $Z_\CV$ on the following torus:
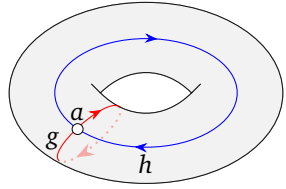
\[
\begin{tikzpicture}[scale=1.2]
\draw[even odd rule,fill=gray!10,opacity=0.7] (0,0) ellipse (1.5 and 1) (0.5,0) .. controls (0.2,0.3) and (-0.2,0.3) .. (-0.5,0) .. controls (-0.2,-0.3) and (0.2,-0.3) .. cycle ;
%
%\draw[thin,red] (-1,-0.8)--(-0.3,-0.2) ;
\draw[->-=0.7,thick,red,dotted] (-0.3,-0.15) .. controls (0,-0.35) and (-0.7,-0.95) .. (-1,-0.75) ;
\draw[even odd rule,fill=gray!10,opacity=0.7] (0,0) ellipse (1.5 and 1) (0.5,0) .. controls (0.2,0.3) and (-0.2,0.3) .. (-0.5,0) .. controls (-0.2,-0.3) and (0.2,-0.3) .. cycle ;
\draw (-0.6,0.1)--(-0.5,0) (0.5,0)--(0.6,0.1) ;
%
%\draw[->-,thin,blue] (0,0) ellipse (1 and 0.6) ;
\draw[->-,thin,blue] (1,0) arc (0:-180:1 and 0.6) ;
\draw[->-,thin,blue] (-1,0) arc (180:0:1 and 0.6) ;
\draw[->-=0.8,thin,red] (-1,-0.75) .. controls (-1.3,-0.55) and (-0.6,0.05) .. (-0.3,-0.15) ;
\draw[very thin,fill=white] (-0.86,-0.29) circle (0.06) node[above] {$\,a$} ;
\node at (-1.14,-0.45) {$g$} ;
\node at (0,-0.78) {$h$} ;
\end{tikzpicture}
\]
This torus can also be represented as a square whose edges are identified with their opposite edges:
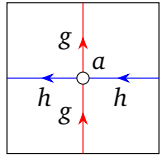
\[
\begin{tikzpicture}
\draw (0,0) rectangle (2,2) ;
\draw[->-,thin,red] (1,0)--(1,1) node[midway,left,black] {$g$} ;
\draw[->-,thin,red] (1,1)--(1,2) node[midway,left,black] {$g$} ;
\draw[->-,thin,blue] (2,1)--(1,1) node[midway,below,black] {$h$} ;
\draw[->-,thin,blue] (1,1)--(0,1) node[midway,below,black] {$h$} ;
\draw[very thin,fill=white] (1,1) circle (0.08) node[above right] {$a$} ;
\end{tikzpicture}
\]
Let $a' \colon h^L \otimes g \to g \otimes h^L$ be the morphism defined as follows (we omit the associators):
\[
a' \coloneqq \bigl( h^L \otimes g \xrightarrow{1 \otimes \coev_h} h^L \otimes g \otimes h \otimes h^L \xrightarrow{1 \otimes a \otimes 1} h^L \otimes h \otimes g \otimes h^L \xrightarrow{\ev_h \otimes 1} g \otimes h^L \bigr) .
\]
Then the invariance of the partition function under the modular $S$ transformation implies that $\chi_\CV^{(2)}(a) = \chi_\CV^{(2)}(a')$:
\[
Z_\CV \biggl(
\begin{array}{c}
\begin{tikzpicture}
\draw (0,0) rectangle (2,2) ;
\draw[->-,thin,red] (1,0)--(1,1) node[midway,left,black] {$g$} ;
\draw[->-,thin,red] (1,1)--(1,2) node[midway,left,black] {$g$} ;
\draw[->-,thin,blue] (2,1)--(1,1) node[midway,below,black] {$h$} ;
\draw[->-,thin,blue] (1,1)--(0,1) node[midway,below,black] {$h$} ;
\draw[very thin,fill=white] (1,1) circle (0.08) node[above right] {$a$} ;
\end{tikzpicture}
\end{array} \biggr) =
Z_\CV \biggl(
\begin{array}{c}
\begin{tikzpicture}
\draw (0,0) rectangle (2,2) ;
\draw[->-,thin,blue] (1,1)--(1,0) node[midway,left,black] {$h$} ;
\draw[->-,thin,blue] (1,2)--(1,1) node[midway,left,black] {$h$} ;
\draw[->-,thin,red] (2,1)--(1,1) node[midway,below,black] {$g$} ;
\draw[->-,thin,red] (1,1)--(0,1) node[midway,below,black] {$g$} ;
\draw[very thin,fill=white] (1,1) circle (0.08) node[above right] {$a$} ;
\end{tikzpicture}
\end{array} \biggr) =
Z_\CV \biggl(
\begin{array}{c}
\begin{tikzpicture}
\draw (0,0) rectangle (2,2) ;
\draw[->-,thin,blue] (1,0)--(1,1) node[midway,left,black] {$h^L$} ;
\draw[->-,thin,blue] (1,1)--(1,2) node[midway,left,black] {$h^L$} ;
\draw[->-,thin,red] (2,1)--(1,1) node[midway,below,black] {$g$} ;
\draw[->-,thin,red] (1,1)--(0,1) node[midway,below,black] {$g$} ;
\draw[very thin,fill=white] (1,1) circle (0.08) node[above right] {$a'$} ;
\end{tikzpicture}
\end{array} \biggr)
\]
In particular, when $h = e$ and $a \colon g \otimes e \to e \otimes g$ is the obvious isomorphism, we have
\[
\dim \chi_\CV(g) = \chi_{\dim(\CV)}(g) ,
\]
where the right hand side is the character of the $\pi_1(\CG)$-representation $\dim(\CV) = \chi_\CV(e)$ as defined in \eqref{eq_pi_1_action_dim}. If we do the modular $S$ transformation twice, we obtain
\[
\chi_\CV^{(2)}(a) = \chi_\CV^{(2)}(a') = \chi_\CV^{(2)}(a'') = \chi_\CV^{(2)}(a^L \colon g^L \otimes h^L \to h^L \otimes g^L) .
\]
Similarly, the invariance of the partition function under the modular $T$ transformation implies that
\[
\chi_\CV^{(2)}(a) = Z_\CV \biggl(
\begin{array}{c}
\begin{tikzpicture}
\draw (0,0) rectangle (2,2) ;
\draw[->-,thin,red] (1,0)--(1,1) node[midway,left,black] {$g$} ;
\draw[->-,thin,red] (1,1)--(1,2) node[midway,left,black] {$g$} ;
\draw[->-,thin,blue] (2,1)--(1,1) node[midway,below,black] {$h$} ;
\draw[->-,thin,blue] (1,1)--(0,1) node[midway,below,black] {$h$} ;
\draw[very thin,fill=white] (1,1) circle (0.08) node[above right] {$a$} ;
\end{tikzpicture}
\end{array} \biggr) = Z_\CV \biggl(
\begin{array}{c}
\begin{tikzpicture}
\draw (0,0) rectangle (2,2) ;
\draw[->-=0.3,thin,red] (1,0) .. controls (1,1) and (1,0.8) .. (0,0.8) node[midway,below left,black] {$g$} ;
\draw[->-,thin,red] (2,0.8) .. controls (1,0.8) and (1,0.7) .. (1,1) node[midway,below right,black] {$g$} ;
\draw[->-,thin,red] (1,1)--(1,2) node[midway,left,black] {$g$} ;
\draw[->-,thin,blue] (2,1)--(1,1) node[midway,above,black] {$h$} ;
\draw[->-,thin,blue] (1,1)--(0,1) node[midway,above,black] {$h$} ;
\draw[very thin,fill=white] (1,1) circle (0.08) node[above right] {$a$} ;
\end{tikzpicture}
\end{array} \biggr) = \chi_\CV^{(2)} \bigl( g \otimes (g \otimes h) \xrightarrow{1 \otimes a} g \otimes (h \otimes g) \simeq (g \otimes h) \otimes g \bigr) .
\]
In particular, when $h = e$ and $a \colon g \otimes e \to e \otimes g$ is the obvious isomorphism, the right hand side is
\[
\chi_\CV^{(2)}(1_{g \otimes g}) = \tr\bigl( \chi_\CV(g) \xrightarrow{\psi_{g,g}} \chi_\CV((g \otimes g) \otimes g^*) \simeq \chi_\CV(g) \bigr) ,
\]
which is $\dim(\chi_\CV(g))$ by the $S^1$-invariance of $\chi_\CV$. This also shows that $\dim(\chi_\CV(g)) = \dim(\chi_\CV(g \otimes g))$. Moreover, the joint 2-characters also have the conjugation invariance: for every $k \in \CG$ we have
\[
\chi_\CV^{(2)}(a) = \chi_\CV^{(2)} \bigl(1_k \otimes a \otimes 1_k^* \colon (k \otimes g \otimes k^*) \otimes (k \otimes h \otimes k^*) \to (k \otimes h \otimes k^*) \otimes (k \otimes g \otimes k^*) \bigr) .
\]

%% There is a space with a modular group action!

\begin{rem}
Let $\SC$ be a 2-category. Suppose $A,B \colon x \to x$ are 1-endomorphisms in $\SC$ and $\eta \colon A \circ B \Rightarrow B \circ A$ is a 2-isomorphism. Ganter and Kapranov defined the ``joint trace'' for $(A,B,\eta)$ in \cite{GK08}. When $\SC = 2\vect$, the joint trace of $(g \odot -,h \odot -,a \odot -)$ is the joint 2-character $\chi_\CV^{(2)}(a)$ defined as above.
\end{rem}

\begin{rem}
Suppose $\CG = G$ is an ordinary finite group. Then the morphism $g \otimes h \to h \otimes g$ must be identity, and its existence implies the $g$ and $h$ commute. So the joint 2-character $\chi_\CV^{(2)}$ can be viewed as functions on the pairs of commuting elements in the group $G$. This case was studied by Ganter and Kapranov \cite{GK08} (see also \cite{Oso10,RW18} for explicit calculations).

On the other hand, for an object $A \in \FZ_1(\vect_G)$ (i.e., $A$ is a finite-dimensional $G$-graded vector space equipped with a compatible $G$-action), there is a function $\chi_A$, also called the ``character'' of $A$, defined by
\[
\chi_A(g,h) \coloneqq \chi_{A_g}(h) , \quad g,h \in G \text{ commute} ,
\]
and it has a modular invariance when $A$ is a Lagrangian algebra (see for example \cite{Dav10a} and references therein). By Theorem \ref{thm_2-character_Lagrangian_algebra}, the ``character'' $\chi_{Z(\CV^\op)}$ of the Lagrangian algebra $Z(\CV^\op)$ coincides with the joint 2-character $\chi_\CV^{(2)}$. Then the modular invariance of the joint 2-character also gives a geometric understanding of the modular invariance of $\chi_A$ for Lagrangian algebras $A \in \FZ_1(\vect_G)$.
\end{rem}

\begin{rem}
The joint 2-characters can be understood heuristically in the context of the generalized character theory of Hopkins, Kuhn, and Ravenel \cite{HKR00}. For a finite group $G$ and height-$n$ \emph{Morava $E$-theory} $E_n$, the HKR character map \cite[Theorem C]{HKR00} provides an isomorphism
\[
C_0 \otimes_{E_n^0} E_n^0(\mathrm B G) \simeq \mathrm{Cl}_n(G, C_0),
\]
where the right-hand side is the ring of $C_0$-valued generalized class functions on $\Hom((\mathbb{Z}^\wedge_p)^n, G)/G$, the set of conjugacy classes of $n$-tuples of pairwise commuting $p$-power-order elements of $G$, where $\mathbb{Z}^\wedge_p := \lim_{k \to \infty} \mathbb{Z}/p^k\mathbb{Z}$ denotes the $p$-adic completion of $\mathbb{Z}$, and $C_0$ is a certain faithfully flat extension of $E_n^0$.

At height 1 this recovers ordinary character theory; at height 2 the characters are class functions on commuting pairs.

As observed by Ganter and Kapranov \cite{GK08}, the 2-character of a categorical representation of a finite group $G$, evaluated on commuting pairs, gives precisely a height-2 generalized class function: for $G$ an ordinary group, the joint 2-character $\chi_V^{(2)}$ is a function on $\Hom(\mathbb{Z}^2, G)/G$, which becomes the domain of HKR characters at height 2 after $p$-adic completion. The induced 2-character formula of \cite[Theorem 4.1]{GK08} reflects the HKR transfer map. This was the original motivation for the categorification program of \cite{GK08}.

For a genuine 2-group $\CG$ whose classifying space $\lvert \mathrm{B}\CG \rvert$ is a homotopy 2-type, the correct generalization is to replace $\Hom(\mathbb{Z}^2, \CG)/\CG$ by the mapping 2-groupoid $\map(T^2,\lvert \mathrm{B}\CG \rvert)$, which parametrizes maps from the torus to $|\mathrm{B}\CG|$. A map $T^2 \to \lvert \mathrm{B}\CG \rvert$ is determined (up to homotopy) by a pair of objects $g, h \in \CG$ together with a morphism $a \colon g \otimes h \to h \otimes g$, modulo conjugation — which is precisely the data on which our joint 2-character $\chi_V^{(2)}(a)$ is evaluated in this section. The modular invariance of $\chi_V^{(2)}$ then reflects the action of the mapping class group $SL_2(\mathbb{Z})$ on $\map(T^2,\lvert \mathrm{B}\CG \rvert)$.

From this perspective, the passage from ordinary characters to 2-characters categorifies the passage from height 1 to height 2 in HKR theory: the 2-character $\chi_V \colon G \to \vect$ is a categorified class function (a functor on the loop groupoid $\map(S^1,\lvert \mathrm{B}G \rvert)$), while the joint 2-character $\chi_V^{(2)}$ is a numerical invariant on $\map(T^2,\lvert \mathrm{B}G \rvert)$. The $K(2)$-Euler characteristic formula $\chi_{K(2)}(BG) = \lvert \Hom((\mathbb{Z}^\wedge_p)^2, G)/G \rvert$ of \cite[Theorem B]{HKR00}  counts exactly the conjugacy classes of commuting pairs, which is the number of isomorphism classes in the domain of $\chi_V^{(2)}$ when $G$ is a finite group.
\end{rem}

\subsection{Orthogonality of 2-characters} \label{sec_orthogonality}

Let $\CB$ be a nondegenerate braided fusion category and $A,B \in \CB$ be Lagrangian algebras. Then $\RMod_A(\CB)$ and $\RMod_B(\CB)$ are fusion categories whose Drinfeld centers are equivalent to $\CB$ as braided fusion categories (the equivalences are given by the $\alpha$-induction). By \cite[Proposition 4.8]{DMNO13}, the Lagrangian algebra $B \in \CB \simeq \FZ_1(\RMod_A(\CB))$ determines an indecomposable finite semisimple $\RMod_A(\CB)$-module, denoted by $\CM$. By \cite{ENO11,ENO10}, this category $\CM$ is an invertible $\RMod_A(\CB)$-$\RMod_B(\CB)$-bimodule and can be described as follows. First, we view $B$ as a Lagrangian algebra in $\FZ_1(\RMod_A(\CB))$. Then the image of $B$ under the forgetful functor $\forget \colon \FZ_1(\RMod_A(\CB)) \to \RMod_A(\CB)$ is a separable algebra in $\RMod_A(\CB)$. Note that we have the following commutative diagram:
\[
\xymatrix{
\CB \ar[rr]^-{\simeq} \ar[dr]_{- \otimes A} & & \FZ_1(\RMod_A(\CB)) \ar[dl]^{\forget} \\
 & \RMod_A(\CB)
}
\]
So this algebra in $\RMod_A(\CB)$ is $B \otimes A$. In general, this is not an indecomposable algebra. If we decompose it as the direct sum of indecomposable subalgebras $B \otimes A \simeq \bigoplus_{i=1}^n B_i$ in $\RMod_A(\CB)$, then all direct summands $B_i$ are Morita equivalent, and we have $\RMod_{B_i}(\RMod_A(\CB)) \simeq \CM$.

On the other hand, the image of the Lagrangian algebra $B$ in $\FZ_1(\RMod_A(\CB))$ is the full center of $\CM$. Thus we have
\[
B \otimes A = \forget(Z(\CM)) = \int_{m \in \CM} [m,m]_{\RMod_A(\CB)} \simeq \bigoplus_{m \in \Irr(\CM)} [m,m]_{\RMod_A(\CB)} .
\]
Therefore, the indecomposable subalgebras $B_i$ are the internal homs $[m,m] \in \RMod_A(\CB)$, and the number of direct summands is equal to the number of simple objects in $\CM$. This equality of two numbers can also be categorified:
\begin{multline} \label{eq_intersection_Lagrangian_algebra}
\CB(\one,B \otimes A) \simeq \RMod_A(\CB)(A,B \otimes A) \simeq \int_{m \in \CM} \RMod_A(\CB)(A,[m,m]) \\
\simeq \int_{m \in \CM} \CM(A \odot m,m) \simeq \int_{m \in \CM} \CM(m,m) = \dim(\CM) .
\end{multline}
Note that $\CB(\one,-) \colon \CB \to \vect$ is a lax monoidal functor, thus preserves algebras. It is not hard to see that \eqref{eq_intersection_Lagrangian_algebra} is an isomorphism of $\bk$-algebras.

\begin{rem}
The nondegenerate braided fusion category $\CB$ can be viewed as the category of particle-like topological defects in a 2d (spatial dimension) topological order. The Lagrangian algebras $A,B \in \CB$ lead to two anyon condensation processes \cite{Kon14e} and produce two 1d boundary topological orders. The fusion category of particle-like topological defects on these two boundaries are $\RMod_A(\CB)$ and $\RMod_B(\CB)$, respectively. The gapped 0d domain walls between these two boundaries also form a category, which is exactly the invertible bimodule $\CM$ given above. We depict these topological orders in the following figure:
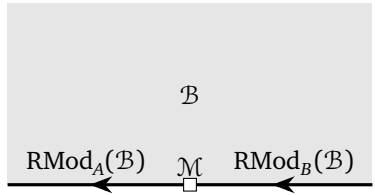
\[
\begin{tikzpicture}[scale=1.2]
\fill[gray!20] (-2,0) rectangle (2,2) node[midway,black] {$\CB$} ;
\draw[very thick,->-] (0,0)--(-2,0) node[midway,above] {$\RMod_A(\CB) \quad$} ;
\draw[very thick,->-] (2,0)--(0,0) node[midway,above] {$\quad \RMod_B(\CB)$} ;
\draw[fill=white] (-0.07,-0.07) rectangle (0.07,0.07) node[midway,above] {$\CM$} ;
\end{tikzpicture}
\]
For more details and intuitions, we refer readers to \cite[Section 4.4]{KZ22a}.

Then we provide a physical motivation and intuition of the isomorphism \eqref{eq_intersection_Lagrangian_algebra}. Denote $\CC \coloneqq \RMod_A(\CB)$, $\CD \coloneqq \RMod_B(\CB)$, $F \coloneqq - \otimes A \colon \CB \to \RMod_A(\CB)$ and $G \coloneqq - \otimes B \colon \CB \to \RMod_B(\CB)$. Let us compute the ground state degeneracy of the following system on a cylinder (which is the spatial manifold):
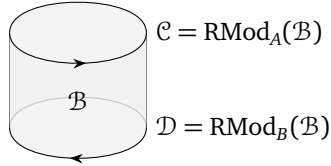
\[
\begin{tikzpicture}[scale=1.8]
%% cylinder
\draw (0,0)--(0,0.7) ;
\draw (1,0)--(1,0.7) ;
\fill[gray!10,opacity=0.8] (0,0)--(0,0.7) .. controls (0,1) and (1,1) .. (1,0.7)--(1,0) .. controls (1,0.3) and (0,0.3) .. cycle ;
\draw (0,0.7) .. controls (0,1) and (1,1) .. (1,0.7) node[right] {$\CC = \RMod_A(\CB)$} ;
\draw (0,0) .. controls (0,0.3) and (1,0.3) .. (1,0) node[right] {$\CD = \RMod_B(\CB)$};
\fill[gray!10,opacity=0.8] (0,0)--(0,0.7) .. controls (0,0.4) and (1,0.4) .. (1,0.7)--(1,0) .. controls (1,-0.3) and (0,-0.3) .. cycle ;
\draw[->-] (0,0.7) .. controls (0,0.4) and (1,0.4) .. (1,0.7) ;
\draw[->-] (1,0) .. controls (1,-0.3) and (0,-0.3) .. (0,0) ;
\node at (0.5,0.2) {$\CB$} ;
\end{tikzpicture}
\]
The 2d topological order $\CB$ is put on the 2d cylinder, and two boundary topological orders $\CC,\CD$ are put on the top and bottom boundaries of the cylinder, respectively. For mathematicians, the ground state degeneracy can be understood as the factorization homology.
\bit
\item If we shrink two circles to points, we obtain a sphere labeled by $\CB$ with two punctures labeled by $A$ and $B$, respectively (see \cite[Corollary 5.15]{AKZ17}). The ground state degeneracy of the system on this sphere is $\CB(\one_\CB,A \otimes B)$ (this is a special case of \cite[Corollary 5.18]{AKZ17}).
\item If we shrink this cylinder along the vertical direction, we obtain a 1d topological order on a circle. This 1d topological order is not stable in general, and the particle-like topological defects form the multi-fusion category $\fun_\bk(\CM,\CM)$ (this is a corollary of the boundary-bulk relation \cite{KWZ15,KWZ17}; see \cite[Section 4.4]{KZ22a} for example). Then its ground state degeneracy on a circle is the endomorphism space of its tensor unit, that is, $\dim(\CM)$.
\eit
These two different ways to computing the ground state degeneracy lead to the isomorphism \eqref{eq_intersection_Lagrangian_algebra}.

Moreover, we also have the following isomorphism:
\begin{multline*}
\int^{x \in \CB} \CC(F(x),\one_\CC) \otimes \CD(\one_\CD,G(x)) \simeq \int^{x \in \CB} \CC(F(x),\one_\CC) \otimes \CD(G(x^R),\one_\CD) \\
\simeq \int^{x \in \CB} \CB(x,F^R(\one_\CC)) \otimes \CB(x^R,G^R(\one_\CD)) \simeq \int^{x \in \CB} \CB(x,A) \otimes \CB(x^R,B) \\
\simeq \CB(A^R,B) \simeq \CB(\one,B \otimes A) \simeq \dim(\CM) .
\end{multline*}
By Corollary \ref{cor_semisimple_end_direct_sum}, the left hand side is isomorphic to
\[
\bigoplus_{x \in \Irr(\CB)} \CC(F(x),\one_\CC) \otimes \CD(\one_\CD,G(x)) .
\]
By counting the dimension, it follows that the number of simple objects in $\CM$ is equal to
\[
\sum_{x \in \Irr(\CB)} [F(x):\one_\CC] \cdot [G(x):\one_\CD] ,
\]
where $[a:b]$ is the multiplicity of the simple object $b$ in $a$. This equality is well-known for physicists (see for example \cite{LWW15,SH19}). 
\end{rem}

When $\CB$ is a nondegenerate braided multi-fusion category, it is possible that two Lagrangian algebras $A$ and $B$ are contained in different indecomposable direct summands of $\CB$. In this case we understand that $\CM = 0$ and \eqref{eq_intersection_Lagrangian_algebra} still holds. This convention is convenient when we consider from the `boundary point of view' rather than `bulk point of view'.

\begin{prop} \label{prop_Lagrangian_algebra_intersection}
Let $\CC$ be a multi-fusion category. Suppose $\CV,\CW$ are indecomposable finite semisimple left $\CC$-modules. Then there is an isomorphism of $\bk$-algebras
\[
\FZ_1(\CC)(\one,Z(\CW) \otimes Z(\CV)) \simeq \dim(\fun_\CC(\CV,\CW)) .
\]
\end{prop}

\pf
Take $\CB \coloneqq \FZ_1(\CC)$, $A \coloneqq Z(\CV)$ and $B \coloneqq Z(\CW)$. Then $\RMod_A(\CB) \simeq \fun_\CC(\CV,\CV)^\rev$ and $\RMod_B(\CB) \simeq \fun_\CC(\CW,\CW)^\rev$. When $\CV,\CW \in \RMod_\CC(2\vect)$ are connected (i.e., $\fun_\CC(\CV,\CW) \neq 0$), these two fusion categories are Morita equivalent, and $\fun_\CC(\CV,\CW)$ and $\fun_\CC(\CW,\CV)$ are invertible bimodules. Then use \eqref{eq_intersection_Lagrangian_algebra} we obtain the result. When $\CV,\CW \in \RMod_\CC(2\vect)$ are not connected, both sides are zero.
\epf

\medskip
Let $\CG$ be a finite 2-group. Suppose $(V,u = \{u_g \colon (g \otimes V) \otimes g^* \to V\}_{g \in G}) \in (\vect_\CG)^\CG \simeq \FZ_1(\vect_\CG)$ and $(F,\phi = \{\phi_g \colon g \odot F \Rightarrow F\}_{g \in \CG}) \in \fun(\CG,\vect)^\CG)$. Then for every $g \in \CG$ there is an automorphism
\begin{equation} \label{eq_pi_1_action_equivariant_evaluation}
F(V) \xrightarrow{(\phi_g)_V} F((g \otimes V) \otimes g^*) \xrightarrow{F(u_g)} F(V) .
\end{equation}
It is not hard to see that this defines a $\pi_1(\CG)$-action on the vector space $F(V)$. When $(V,u) \in \FZ_1(\vect_\CG)$ is the tensor unit, this $\pi_1(\CG)$-action is essentially defined in \eqref{eq_pi_1_action_dim}. We denote the $\pi_1(\CG)$-invariant subspace by $F[V] \coloneqq F(V)^{\pi_1(\CG)}$.

\begin{lem} \label{lem_invariant_pairing_hom}
Suppose $F = (F,\phi = \{\phi_g \colon g \odot F \Rightarrow F\}_{g \in \CG}) \in \fun(\CG,\vect)^\CG)$. Then $F[e]$ is isomorphic to the hom space
\[
\fun(\CG,\vect)^\CG(\Phi(e),F) .
\]
\end{lem}

\pf
By definition, $\fun(\CG,\vect)^\CG(\Phi(e),F)$ is a subspace of
\[
\fun(\CG,\vect)(\Phi(e),F) = \fun(\CG,\vect)(\vect_\CG(I(-),e),F) \simeq F(e) ,
\]
This isomorphism is given by the Yoneda lemma: it sends a natural transformation $\alpha \colon \Phi(e) \Rightarrow F$ to $\alpha_e(1_e) \in F(e)$, and its inverse sends an element $s \in F(e)$ to the natural transformation $\alpha^s \colon \Phi(e) \Rightarrow F$ defined by
\begin{align*}
\alpha^s_g \colon \Phi(e)(g) = \vect_\CG(I(g),e) & \to F(g) \\
f & \mapsto F(I(f))(s) = F(f)^{-1}(s) .
\end{align*}
By definition, a natural transformation $\alpha \colon \Phi(e) \Rightarrow F$ is a morphism in $\fun(\CG,\vect)^\CG$ if and only if the following diagram commutes for every $g,h \in \CG$:
\[
\xymatrix@C=5em{
\Phi(e)(h) \ar[r]^-{\alpha_h} \ar[d]_{\simeq} & F(h) \ar[d]^{\phi_{g,h}} \\
\Phi(e)((g \otimes h) \otimes g^*) \ar[r]^-{\alpha_{(g \otimes h) \otimes g^*}} & F((g \otimes h) \otimes g^*)
}
\]
Now we take $h = e$ and $\alpha = \alpha^s$ for some $s \in F(e)$. It is clear from the above commutative diagram that $\alpha^s$ is a morphism in $\fun(\CG,\vect)^\CG$ if and only if $s = \alpha^s_e(1_e) \in F(e)^{\pi_1(\CG)}$ is contained in the $\pi_1(\CG)$-invariant subspace of $F(e)$.
\epf

\begin{cor} \label{cor_algebra_2-character_pairing}
Suppose $F \in \fun(\CG,\vect)^\CG$ is an algebra (i.e., equipped with a lax monoidal functor structure). Then $F[e] \in \vect$ is naturally a $\bk$-algebra.
\end{cor}

\pf
By Lemma \ref{lem_invariant_pairing_hom}, $F[e] \simeq \fun(\CG,\vect)^\CG(\Phi(e),F)$. Since $\Phi(e) \in \fun(\CG,\vect)^\CG$ is the tensor unit, the hom functor $\fun(\CG,\vect)^\CG(\Phi(e),-)$ is a lax monoidal functor and hence preserves algebras.
\epf

\begin{defn}
For $F,G \in \fun(\CG,\vect)^\CG$, their \emph{inner product} is defined to be the invariant space
\[
\langle F,G \rangle \coloneqq (F \circledast G)[e] = (F \circledast G)(e)^{\pi_1(\CG)} ,
\]
where the $\pi_1(\CG)$-action on $(F \circledast G)(e)$ is defined by \eqref{eq_pi_1_action_equivariant_evaluation}.
\end{defn}

\begin{thm} \label{thm_2-character_pairing_dim_hom}
For $\CV,\CW \in 2\rep(\CG)$, there is a canonical isomorphism of $\bk$-algebras
\[
\langle \chi_\CV,\chi_\CW \rangle \simeq \dim(\fun_\CG(\CV,\CW)) ,
\]
where the $\bk$-algebra structure of the left hand side is given by Corollary \ref{cor_algebra_2-character_pairing}.
\end{thm}

\pf
The vector spaces on both sides are bi-additive with respect to $\CV$ and $\CW$. So it suffices to prove this isomorphism when $\CV,\CW \in 2\rep(\CG)$ are simple.

By Lemma \ref{lem_invariant_pairing_hom} and Theorem \ref{thm_2-character_Lagrangian_algebra}, the left hand side is
\[
\langle \chi_\CV,\chi_\CW \rangle \simeq \fun(\CG,\vect)^\CG(\Phi(e),\chi_\CV \circledast \chi_\CW) \simeq \fun(\CG,\vect)^\CG(\Phi(e),\Phi(Z(\CV^\op)) \circledast \Phi(Z(\CW^\op))) .
\]
Since the Fourier 2-transform $\Phi \colon \FZ_1(\vect_\CG) \to \fun(\CG,\vect)^\CG$ is an equivalence, this is isomorphic to
\[
\FZ_1(\vect_\CG)(e,Z(\CV^\op) \otimes Z(\CW^\op)) .
\]
By Proposition \ref{prop_Lagrangian_algebra_intersection}, this is isomorphic to
\[
\dim(\fun_{\vect_\CG}(\CW^\op,\CV^\op)) \simeq \dim(\fun_{\vect_\CG}(\CV,\CW)) \simeq \dim(\fun_\CG(\CV,\CW)) ,
\]
where the first isomorphism is induced by taking right adjoints.
\epf

\begin{rem}
Theorem \ref{thm_2-character_pairing_dim_hom} is a categorification of the classical orthogonality of characters. For a finite group $G$ and $V,W \in \rep(G)$, we have the following orthogonality of characters:
\[
\langle \chi_V,\chi_W \rangle \coloneqq \frac{1}{\lvert G \rvert} (\chi_V \ast \chi_W)(e) = \frac{1}{\lvert G \rvert} \sum_{g \in G} \chi_V(g^{-1}) \chi_W(g) = \dim_\bk \Hom_G(V,W) .
\]
So we call Theorem \ref{thm_2-character_pairing_dim_hom} the \emph{orthogonality of 2-characters}.
\end{rem}

\begin{expl}
Let $\vect$ be the trivial 2-representation of $\CG$. Then for any $\CV \in 2\rep(\CG)$, by the definition of the inner product of class functors, we have
\[
\langle \chi_\vect,\chi_\CV \rangle \simeq \biggl( \int^{g \in \CG} \chi_\vect(g^L) \otimes \chi_\CV(g) \biggr)^{\pi_1(\CG)} \simeq \bigl( \bigoplus_{g \in \pi_1(\CG)} \bk \otimes \chi_\CV(g) \bigr)^{\pi_1(\CG)} \simeq \bigl( \bigoplus_{g \in \pi_1(\CG)} \chi_\CV(g) \bigr)^{\pi_1(\CG)} .
\]
By Theorem \ref{thm_2-character_pairing_dim_hom}, the left hand side is isomorphic to
\[
\dim(\fun_\CG(\vect,\CV)) \simeq \dim(\fun(\vect,\CV)^\CG) \simeq \dim(\CV^\CG) .
\]
So we obtain an isomorphism of $\bk$-algebras
\[
\dim(\CV^\CG) \simeq \bigl( \bigoplus_{g \in \pi_1(\CG)} \chi_\CV(g) \bigr)^{\pi_1(\CG)} .
\]
When $\CG$ is a 1-group, this is a special case of \cite[Theorem 5.17]{Gan15}. When $\CV = \vect$ is trivial, we have
\[
\dim(\vect^\CG) \simeq \dim(\rep(\pi_1(\CG))) = \nat(1_{\rep(\pi_1(\CG))},1_{\rep(\pi_1(\CG))}) ,
\]
and
\[
\langle \chi_\vect,\chi_\vect \rangle \simeq \bk[\pi_1(\CG)]^{\pi_1(\CG)} ,
\]
where the $\pi_1(\CG)$-action on $\bk[\pi_1(\CG)]$ is the conjugation action. It is well-known that both two algebras are isomorphic to the center $Z(\bk[\pi_1(\CG)])$.
\end{expl}

\begin{rem}
In \cite[Definition 5.21]{Gan15}, Ganter defined an inner product of two 2-representations of a 1-group $G$ in a lax monoidal linear 2-category. In the case we are concerned with, Ganter's inner product of two finite semisimple 2-representations $\CV,\CW$ of $G$ is the dimension space of the equivariantization $(\CV \boxtimes \CW)^G$. By Theorem \ref{thm_2-character_pairing_dim_hom}, our inner product of 2-characters
\[
\langle \chi_\CV,\chi_\CW \rangle \simeq \dim(\fun_G(\CV,\CW)) \simeq \dim(\fun_\bk(\CV,\CW)^G) \simeq \dim((\CW \boxtimes \CV^\op)^G)
\]
is isomorphic to the Ganter's inner product of $\CV^\op$ and $\CW$.
\end{rem}

\bibliography{Top}

\end{document}